\documentclass[12pt]{article}
\usepackage{amsmath}
\usepackage{mathdots}
\usepackage{amsfonts}
\usepackage{amssymb}
\usepackage{latexsym}
\usepackage{graphicx}
\usepackage{pgfplots}
\usepackage{todonotes}

\usepackage[ruled,lined,linesnumbered]{algorithm2e}

\usepackage{color}
\usepackage{fullpage}
\usepackage{verbatim}

\usepackage{color}
\usepackage{verbatim}

\def\P{\mathbb{P}}
\def\E{\mathbb{E}}
\def\EE{\mathbb{E}_{\frac{1}{w}}}
\def\V{{\bf V}}\def\I{{\bf I}}
\def\Q{{\bf Q}}\def\R{{\bf R}}
\def\Z{{\bf Z}}\def\X{{\bf X}}
\def\c{{\bf c}}\def\f{{\bf f}}\def\W{{\bf W}}\def\x{{\bf x}}\def\g{{\bf g}}\def\v{{\bf v}}\def\e{{\bf e}}\def\b{{\bf b}}

\DeclareMathOperator{\argmin}{argmin}

\newcommand{\yntodo}{\todo[color=green]}

\newtheorem{theorem}{Theorem}[section]
\newtheorem{lemma}{Lemma}[section]
\newtheorem{proposition}{Proposition}[section]

\numberwithin{equation}{section}
\numberwithin{theorem}{section} 
\numberwithin{proposition}{section} 
\numberwithin{lemma}{section} 
\numberwithin{table}{section} 
\numberwithin{figure}{section}

\newcommand{\rr}[1]{\textcolor{red}{#1}}

\newcommand{\ignore}[1]{}

\newcommand{\bb}[1]{\textcolor{blue}{#1}}

\title{
Approximate and integrate: 
Variance reduction in Monte Carlo integration via function approximation
}

\author{Yuji Nakatsukasa
\thanks{National Institute of Informatics, 2-1-2 Hitotsubashi, Chiyoda-ku, Tokyo 101-8430, Japan. Email: {\tt nakatsukasa@nii.ac.jp}}
}

\begin{document}
\maketitle

\begin{abstract}
Classical algorithms in numerical analysis for numerical integration (quadrature/cubature) follow the principle of approximate and integrate: the integrand is approximated by a simple function (e.g. a polynomial), which is then integrated exactly. In high-dimensional integration, such methods quickly become infeasible due to the curse of dimensionality. A common alternative is the Monte Carlo method (MC), which simply takes the average of random samples, improving the estimate as more and more samples are taken. The main issue with MC is its slow (though dimension-independent) convergence, and various techniques have been proposed to reduce the variance. In this work we suggest a numerical analyst's interpretation of MC: it approximates the integrand with a constant function, and integrates that constant exactly. This observation leads naturally to MC-like methods where the approximant is a non-constant function, for example low-degree polynomials, sparse grids or low-rank functions. We show that these methods have the same $O(1/\sqrt{N})$ asymptotic convergence as in MC, but with reduced variance,  equal to the quality of the underlying function approximation. We also discuss methods that improve the approximation quality as more samples are taken, and thus can converge faster than $O(1/\sqrt{N})$. The main message is that techniques in high-dimensional approximation theory can be combined with Monte Carlo integration to accelerate convergence. 

\end{abstract}

\section{Introduction}
This paper deals with the numerical evaluation (approximation) of the definite integral 
\begin{equation}
  \label{eq:goal}
I:=  \int_{\Omega} f(\x)d\x,
\end{equation}
for $f:\Omega\rightarrow \mathbb{R}$. 
For simplicity, we assume 
$\Omega$ is the $d$-dimensional cube 
$\Omega=[0,1]^d$ unless otherwise specified, although little of what follows relies crucially on this assumption. 
Our goal is to deal with the (moderately) high-dimensional case $d\gg 1$. 
The need to approximately evaluate integrals of the form~\eqref{eq:goal} arises in a number of applications, which are too numerous to list fully, but prominent examples include finance~\cite{glasserman2013monte}, machine learning~\cite{murphy2012machine}, biology~\cite{manly2006randomization}, and stochastic differential equations~\cite{Kloedenbook}. In many of these applications, the integral~\eqref{eq:goal} often represents the expected value of a certain quantity of interest. 

Integration is a classical subject in numerical analysis (e.g.,~\cite{davis2007methods,suli2003introduction,trefethenatap}). 
When $d=1$, effective integration (quadrature) rules  are available that converge extremely fast: for example, Gauss and Clenshaw-Curtis quadrature rules converge exponentially if $f$ is analytic. These formulas extend to $d>1$ by taking tensor products, but the computational complexity grows like $N=O(n^d)$ where $n$ is the number of sample points in one direction (equal to the maximum degree of the polynomial approximation underlying the integration; we give more details in Section~\ref{sec:mcasquad}). 


Among the alternatives for approximating~\eqref{eq:goal} when $d$ is large, Monte Carlo (MC) integration is one of the most widely used (another is sparse grids, which we treat briefly in Section~\ref{sec:spgridmc}). 
 In its simplest form, Monte Carlo integration~\cite{caflisch1998monte,robert2004monte} approximates the integral $\int_{[0,1]^d} f(\x)d\x$ by the average of $N(\gg 1)$ random samples 
 \begin{equation}   \label{eq:mcdef}
\int_{[0,1]^d} f(\x)d\x \approx \frac{1}{N}\sum_{i=1}^N f(\x_i)=:c_0.
 \end{equation}
 Here $\{\x_i\}_{i=}^N$ are sample points of $f$, chosen uniformly  at random in $[0,1]^d$. 
When $\Omega\neq [0,1]^d$, the MC estimate becomes $c_0|\Omega|$, where $|\Omega|$ is the volume of $\Omega$.  
MC is usually derived and analyzed using probability theory and statistics. In particular, 
 the central limit theorem
shows that for sufficiently large $N$ the MC error scales like
$\frac{\sigma(f)}{\sqrt{N}}=
\frac{1}{\sqrt{N}}\|f-\bar f\|_2$~\cite{caflisch1998monte}, where $(\sigma(f))^2$ is the variance of $f$. 
We prefer to rewrite this as (the reason will be explained in Section~\ref{sec:careful}) 
\begin{equation}
  \label{eq:MCvar}
\frac{1}{\sqrt{N}}\min_c\|f-c\|_2  .
\end{equation}
The estimate comes with a confidence interval: for example $c_0\pm 2
\min_c\|f-c\|_2/\sqrt{N} $ (or $2$ replaced with $1.96$) gives a $95\%$ confidence interval. In practice, the variance is estimated by 
\begin{equation}  \label{eq:sigmaestimate}
\widetilde\sigma^2 = \frac{1}{N-1}\sum_{i=1}^N(f(\x_i)-c_0)^2. 
\end{equation}


A significant advantage of MC is that its convergence~\eqref{eq:MCvar}
is independent of the dimension $d$. The disadvantage is the `slow'
$\sigma(f)/\sqrt{N}$ convergence. Great efforts have been made to either (i) improve the convergence $O(1/\sqrt{N})$, as in quasi-Monte Carlo (which often achieves $O(1/N)$ convergence), or (ii) reduce the constant $\sigma(f)$, as in a variety of techniques for \emph{variance reduction} (e.g.~\cite{glasserman2013monte,owenmcbook}). Another important line of recent work is multilevel~\cite{giles2015multilevel} and multifidelity~\cite{peherstorfer2016optimal} Monte Carlo methods. 



This work is aligned more closely with (ii), but later we discuss MC-like methods that can converge faster than $O(1/\sqrt{N})$ or $O(1/N)$. 
Our approach is nonstandard: we revisit MC from a numerical analyst's viewpoint, and we interpret MC as a classical quadrature scheme in numerical analysis.

A quadrature rule for integration follows 
 the principle: \emph{approximate, and integrate} (this will be the guiding priciple throughout the paper):
\begin{enumerate}
\item  Approximate the integrand with a simple function (typically a polynomial)
$p(\x)\approx f(\x)$,
\item Integrate $I_p = \int_\Omega p(\x)d\x$ exactly. 
\end{enumerate}
$I_p$ is then the approximation to the desired integral $I$. 
Since $p$ is a simple function such as a polynomial, step 2 is usually straightforward. 
In standard quadrature rules (including Gauss and Clenshaw-Curtis quadrature),
 the first step is achieved by finding a polynomial interpolant s.t. $p(\x_i)=f(\x_i)$ at carefully chosen sample points $\{\x_i\}_{i=1}^N$. 
For the forthcoming argument, 
we note that one way to 
obtain the polynomial interpolant $p(\x)=\sum_{j=0}^nc_j\phi_j(\x)$, where 
$\{\phi_j\}_{j=0}^n$ is a basis for polynomials of degree $n$ (for example Chebyshev or Legendre polynomials), is 
to solve the linear system\footnote{\label{foot1}Of course,~\eqref{eq:vanderlinsys} is not solved explicitly; 
the solution $\c$ (or the desired integral $\int_\Omega p(\x)d\x$) 
is given as an explicit linear combination of $f(\x_i)$'s, exploiting the structure of $\V$. 
For example, in Gauss quadrature $\V$ satisfies a 
 weighted orthogonality condition, namely $\V\V^T=\W^{-1}$ is diagonal. Thus the solution $\c$ of~\eqref{eq:vanderlinsys} is equal to that of the least-squares problem $\min_\c\|\sqrt{\W}(\V_{(:,1)}-\f)\|_2$ where $\V_{(:,1)}$ is the first column of $\V$ (i.e., the vector of ones); this allows for a fast and explicit computation of $\c$ and hence $\int_\Omega p(\x)d\x$. 
}
\begin{equation}
  \label{eq:vanderlinsys}
  \begin{bmatrix}
1&\phi_1(\x_1)&\phi_2(\x_1)&\ldots &\phi_{n}(\x_1)\\
1&\phi_1(\x_2)&\phi_2(\x_2)&\ldots&\phi_n(\x_2)\\
\vdots&\vdots&\vdots\\
1&\phi_1(\x_{N})&\phi_2(\x_{N})&\ldots&\phi_n(\x_{N})\\
  \end{bmatrix}
  \begin{bmatrix}
c_0\\c_1\\\vdots \\c_{n}
  \end{bmatrix}
=  \begin{bmatrix}
    f(\x_1)\\f(\x_2)\\\vdots\\ f(\x_{N})
  \end{bmatrix},
\end{equation}
where $N=n+1$. We write~\eqref{eq:vanderlinsys} as $\V \c = \f$, where $\V\in\mathbb{R}^{N\times N}$ with  $\V_{i,j} = \phi_{j-1}(\x_i)$ is the Vandermonde matrix. We emphasize that the sample points $\{\x_i\}_{i=1}^{N}$ 
are chosen deterministically (extrema of the $n$th Chebyshev polynomial in Clenshaw-Curtis, and roots of $(n+1)$th Legendre polynomial in Gauss quadrature). 


In the next section we reveal the interpretation of Monte Carlo integration as a classical numerical integration scheme, which exactly follows the principle displayed above: approximate $f\approx p$, and integrate $p$ exactly. 
This simple but fundamental observation lets us naturally develop integration methods that blend statistics with approximation theory. 
We show that by employing a better approximant $p$, the MC convergence can be improved accordingly, 
resulting in  a reduced variance. 
Furthermore, by 
employing a function approximation method that improves 
the approximation quality $\|f-p\|_2$ as more samples are taken, 
we can achieve asymptotic convergence faster than $O(1/\sqrt{N})$ or even $O(1/N)$. 
Overall, this paper shows that MC can be combined with any method in (the very active field of) high-dimensional approximation theory, resulting in an MC-like integrator that gets the best of both worlds.

This paper is organized as follows. 
In Section~\ref{sec:mcasquad} we present a numerical analyst's interpretation of MC and introduce a new MC-like algorithm which we call Monte Carlo with least-squares (MCLS for short). 
In Section~\ref{sec:MCLSanalysis} we explore the properties of MCLS. 
We then introduce MCLSA, an adaptive version of MCLS in Section~\ref{sec:lsmca} that can converge faster than $O(1/\sqrt{N})$. 
Sections~\ref{sec:leg} and~\ref{sec:spgridmc} present specific examples of MCLS integrators, namely combining MC with  approximation by polynomials and sparse grids. We outline other approximation strategies in Section~\ref{sec:mainalg} and discuss connections to classical MC techniques in Section~\ref{sec:discuss}. 

{\em Notation.}
$N$ denotes the number of samples, and $n$ is the number of basis functions. We always take $N\geq n$. Our MC(LS) estimate using $N$ sample points is denoted by $\hat I_N$. 
Bold-face captial letters denote matrices, 
and  $\I_n$ denotes the $n\times n$ identity matrix. 
Bold-face lower-case letters represent vectors: in particular, 
$\x_1,\ldots,\x_N\in\mathbb{R}^N$ are the sample points, 
and $\f=[f(\x_1),\ldots,f(\x_N)]^T$ is the vector of sample values. 

To avoid technical difficulties we assume throughout that $f:\Omega\rightarrow \mathbb{R}$ is $L_2$-integrable, i.e.,
$\int_{\Omega} f(\x)^2d\x<\infty$. Thus we are justified in taking expectations and variances. All experiments were carried out using MATLAB version 2016b.

\section{Monte Carlo integration as a quadrature rule}\label{sec:mcasquad}
Here we reveal the most important observation in this work, a numerical analyst's interpretation of Monte Carlo: again, approximate and integrate. What is the approximant here? To answer this we state a trivial but important result. 
\begin{lemma}  
\label{lem0}
The solution $c_0$ to the least-squares problem
\begin{equation}
  \label{eq:ls}
\min_{c_0\in\mathbb{R}}
\left\|
  \begin{bmatrix}
    1\\1\\\vdots \\1 
  \end{bmatrix}c_0-
  \begin{bmatrix}
    f(\x_1)\\f(\x_2)\\\vdots\\ f(\x_N)
  \end{bmatrix}
\right\|_2,
\end{equation}
which we write $\min_{c_0}\|\V c_0- \f\|$, 
is $c_0=\frac{1}{N}\sum_{i=1}^N f(\x_i)$, the MC estimator~\eqref{eq:mcdef}. 
\end{lemma}
{\sc proof.} 
Since $\V=[1,...,1]^T\in\mathbb{R}^{N}$ is clearly of full column rank, 
the solution $c_0$ can be expressed  via the normal equation~\cite[Ch.~5]{golubbook4th} as 
\begin{equation}
  \label{eq:equivalence}
c_0=(\V^T\V)^{-1}\V^T\f=\frac{1}{N}\sum_{i=1}^N f(\x_i). 
\end{equation}
\hfill$\square$

This shows that the MC estimator is equal to $c_0$, the integral of the constant function $c_0$, obtained by solving a (extremely tall and skinny) least-squares problem~\eqref{eq:ls} that attempts to approximate $f$ by the constant $c_0$. 
In other words, MC can be understood as follows: 
\begin{enumerate}
\item  Approximate the integrand with a constant function 
$c_0\approx f(\x)$, 
\item Integrate $\hat I_N = \int_\Omega c_0d\x=c_0$ exactly. 
\end{enumerate}
$\hat I_N$ is then taken as an approximation to $I$. 
Contrast this with the principle of quadrature presented above. 
As before, when $\Omega$ is not the unit hypercube we take 
$\hat I_N= c_0|\Omega|$. 

\ignore{
In order to justify the MC integrator~\eqref{eq:mcdef} can be understood the above way, 
we consider the least-squares (LS) problem
\begin{equation}
  \label{eq:ls}
\min_{c_0\in\mathbb{R}}
\left\|
  \begin{bmatrix}
    1\\1\\\vdots \\1 
  \end{bmatrix}c_0-
  \begin{bmatrix}
    f(\x_1)\\f(\x_2)\\\vdots\\ f(\x_N)
  \end{bmatrix}
\right\|_2,
\end{equation}
which is a classical least-squares problem 
in numerical linear algebra: 
$\min_{c_0}\|\V c_0- \f\|$ with 
 $\V=[1,...,1]^T\in\mathbb{R}^{N}$, and $\f=[f(\x_1),\ldots,f(\x_N)]^T$. 
The solution can be expressed  via the normal equation~\cite[Ch.~5]{golubbook4th}
\begin{equation}
  \label{eq:equivalence}
c_0=(\V^T\V)^{-1}\V^T\f=\frac{1}{N}\sum_{i=0}^N f(\x_i),
\end{equation}
which is exactly equal to~\eqref{eq:mcdef}. 
} 

\ignore{
To make the parallel with numerical integration rules more explicit, 
we note that a classical numerical integration rule based on interpolation 
(including Gauss and Clenshaw-Curtis)
can be understood as solving the linear system 
\begin{equation}
  \label{eq:vanderlinsys}
\V \c = \f  
\end{equation}
to obtain the polynomial interpolant $p(x)=\sum_{j=0}^nc_j\phi_j(x)$, which approximates $f(x)$. Here $\V_{i,j} = \phi_{j-1}(x_i)$ is the Vandermonde matrix where
$\{x_i\}_{i=1}^{N}$ with $N=n+1$ are deterministic points (extrema of the $n$th Chebyshev polynomial in Clenshaw-Curtis, and roots of $(n+1)$th Legendre polynomial in Gauss quadrature) and
 $\{\phi_j\}_{j=0}^n$ is a basis for polynomials of degree $n$, for example Chebyshev or Legendre polynomials\footnote{Of course, the linear system~\eqref{eq:vanderlinsys} is not solved explicitly; the solution $p$ (or its integral $\int p$) is given as an explicit linear combination of $f(x_i)$. 
The points $\{x_i\}_{i=1}^{n+1}$ and basis $\{\phi_j\}_{j=0}^n$ are chosen in a careful way such that the Vandermonde matrix has a (weighted) orthogonality structure (namely $\V\V^T$ is diagonal), which allows for a fast and explicit computation of the first element of $\V^{-1}\f$.}. 
}


Comparing~\eqref{eq:lsdisplay} with \eqref{eq:vanderlinsys}, 
one can therefore view MC as a classical quadrature rule where (i) the sample points are chosen randomly, and (ii) the approximant is obtained by a least-squares problem $\min_{\c}\|\V \c - \f\|_2 $ rather than a linear system $\V \c = \f $. Note that linear systems can be regarded as a special case of least-squares problems where $\V$ is square\footnote{Moreover, as mentioned in footnote \ref{foot1}, the linear system \eqref{eq:vanderlinsys} is equivalent to a certain least-squares problem with $\V\in\mathbb{R}^{N\times 1}$.}.  
Also note that the variance estimate~\eqref{eq:sigmaestimate} 
can be written as $\frac{1}{N-1}\|\V c_0-\f\|^2_2$, 
which is essentially the squared residual norm in the least-squares fit.



\subsection{MCLS: Monte Carlo with least-squares}
The significance of this new interpretation of MC is that it naturally suggests an extension, in which the approximant is taken to be non-constant, for example a low-degree polynomial. Put another way, we sample as in Monte Carlo, but approximate as in quadrature rules.
Namely, here is our prototype algorithm, which we call Monte Carlo with least-squares (MCLS)\footnote{
The first name the author thought of is least-squares Monte Carlo, but 
there is a popular method bearing this name 
 in finance for American options~\cite{longstaff2001valuing}. 
}. 
\begin{enumerate}
\item  Approximate the integrand with a 
simple function
$p(\x):=\sum_{j=0}^nc_j\phi_j(\x)\approx f(\x)$,
\item Integrate $\hat I_N = \int_\Omega p(\x)d\x$ exactly. 
\end{enumerate}
Here, $\{\phi_j\}_{j=0}^n$ are prescribed functions, referred to as \emph{basis functions}. 
We always take $\phi_0(\x)=1$, the constant function, so that the process reduces to MC when $n=0$. 
We also assume that 
$\{\phi_j\}_{j=0}^n$ are linearly independent, as otherwise $\V$ is rank deficient. 
The main question is how to obtain $\c=[c_0,c_1,\ldots,c_n]^T$. 
In view of~\eqref{eq:ls} in MC, we do this by solving 
\begin{equation}
  \label{eq:lsdisplay}
\min_{\c\in\mathbb{R}^{n+1}}
\left\|
  \begin{bmatrix}
1&\phi_1(\x_1)&\phi_2(\x_1)&\ldots &\phi_n(\x_1)\\
1&\phi_1(\x_2)&\phi_2(\x_2)&\ldots&\phi_n(\x_2)\\
\vdots&\vdots&\vdots\\
1&\phi_1(\x_N)&\phi_2(\x_N)&\ldots&\phi_n(\x_N)\\
  \end{bmatrix}
  \begin{bmatrix}
c_0\\c_1\\\vdots \\c_n    
  \end{bmatrix}
-  \begin{bmatrix}
    f(\x_1)\\f(\x_2)\\\vdots\\ f(\x_N)
  \end{bmatrix}
\right\|_2,
\end{equation}
which is an $N\times (n+1)$ 
($N>n$, often $N\gg n$)
linear least-squares problem, which as in~\eqref{eq:ls} we express as 
\begin{equation}  \label{eq:Vc-Fmain}
\min_\c\|\V \c-\f\|_2  , 
\end{equation}
 but with (many) more columns than one, employing more basis functions than just the constant function. 

The solution for $\min\|\V\c-\f\|_2$ is again 
$\c = (\V^T\V)^{-1}\V^T\f$. The cost of solving~\eqref{eq:lsdisplay} is $O(Nn^2)$ using a standard QR-based least-squares solver (this can usually be reduced to $O(Nn)$ using a Krylov subspace method; see Appendix~\ref{sec:comp}). This is  the main computational task in MCLS. 

Once the solution $\c =[c_0,\ldots,c_n]^T$ is obtained, the approximate integral is computed as
\begin{equation}  \label{eq:lsmcval}
\hat I_N:=   \int_\Omega p(\x)d\x
 =c_0|\Omega| +\sum_{j=0}^nc_j\int_{\Omega} \phi_j(\x)d\x. 
\end{equation}
Assuming $\int_\Omega\phi_j(\x)d\x$ is known for every $j$, integrating $p$ is an easy task once $\c$ is determined.

\begin{figure}[htbp]
  \begin{minipage}[t]{0.325\hsize}
  \centering
\includegraphics[width=1.0\textwidth]{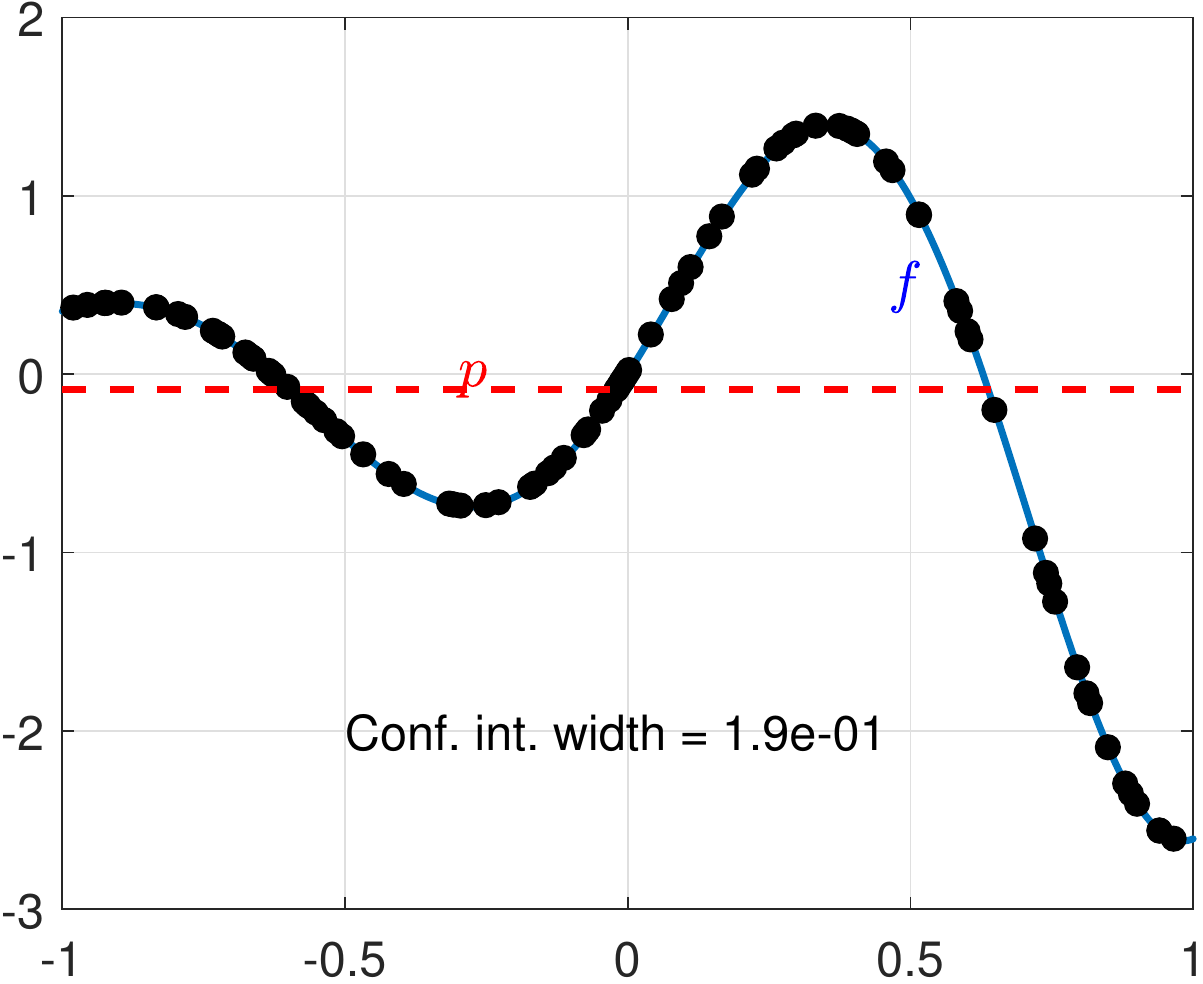}
  \end{minipage}
  \begin{minipage}[t]{0.33\hsize}
  \centering
\includegraphics[width=.98\textwidth]{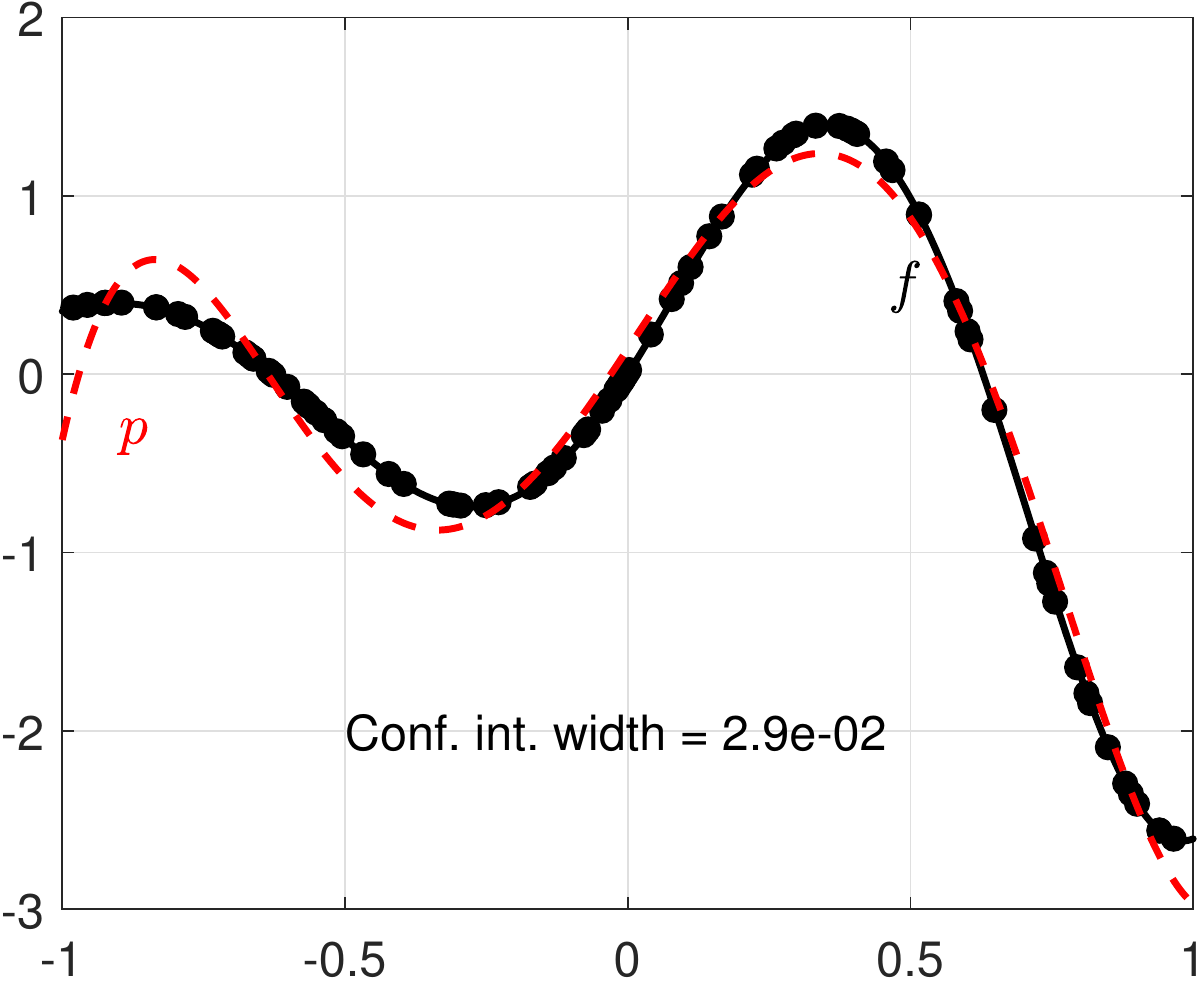}
  \end{minipage}
  \begin{minipage}[t]{0.325\hsize}
  \centering
\includegraphics[width=1.0\textwidth]{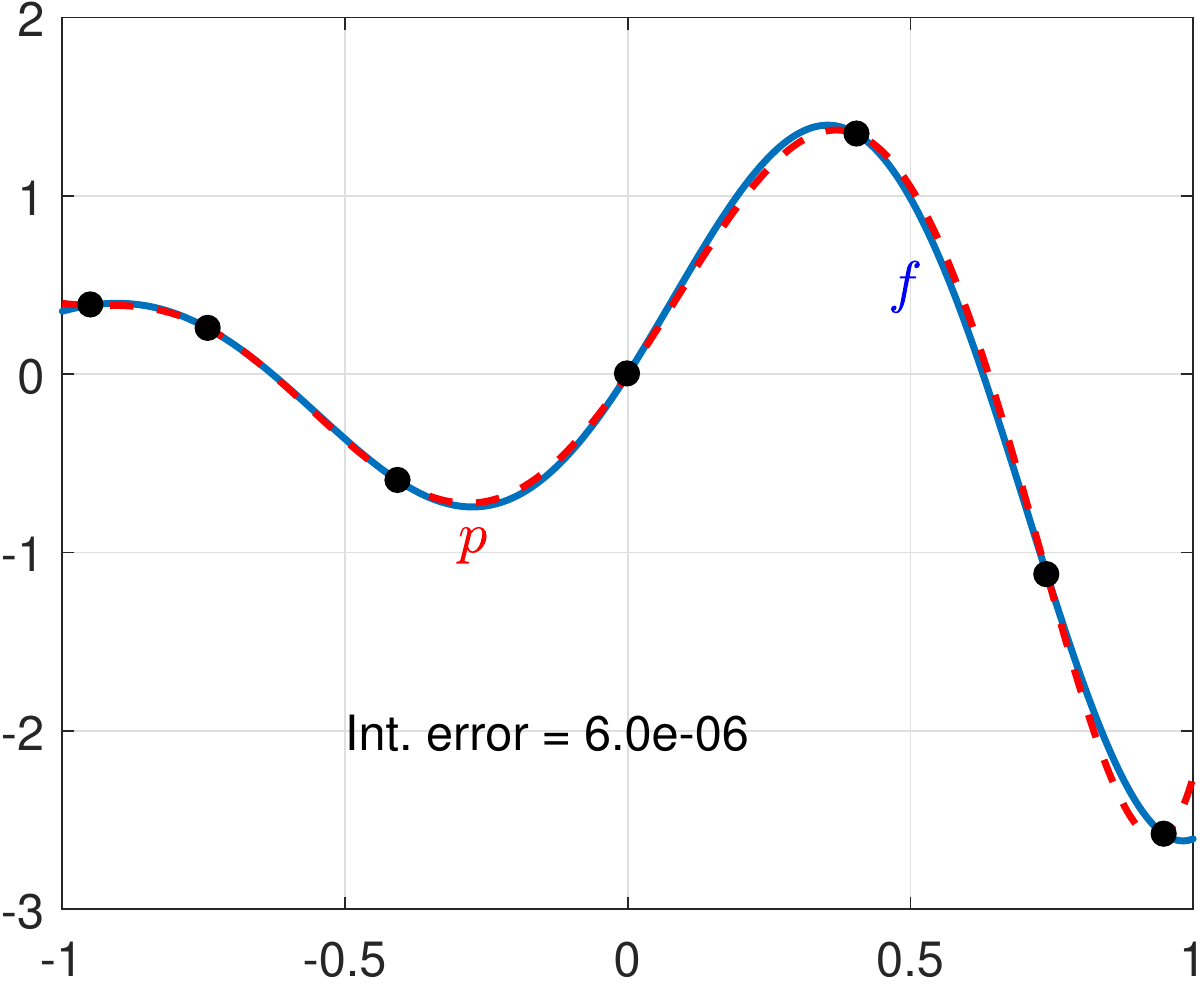}
  \end{minipage}
  \caption{Illustration of approximants (dashed red) underlying integration methods for computing $\int_{-1}^1f(x)dx$ for $f(x)=\exp(x)\sin(5x)$ shown as blue solid curve. Black dots are the sample points. 
MC (left) employs a constant approximation (100 sample points), 
MCLS (center) here uses a polynomial approximation of degree 5,  and  Gauss quadrature (right)  a polynomial interpolant, here degree $6$. 
The inline texts show the widths of the $95\%$ confidence intervals  for MC and MCLS. No confidence interval is available for Gauss, so its integration error is shown. 
}
  \label{fig:plotkeyidea}
\end{figure}

To illustrate the main idea, Figure~\ref{fig:plotkeyidea} shows the sample points and underlying approximation for MC, MCLS and Gauss quadrature for integrating a 1-dimensional function. 
Observe that MCLS uses a better function approximation than MC, and has a smaller confidence interval. 

Algorithm~\ref{alg:MCLS} summarizes the process in pseudocode. 
Clearly, a key component is the choice of basis functions $\{\phi_j\}_{j=1}^n$. 
A natural choice from a numerical analysis viewpoint is 
 polynomials up to some fixed degree $k\leq N$. 
We explore this in Section~\ref{sec:leg} and consider other choices in later sections. 
\begin{algorithm}[h!]
\linespread{1.2}\selectfont
\DontPrintSemicolon
\KwIn{Function $f$, basis functions $\{\phi_j\}_{j=0}^n$, $\phi_0\equiv 1$, integer $N (>n)$.}
\KwOut{Approximate integral $\hat I_N\approx \int_\Omega f(\x)d\x$.}
Generate sample points $\{\x_i\}_{i=1}^N\in\Omega$, uniformly at random. \;
Evaluate $f(\x_i)$, $i=1,\ldots,N$. \;
Solve the least-squares problem~\eqref{eq:lsdisplay} for $\c=[c_0,c_1,\ldots,c_n]^T$.\;
Compute $\hat I_N=c_0|\Omega|+\sum_{j=1}^nc_j\int_\Omega \phi_j(\x)d\x$.\;
\caption{MCLS: Monte Carlo with Least-Squares for approximating 
$\int_\Omega f(\x)d\x$}\label{alg:MCLS}
\end{algorithm}

\ignore{
\bb{In $d$-dimensions, the complexity becomes $O(Nd^{2k})$, which  quickly gets prohibitive, 
except when, say $k=1,2,3$ (this method can nonetheless be quite powerful). 
}

For numerical stability, the recommended basis is the Chebyshev polynomials $\phi_i(x)=T_i(x)$~\cite{trefethenatap} (or an orthogonal polynomial basis such as Legendre or Jacobi).  However, we expect that MCLS would be useful mostly when the basis functions are taken to be low degree (such as $2$ or $3$), in which case this remark is less relevant and the simplest monomials cause little numerical issues.  
The method we advocate chooses $\phi_i$ using the sample values $f(x_i)$, thus $\phi_i$ depends on the function $f$; see Section~\ref{sec:mainalg}. 
}


\subsection{One-dimensional integration by MCLS}\label{sec:1d}
To illustrate the idea further and give a proof of concept, here we 
explore MCLS for one-dimensional integration. 
Although this is not a competitive method compared with e.g., Gauss or Clenshaw-Curtis quadrature, the experiment reveals aspects of MCLS that remain relevant in higher dimensions. 

\begin{figure}[htbp]
  \begin{minipage}[t]{0.5\hsize}
      \includegraphics[width=0.9\textwidth]{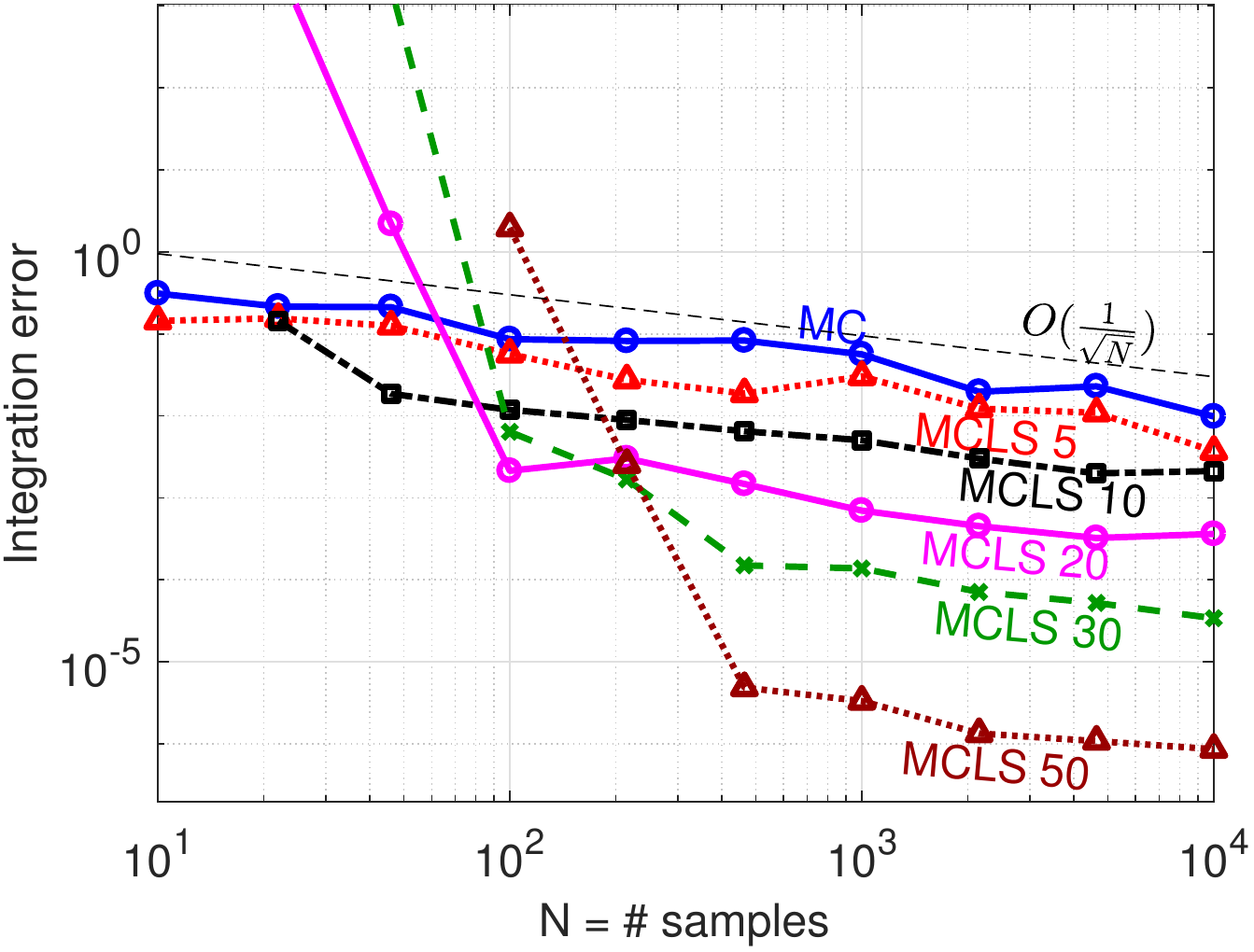}    
  \end{minipage}   
  \begin{minipage}[t]{0.5\hsize}
      \includegraphics[width=0.9\textwidth]{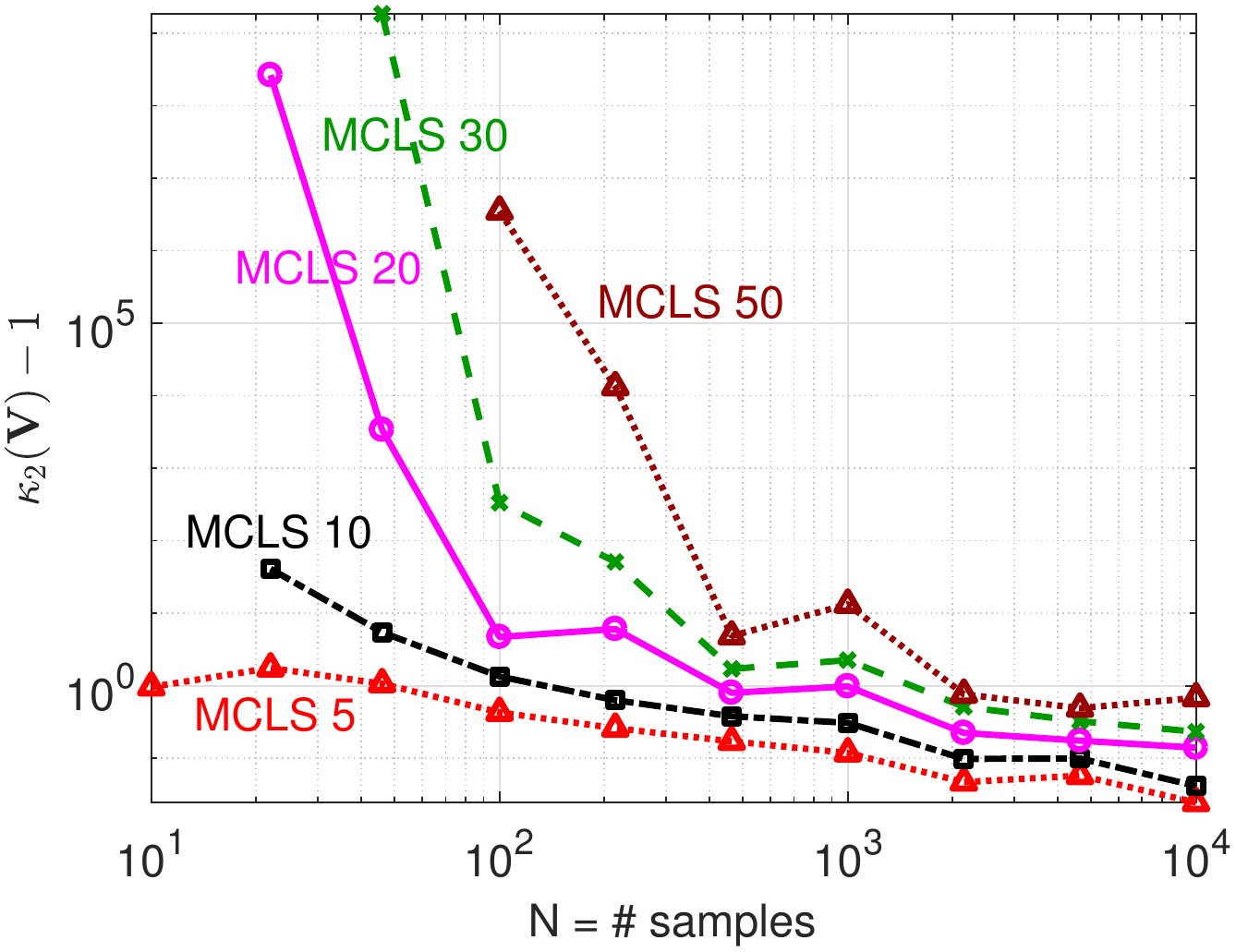}    
  \end{minipage}
  \caption{Left: MC vs. MCLS errors for approximating $\int_{-1}^1\frac{1}{25x^2+1}$, for varying MCLS degrees (``MCLS x'' indicates MCLS with polynomial degree x).  
The curves show the $95\%$ confidence intervals. Right: Conditioning $\kappa_2(\V)-1$ in the least-squares fit (note the minus 1). 
}
  \label{fig:1d}
\end{figure}

Here we consider computing $\int_{-1}^1\frac{1}{25x^2+1}dx$, a classical problem of integrating Runge's function on $[-1,1]$. 
We take the basis functions $\phi_j$ to be Legendre polynomials up to a prescribed degree.  Figure~\ref{fig:1d} (left) shows the error estimates, the width of the $95\%$ confidence intervals (given by~Theorem~\ref{thm:mainvar} in Section~\ref{sec:careful}; here and throughout we plot the confidence interval rather than the error, because confidence intervals converge much more regularly). 
Note that data are shown only when the number of sample points is larger than the degree $N>n$; this is because otherwise the least-squares problem~\eqref{eq:lsdisplay} is underdetermined. 

We make three observations. 
\begin{itemize}
\item Asymptotically for  large $N$, the convergence curves all have similar slopes, exhibiting the same $O(1/\sqrt{N})$ convergence as in standard Monte Carlo. 
\item The constant in front of $1/\sqrt{N}$ decreases as we increase the polynomial degree. 
\item In the pre-asymptotic stage, the error for MCLS behaves erratically. This effect is pronounced for higher degrees. 
\end{itemize}

Each has important ramifications in higher dimensions, so we explain them in more detail. 

The first two observations
are explained in Section~\ref{sec:careful}: MCLS converges like $O(1/\sqrt{N})$ as in MC, but the constant is 
precisely the quality of the approximation $\|p-f\|_2$. 
From a numerical analyst's viewpoint, there is no surprise that the integration accuracy improves as the integrand is approximated better. The MC interpretation is that the variance has been reduced. 

The third observation is a consequence of \emph{ill-conditioning}; namely the coefficient matrix $\V$ in the least-squares problem is ill-conditioned 
if the number of sample points is insufficient. 
 To further illustrate this, in Figure~\ref{fig:1d} (right) we plot $\kappa_2(\V)-1$; the subtraction by $1$ is done to let us see how well-conditioned they become. We see the erratic behavior in convergence is observed only when $\kappa_2(\V)\gg 1$. 

The effect of ill-conditioning is never present in standard MC, since $\V\in\mathbb{R}^{N\times 1}$ always has condition number 1; 
by contrast, it is omnipresent  in numerical analysis~\cite{Higham:2002:ASNA}, and is often a central object to examine. We discuss alleviating conditioning of $\V$ in Section~\ref{sec:leghigh}. 


\subsection{Monte Carlo integration: two viewpoints}\label{sec:desire}
Monte Carlo is an extremely practical and popular method for high-dimensional integration. 
Above, we have interpreted Monte Carlo as a quadrature rule, which, while elementary, 
is at the very heart of this work. 
From the viewpoint of classical MC, the novelty lies in the observation that there is a very simple function approximation underlying the method. 

The hallmark and big surprise of MC---
from a numerical analyst's viewpoint---is that, even though $p$ is poor (in fact not converging at all) as an approximant to $f$, its integral still converges to that of $f$\footnote{\label{foot:4}In some sense this phenomenon manifests itself in classical 1-D quadrature rules. Gauss quadrature, for example, finds 
the polynomial $p$ that interpolates $f$ at the roots of the Legendre polynomial of degree $n+1$. We clearly have $p=f$ if $f$ is a polynomial of degree up to $n$; clearly $p\neq f$ if $f$ is of higher degree. Nonetheless, the integral $\int_\Omega pdx$ is equal to $\int_\Omega fdx$ if $f$ is a polynomial of double degree (plus one), $2n+1$. What is taking effect here is not central limit theorem but a phenomenon called aliasing~\cite[Sec.~8]{trefethen2014exponentially}. 
We will encounter this ``integral is more accurate than $\|f-p\|_2$'' effect several times. 
}. 
One might also wonder if a better integration can be obtained by using a more sophisticated approximation. For example, for functions analytic in $\Omega$, polynomials are (at least in the limit $N\rightarrow\infty$) able to approximate with superalgebraic convergence (which means the asymptotic convergence is faster than $O(N^{-c})$ for any constant $c>0$). 

In light of these, 
one might wish 
to have a numerical integration method that 
achieves the following desiderata: 
\begin{itemize}
\item It takes advantage of the $O(1/\sqrt{N})$ convergence as in Monte Carlo. 
\item It approximates the function $f$ as much as possible, and the quality of the estimate reflects the approximation quality. 
\item When used for integrating analytic functions, the asymptotic convergence is superalgebraic. 
\end{itemize}

MCLS as presented in Algorithm~\ref{alg:MCLS} satisfies the first two conditions, but not quite the third. 
We investigate the convergence properties of MCLS in the next section, 
and then introduce variants of MCLS that achieve all three requirements. 

\section{Properties of MCLS}\label{sec:MCLSanalysis}
In this section we explore the theoretical 
 properties of MCLS. 
We first investigate the error in the MCLS estimator and reveal the intimate connection between the variance and the approximation quality $\|f-p\|_2$. 
In Section~\ref{sec:bias} we examine the bias of 
the MCLS estimator, which is nonzero but negligible relative to the integration error. 
To keep the focus on the statistical aspects of MCLS, the computational aspects are deferred to Appendix~\ref{sec:comp}, where we discuss updating the solution as more samples are taken, and fast $O(Nn)$ solvers. 


\subsection{Convergence of MCLS}\label{sec:careful}
To gain insight, consider the following: apply classical MC to evaluate  $\int_\Omega (f(\x)-p(\x))d\x$, with the same sample points $\{\x_i\}_{i=1}^N$ and $p$ is the solution for~\eqref{eq:lsdisplay}. The outcome is $0$, because of the least-squares fit, which imposes $\V^T(\V\c-\f)=0$; the first element is precisely $\sum_{i=1}^N(p(\x_i)-f(\x_i))=0$. Classical MC analysis~\eqref{eq:sigmaestimate} shows that the variance is that of the function $f-p$, divided by $N$. Thus we expect the estimator $\hat I_N$ to converge to $I$ with error $O(\sigma(f-p)/\sqrt{N})$. Since $\sigma(f-p)\leq \|f-p\|_2$, we expect convergence like $\|f-p\|_2/\sqrt{N}$  (we expect the inequality here to be sharp since $\phi_0=1$ is among the basis functions).

Care is needed to make a rigorous argument, because the 
approximant  $p$ clearly depends on the sample points $\{\x_i\}_{i=1}^N$. 
Nevertheless, the above informal argument gives the correct asymptotics: 
  \begin{theorem}\label{thm:mainvar}
Fix $n$ and the basis functions $\{\phi_j\}_{j=0}^n$. 
Then with the MCLS estimator $\hat I_N$ in~\eqref{eq:lsmcval}, 
as $N\rightarrow \infty$ we have  
\begin{equation}
  \label{eq:MCLSestimate}
\sqrt{N}(\hat I_N-I)\xrightarrow{d} \mathcal{N}(0,    \min_\c\|f-\sum_{j=0}^nc_j\phi_j\|_2^2) ,
\end{equation}
where 
$\xrightarrow{d}$ denotes convergence in distribution. 
  \end{theorem}
{Proof.} 
We assume that the basis functions $\{\phi_j\}_{j=0}^n$ form an orthonormal basis, i.e., 
$\int_\Omega \phi_i(x)\phi_j(x)dx=\delta_{ij}$, the Kronecker delta function. This 
simplifies the analysis, and can be done without loss of generality, as the least-squares problem~\eqref{eq:lsdisplay} gives the same solution $p$ 
for another basis $\{\tilde\phi_j\}_{j=0}^n$
as long as $\mbox{span}(\{\phi_j\}_{j=0}^n)=\mbox{span}(\{\tilde\phi_j\}_{j=0}^n)$. An orthonormal basis can be obtained  using e.g. Gram-Schmidt orthogonalization. 

We decompose $f$ into a sum of orthogonal terms
\begin{equation}
  \label{eq:fdecomp}
f = 
\sum_{j=0}^nc_j^\ast\phi_j+g
=:f_1+g, 
\end{equation}
where $g$ is a function orthogonal to all the basis functions $\{\phi_j\}_{j=0}^n$, including the constant $\phi_0=1$, that is, $\int_\Omega g(\x)\phi_j(\x)d\x=0$. Note that $\|g\|_2=\min_\c\|f-\sum_{j=0}^nc_j\phi_j\|_2$. 
The vector of sample values is 
\[
\f= 
[\sum_{j=0}^nc_j^\ast\phi_j(\x_1)+g(\x_1),\ldots,\sum_{j=0}^nc_j^\ast\phi_j(\x_N)+g(\x_N)]^T
=\V\c^\ast+\g. 
\]
Denoting by $\hat \c$ the least-squares solution to~\eqref{eq:lsmcval}, 
we have 
\[\hat \c=\displaystyle\argmin_{\c} \|\V\c - (\V\c^\ast+\g)\|_2 =
\argmin_{\c} \|\V(\c-\c^\ast) -\g \|_2 =
 \c^\ast+ \c_g,\]
 where 
$\c_g=\mbox{argmin}_{\c} \|\V\c - \g\|_2 = (\V^T\V)^{-1}\V^T\g$. 
It thus follows that 
$\hat I_n-I = c_{g,0}=[1,0,\ldots,0](\V^T\V)^{-1}\V^T\g$. 

Now  by the strong law of large numbers  we have 
\[\frac{1}{N}(\V^T\V)_{i+1,j+1} = 
\frac{1}{N}\sum_{\ell=1}^N\phi_i(\x_\ell)\phi_j(\x_\ell)
\rightarrow 
\int_{\Omega}\phi_i(\x)\phi_j(\x)d\x
=\delta_{ij}\]
almost surely as $N\rightarrow \infty$, by the orthonormality of $\{\phi_j\}_{j=0}^n$. Therefore we have 
$\frac{1}{N}\V^T\V\rightarrow \I_{n+1}$ 
 as $N\rightarrow \infty$, so 
 \begin{equation}
   \label{eq:cg0}
\sqrt{N}c_{g,0}=
\sqrt{N}[1,0,\ldots,0]^T
\frac{1}{N}(\frac{1}{N}\V^T\V)^{-1}\V^T\g
\rightarrow \sqrt{N}\left(\frac{1}{N}\sum_{i=1}^Ng(\x_i)\right)
\xrightarrow{d} \mathcal{N}(0,\|g\|^2)   
 \end{equation}
by the central limit theorem. 
\hfill$\square$



Theorem~\ref{thm:mainvar} shows the MCLS estimator gives an approximate integral $\hat I_N$ such that 
\begin{equation}  \label{eq:err}
\E(|\hat I_N-I|) \approx \frac{\min_\c\|f-\sum_{j=0}^nc_j\phi_j\|_2}{\sqrt{N}} . 
\end{equation}
Contrast this with the error with classical Monte Carlo
$\frac{1}{\sqrt{N}}\min_c\|f-c\|_2$ in~\eqref{eq:MCvar}: 
the asymptotic error is still $O(1/\sqrt{N})$, but 
the variance has been reduced 
from $\min_c\|f-c\|_2^2$ 
to $\min_\c\|f-\sum_{j=0}^nc_j\phi_j\|_2^2$. 
In both cases, the variance 
is precisely the squared $L_2$-norm of the error in the approximation 
$f\approx \sum_{j=0}^nc_j\phi_j$ (in  MC, $f\approx c_0$). In other words, \emph{the constant in front of the $O(1/\sqrt{N})$ convergence in MC(LS) is equal to the function approximation error} (in the $L_2$ norm).

Note that 
if $f$ lies in $\mbox{span}\{\phi_j\}_{j=0}^n$, 
the variance 
$\min_\c\|f-\sum_{j=0}^nc_j\phi_j\|_2^2$ becomes zero. 
This means the MCLS estimate will be exact (assuming $\V$ is of full column rank; this holds almost surely if $N>n$). 
This claim can be verified in a straightforward manner by noting that 
$\f = \V\c^*$ for the exact coefficient vector $\c^*$, regardless of the sample points $\{\x\}_{i=1}^N$, and $\c^*$ is the unique solution for $\f = \V\c^*$ if $\V$ is full rank, giving $p=\sum_{j=0}^nc^*_j\phi_j=f$. 

In practice in  MCLS, once~\eqref{eq:lsdisplay} has been solved, the variance
$\min_\c\|f-\sum_{j=0}^nc_j\phi_j\|_2^2$
 can be estimated via
\begin{equation}  \label{eq:sigmaestimatels}
\widetilde\sigma_{LS}^2 = \frac{1}{N-n-1}\sum_{i=1}^N(f(\x_i)-p(\x_i))^2=\frac{1}{N-n}\|\V \c-\f\|^2_2, 
\end{equation}
as is commonly done in linear regression. As in MC, the variance estimate is proportional to the squared residual norm in the least-squares fit $\min_\c\|\V\c-\f\|_2$. 

\subsubsection{Convergence comparison}
We briefly return to classical cubature rules. 
As mentioned previously, these are based on approximating the integrand $p\approx f$. The error can thus be bounded as 
\begin{equation}  \label{eq:interr}
\big|\int_\Omega f(\x)d\x-\int_\Omega p(\x)d\x  \big| \leq \int_\Omega |f(\x)- g(\x)|d\x  =\|f-g\|_1\leq\|f-g\|_2,
\end{equation}
where we used the Cauchy-Schwarz  inequality and the fact $\Omega=[0,1]^d$ for the final  inequality. 
However, it is important to keep in mind that in cubature methods, the integration error often comes out much better than the approximation error. We mentioned this for Gauss quadrature in footnote~\ref{foot:4}, and we revisit this phenomenon in Section~\ref{sec:spgridmc}. 

We summarize the comparison between MC, MCLS and classical cubature in Table~\ref{tab:compare}. For the reason just described, the bottom-right entry is an oversimplification, and should be regarded as an (often crude) upper bound. 

\begin{table}[htbp]
  \centering
  \caption{Comparison of integration methods: Monte Carlo (MC), Monte Carlo with least-squares  (MCLS), and quadrature/cubature.
$N$ is the number of sample points, and 
$C_f$ denotes the cost for evaluating $f$ at a single point. 
}
  \label{tab:compare}
  \begin{tabular}{c|c|c|c}
& Computation & Cost & Convergence \\\hline
MC     & 
$\min_{c}\|\V c-\f\|_2$
 & $C_fN$ & $\displaystyle\frac{1}{\sqrt{N}}\min_{c}\|f-c\|_2$\\
MCLS &  $\min_{\c}\|\V\c-\f\|_2$ & $C_fN+O(Nn)$ & $\displaystyle\frac{1}{\sqrt{N}}\min_{\c}\|f-\sum_{j=0}^nc_j\phi_j\|_2$\\
 cubature & $\V\c=\f$ & $C_fN$ 
& $\displaystyle\min_{\c}\|f-\sum_{j=0}^nc_j\phi_j\|_2$
  \end{tabular}
\end{table}


Table~\ref{tab:compare} highlights a number of aspects worth mentioning. First, all three methods (explicitly or implicitly) perform a basic linear algebra operation: least-squares problem or linear system. Hence in a broad sense, they can all be understood as members of the same family of methods that perform linear approximation followed by integration. 
Second, the cost for MCLS is higher than for MC with the same number of sample points $N$. However, in many applications, evaluating $f$ is expensive $C_f\gg 1$, and the error reduction with MCLS may well justify the extra $O(Nn)$ cost (for example, when $n\ll C_f$ there is effectively no overhead). 
Last but not least, MCLS performs the ``function approximation'' as in cubature, while maintaining the $1/\sqrt{N}$ convergence in MC.  In this sense it gets the best of both worlds. 




\subsection{Bias of the MCLS estimator}\label{sec:bias}
The MCLS estimator is actually biased, that is, $\E(I_p)\neq I$, 
where the expectation is taken over the random samples $\{\x_{i}\}_{i=1}^N$, where $N$ 
and $\{\phi_j\}_{j=0}^n$
are  fixed. 
The bias is nonetheless harmless, as we show below. The line of argument closely follows that of Glynn and R. Szechtman~\cite[Thm.~1]{glynn2002some} and Owen~\cite[Eqn.~(8.34)]{owenmcbook}, who analyze MC with control variates (see Section~\ref{sec:controlvariates} for its connection to MCLS). 


\ignore{
The error of the MCLS estimator is 
\[
\bar f - c_0 = \bar f - [1,0,0,\ldots,0](V^TV)^{-1}V^TF
=\bar f - [1,0,0,\ldots,0]\c.
\]
Taking expectations (with respect to the sample points), we obtain 
\[
\E(\bar f - c_0) = \bar f - [1,0,0,\ldots,0]\E((V^TV)^{-1}V^TF). 
\]
Now consider the MC estimate obtained from the ``correct'' values $c^*$, as in~\eqref{eq:fdecomp}. 
Then we clearly have $\bar f = c_0^*$, and the resulting MC estimate is unbiased: $\E(\bar f-\V c^*)=0$. And the difference from 
\[\mbox{mean}(F-Vc_*)-\mbox{mean}(F-Vc) = \mbox{mean}(V)(c_*-c). \]
Here $\mbox{mean}(V)$ denotes the row of column-wise mean values of $V$. 
Now note that $\mbox{mean}(V)$ and $(c_*-c)$ both scale like $O(1/\sqrt{N})$ in norm. Therefore the product $\mbox{mean}(V)(c_*-c)$ scales like $O(1/N)$; thus the bias is $O(1/N)$, which is negligible compared with the expected error in the Monte Carlo estimate $O(1/\sqrt{N})$. This carries over to faster converging methods such as quasi Monte Carlo; essentially $c_*-c$ is converging at the same rate as $\|F-Vc\|$, and the bias is smaller by $\|\mbox{mean}(V)\|$, which is $O(1/\sqrt{N})$ or smaller. 
}

\begin{proposition}
With the MCLS estimator $\hat I_N$ with $n$ and $\{\phi_j\}_{j=0}^n$ fixed,
\begin{equation}
  \label{eq:bias}
|I-\E(\hat I_N)|=O\bigg(\frac{\|f-\sum_{j=0}^nc_j^*\phi_j\|_2}{N}\bigg),
\end{equation}
where $\c^*=[c_0^*,\ldots,c_n^*]^T$ are the exact coefficients as in~\eqref{eq:fdecomp}. 
\end{proposition}
{\sc proof}. 
As in Theorem~\ref{thm:mainvar}, for simplicity we assume 
$\int_\Omega\phi_i(\x)\phi_j(\x)d\x=\delta_{ij}$. 
For a fixed set of sample points $\{\x_i\}_{i=1}^N$, 
once the solution $\hat \c$ of~\eqref{eq:lsdisplay} is computed, the MCLS estimator is 
equal to the MC estimator for $f-\sum_{j=1}^n\hat c_j\phi_j$ 
(note the summand starts from $1$ not $0$), that is, 
the MCLS estimator is 
\[
\hat I_N=
\frac{1}{N}\sum_{i=1}^N(f(\x_i)-\sum_{j=1}^n\hat c_j\phi_j(\x_i))
=\frac{1}{N}\sum_{i=1}^Nf(\x_i)-\sum_{j=1}^n\hat c_j\frac{1}{N}\sum_{i=1}^N\phi_j(\x_i)=:\bar{\f}-\sum_{j=1}^n\hat c_j\bar\phi_j, 
\]
where $\bar\phi_j=\frac{1}{N}\sum_{i=1}^N\phi_j(\x_i)$ is the average value of $\phi_j(\x_i)$. By the central limit theorem we have $\sqrt{N}\bar\phi_j\xrightarrow{d} \mathcal{N}(0,\|\phi_j\|^2)$. 
Now consider a standard MC estimator 
applied to $f-\sum_{j=1}^nc_j^*\phi_j$ with  the exact coefficients $c_j^*$, 
i.e., $\hat I_N^*=\bar{\f}-\sum_{j=1}^nc_j^*\bar\phi_j$. 
This is clearly an unbiased estimator $\E[\hat I_N^*]=I$. 
Therefore the difference between the two estimates is 
\begin{equation}  \label{eq:INdiff}
\hat I_N-\hat I_N^*=\sum_{j=1}^n(\hat c_j- c_j^*)\bar\phi_j.   
\end{equation}
Now arguing as in~\eqref{eq:cg0} to examine the $j$th element $\sqrt{N}c_{g,i}=\sqrt{N}(\hat c_j- c_j^*)$, we obtain 
$\sqrt{N}(\hat c_j- c_j^*)\xrightarrow{d} \mathcal{N}(0,\|\phi_jg\|^2)   $
as $N\rightarrow\infty$. 
Thus both terms in the right-hand side of~\eqref{eq:INdiff} are converging to the mean-zero normal distribution with variance $O(1/\sqrt{N})$. 
Thus by Cauchy-Schwarz 
we obtain 
\[
|\E[(\hat c_j- c_j^*)\bar\phi_j]|\leq 
\sqrt{\E[(\hat c_j- c_j^*)^2]}\sqrt{\E[\bar\phi_j^2]}\approx \frac{1}{N}\|\phi_jg\|_2\|\phi_j\|_2. 
\]
we conclude that the expected value of~\eqref{eq:INdiff} is 
\begin{align}
|\E[\hat I_N-\hat I_N^*]| &= I-\mathbb{E}[\hat I_N] = 
\sum_{j=1}^n\frac{1}{N}\|\phi_jg\|_2\|\phi_j\|_2\nonumber\\
&=O\left(\frac{\|g\|_2}{N}\right)
=O\left(\frac{\|f-\sum_{j=0}^nc_j^*\phi_j\|_2}{N}\right).   \label{eq:INdiffexp}
\end{align}
\hfill$\square$


Comparing~\eqref{eq:bias} and~\eqref{eq:MCLSestimate}, we see 
that  the bias of the MCLS estimator is 
 small relative to the MCLS error estimate~\eqref{eq:MCLSestimate}, suggesting that the bias would not cause issues in practice. 


Nevertheless, if an unbiased estimator is of crucial importance, as Owen describes in~\cite[Sec.~8.9]{owenmcbook}, 
one can perform a two-stage approach: in the pilot stage, obtain $p$ as in MCLS using some (say $N_1=N/2$) of the samples, then use standard Monte Carlo to estimate $\int_\Omega (f-p)d\x $ using the remaining $N-N_1$ samples, and add the exact integral of $p$. The error estimate then becomes $\|f-p\|_2/\sqrt{N-N_1}$. 
We do not pursue this further, because our goal is to use as many sample points as possible to obtain a good approximant $p$. 


%



\section{MCLSA: Approximate and integrate}\label{sec:lsmca}
We have shown that MC(LS) converges like 
$\frac{1}{\sqrt{N}}\|f-p\|_2$, 
as summarized in Table~\ref{tab:compare}. 
In MC $p$ is a constant, and in MCLS $p$ is a linear combination of basis functions $\{\phi_j\}_{j=0}^n$. 
For a fixed set of basis functions, the asymptotic convergence of MCLS is 
$O(1/\sqrt{N})$, the same as MC (though crucially with a smaller constant). 

A natural question arises: can we do better? That is, can we develop MCLS-based methods that asymptotically converge faster than $1/\sqrt{N}$ (or even $1/N$ as in quasi-Monte Carlo)? To answer this, we first need to understand why MCLS is limited to the $O(1/\sqrt{N})$ asymptotic convergence. An answer
is that if $\{\phi_j\}_{j=1}^n$ is fixed, 
the approximation quality $\|f-p\|_2$ does not improve beyond $\min_{\c}\|f-\sum_{j=0}^nc_j\phi_j\|_2$, no matter how many samples are taken---despite the fact that, as we sample more, our knowledge of $f$ clearly improves, and so does our ability to approximate it. 




In view of this, here is a natural idea: refine the quality of the approximant $p$ as $N$ grows, so that $\|f-p\|_2$ decays with $N$. 
In MCLS, this means the basis functions $\{\phi_j\}_{j=0}^n$ are chosen adaptively with $N$. One approach is to increase $n$ with $N$ to enrich $\mbox{span}\{\phi_j\}_{j=0}^n$, which we pursue in Section~\ref{sec:leg}. Another is to fix $n$ but refine $\phi_j$, which we describe in Section~\ref{sec:spgridmc}. 

\ignore{
\begin{algorithm}[h!]
\linespread{1.2}\selectfont
\DontPrintSemicolon
\KwIn{Function $f$, 
approximation scheme $\mathcal{L}$ s.t. 
$\mathcal{L}(\{f(\x_i)\}_{i=1}^N)=p\approx f$ 
}
\KwOut{Approximate integral $\hat I_N\approx \int_\Omega f(x)dx$.}
Generate sample points $\{x_i\}_{i=1}^N\in\Omega$, uniformly at random. \;
Evaluate $f(x_i)$, $i=1,\ldots,N$. \;
Solve the least-squares problem~\eqref{eq:ls} for $\c=[c_0,c_1,\ldots,c_n]^T$.\;
Compute $\hat I_N=c_0|\Omega|+\sum_{j=1}^nc_j\int_\Omega \phi_j(x)dx$.\;
\caption{MCLSA: Monte Carlo with Least-Squares with Adaptive approximation }\label{alg:MCLSadapt}
\end{algorithm}
}

We call such algorithms \emph{MCLSA} (A standing for adaptive) to highlight the adaptive nature of the basis functions. 
MCLSA
therefore  takes advantage of  both  
statistical and analytic
improvements: by sampling more, it enjoys the MC-like $1/\sqrt{N}$ convergence, together with an improved approximation quality $\|f-\sum_{j=0}^nc_j\phi_j\|_2$. 
For example, if $\|f-\sum_{j=0}^nc_j\phi_j\|_2$ converges like $O(1/N^\alpha)$ for some $\alpha>0$, then MCLSA converges like $O(1/N^{\alpha+\frac{1}{2}})$. 




\ignore{
We have seen above that with MCLS, the resulting variance is reduced from $\|f-\bar f\|^2\leq \|f\|^2$ to 
\begin{equation}  \label{eq:varredest}
\|\hat f-\bar{\hat f}\|^2\leq \|\hat f\|^2 = \|f-g\|_2^2 .
\end{equation}
 Note that we expect the last inequality to be sharp (i.e., equality nearly holds), because we expect $(\bar{\hat f})^2$ to be small, scaling like $\|\hat f\|^2/N$. 
It therefore follows that the quality in the function approximation
$\|\hat f\|_2^2 = \|f-g\|_2^2$ is directly reflected in the variance in MCLS; 
we verify this in our experiments.} 

\section{MCLS(A) with polynomial approximation}\label{sec:leg}
In the remainder of the paper we describe several MCLS-based algorithms for integrating~\eqref{eq:goal}, differing mainly in how the 
basis functions $\{\phi\}_{j=0}^n$ are chosen. 
In this section, the sampling strategy (non-uniform sampling) will also play an important role. 

\ignore{
Consider the Legendre expansion (assuming it exists)
\begin{equation}
  \label{eq:fleg}
 f(\x)  = \sum_{j_1,\ldots,j_d}c_{j_1\ldots,j_d}^*P_{j_1}(x_1)\ldots P_{j_d}(x_d). 
\end{equation}
Here $\x=[x_1,\ldots,x_d]^T$. } 

In this section we pursue the perhaps most natural idea from a numerical analysis perspective: use  polynomial basis functions, in particular tensor-product Legendre polynomials. 
We aim to approximate $f$ by a polynomial  of total degree $k$ 
\begin{equation}
  \label{eq:gleg}
 p(\x)  = \sum_{j_1+j_2+\cdots+j_d\leq k}c_{j_1\ldots,j_d}P_{j_1}(x_1)\ldots P_{j_d}(x_d). 
\end{equation}
The number of basis functions (i.e., the number of terms in the summand) is  $_{d+k}C_k=:n+1$. 

As before, to obtain the coefficients $c_{j_1\ldots,j_d}$ we solve the $N\times (n+1)$ least-squares problem~\eqref{eq:ls}, 
with the basis functions 
$\phi_{j}(\x)=P_{j_1}(x_1)\ldots P_{j_d}(x_d)$ for $j=1,\ldots,n$; 
the ordering of the $n$ terms affects neither the approximant $p$ 
nor $\hat I_N$. 
This is a (high-dimensional) approximation problem via discrete least-squares approximation, where the columns of the Vandermonde matrix are discrete samples of an orthonormal set of basis functions. 




\subsection{Optimally weighted least-squares polynomial fitting}\label{sec:leghigh}
Function approximation via discrete least-squares~\eqref{eq:ls} is a classical approach in approximation theory, and in particular was analyzed 
 in the significant paper by Cohen, Davenport and Leviatan~\cite{cohen2013stability}. 
They show that if sufficiently many samples are taken, then $\V$ becomes well-conditioned with high probability (owing to a matrix Chernoff bound~\cite{tropp2012user}), and $p$ is close to the best possible approximant: 
$\E\|f-p\|_2^2$ is bounded  by $(1+\epsilon(N))\|f-f_k\|_2^2$ plus a term decaying rapidly with $N$, where $\epsilon(N)=O(1/\log n)$. 
An issue is that with a uniform sampling strategy in $\Omega=[0,1]^d$, $N$ is required to be as large as $N=O(n^2)$ to  obtain $\kappa_2(\V)=O(1)$ with high probability (this bound is sharp when $d=1$ and becomes less so as $d$ grows, assuming the use of the total degree; with the maximum degree it is sharp for all $d$). We investigate the effect of $\kappa_2(\tilde\V)>1$ on MCLS in Section~\ref{sec:kappa}. 

The work~\cite{cohen2013stability} was recently revisited by Cohen and Migliorati~\cite{cohenoptimal}, who show that the $N=O(n^2)$ obstacle can be improved to $N=O(n\log n)$, if the sampling strategy is chosen appropriately. Similar findings have been reported by 
Hampton and Doostan~\cite{hampton2015coherence} 
and Narayan, Jakeman and Zhou \cite{narayan2017christoffel}. 
Specifically, define the nonnegative function 
$w$ via 
\begin{equation}
  \label{eq:optsample}
  \frac{1}{w(\x)}=\frac{\sum_{j=0}^n\phi_{j}(\x)^2}{n+1}. 
\end{equation}
$\frac{1}{w}$ is a well-defined probability distribution, since $w> 0$ on $\Omega$ and $\int_\Omega \frac{1}{w(\x)}d\x=1$. 
The function $\frac{w(\x)}{n+1}=(\sum_{j=0}^n\phi_{j}(\x)^2)^{-1}$ is the so-called Christoffel functions, which are extensively studied in the literature of orthogonal polynomials~(e.g.~\cite{nevai1986geza}).

Having defined $w$, 
 we take samples $\{\tilde\x_i\}_{i=1}^N$ according to $\frac{1}{w}$, that is, we sample more often where $\sum_{i=0}^n\phi_{i}(\x)^2$ takes large values. 
\ignore{
The Monte Carlo estimate for $\int_{\Omega} f(\x)d\x=\int_{\Omega} w(\x)f(\x)\frac{d\x}{w(\x)}$ from these samples  is 
$\frac{1}{N}\sum_{i=1}^Nw(\tilde\x_i)f(\tilde\x_i)$, which is the solution for the least-squares problem 
$\min_{\c}\|\sqrt{\W}(\V \c-\f)\|_2$, where 
$\sqrt{\W} = \mbox{diag}(\sqrt{w(\tilde\x_1)},\sqrt{w(\tilde\x_2)},\ldots, \sqrt{w(\tilde\x_N)})$ and $\V=[1,\ldots,1]^T$ as was in~\eqref{eq:ls}. 
}
Since 
\[\|f-\sum_{j=0}^nc_j\phi_j\|_2=\int_{\Omega} (f(\x)-\sum_{j=0}^nc_j\phi_j(\x))^2d\x=\int_{\Omega} w(\x)(f(\x)-\sum_{j=0}^nc_j\phi_j(\x))^2\frac{d\x}{w(\x)},\]
the discretized least-squares fitting 
for $\min_{\c}\|f-\sum_{j=0}^nc_j\phi_j\|_2$ 
using $\{\x_i\}_{i=1}^N$
(and hence the core part of MCLS) is 
\begin{equation}  \label{eq:MCLSweight}
\min_{\c}\|\sqrt{\W}(\V\c-\f)\|_2,
\end{equation}
where $\sqrt{\W} = \mbox{diag}(\sqrt{w(\tilde\x_1)},\sqrt{w(\tilde\x_2)},\ldots, \sqrt{w(\tilde\x_N)})$, and $\V,\f$ are as before in~\eqref{eq:Vc-Fmain} with $\x\leftarrow \tilde\x$.  
This is again a least-squares problem $\min_\c\|\tilde\V\c-\tilde\f\|$, with coefficient matrix $\tilde\V:=\sqrt{\W}\V$ and right-hand side $\tilde\f := \sqrt{\W}\f$, whose solution is
$\c=(\tilde\V^T\tilde\V)^{-1}\tilde\V^T\sqrt{\W}\f$.
It is the matrix $\tilde\V$ that is well conditioned with high probability, provided that $N\gtrsim n\log n$. Note that the left-multiplication by $\sqrt{\W}$ forces all the rows of $\tilde\V$ to have the same norm (here $\sqrt{n+1}$); this coincides with a commonly employed strategy in numerical analysis to reduce the condition number, which is known to minimize $\kappa_2({\bf D}\V)$ over diagonal ${\bf D}$ up to a factor $\sqrt{N}$~\cite[\S 7.3]{Higham:2002:ASNA}. 
In our context, the significance is that a well-conditioned Vandermonde matrix $\kappa_2(\tilde\V)=O(1)$ can be obtained for essentially as large $n$ as possible, thus reducing the approximation error $\min_{\c}\|f-\sum_{j=0}^nc_j\phi_j\|_2$ (recall Table~\ref{tab:compare}). 


In~\cite{2017arXiv170700026H} a practical method is presented for sampling from the optimal distribution $w$ in~\eqref{eq:optsample}. Essentially, one chooses a basis function $\phi_j$ from $\{\phi_j\}_{j=0}^n$ uniformly at random, and sample from a probability distribution proportional to $\phi_j^2$, and repeat this $N$ times. 
This strategy is very simple to implement, and we have adopted this in our experiments.



We note that MC with IS (sampling from $p=\frac{1}{w}$) does not reduce exactly to MCLS with $n=0$: the estimate by MC with IS is $\frac{1}{N}\sum_{i=1}^Nw(\tilde\x_i)f(\tilde\x_i)$, whereas that of MCLS is $\c=(\tilde\V^T\tilde\V)^{-1}\tilde\V^T\sqrt{\W}\tilde\f=(\tilde\V^T\tilde\V)^{-1}\sum_{i=1}^Nw(\tilde\x_i)f(\tilde\x_i)$. The difference is the factor $N(\tilde\V^T\tilde\V)^{-1}=N/\sum_{i=1}^Nw(\tilde\x_i)$, which is not 1, although it tends to 1 as $N\rightarrow \infty$. Indeed the two methods have different variances, as we show next. 


\ignore{
A significant number of recent studies focus on polynomial approximation via least-squares fitting. Among these, the most relevant are Cohen-Davenport-Leviatan~\cite{cohen2013stability}, Cohen-Migliorati~\cite{cohenoptimal} and 
Narayan,  Jakeman and Zhou \cite{narayan2017christoffel}. 
In the case of polynomial approximation, they show that
if least-squares fitting is done using 
\begin{enumerate}
\item orthogonal polynomials as a basis, 
\item sampling/weighting is employed proportional to the Chebyshev measure $1/\sqrt{1-x^2}$ (i.e., importace sampling with respect to Chebyshev measure)
\end{enumerate}
then the Vandermonde matrix $V\in\mathbb{R}^{m\times n}$ for $Vc\approx f$ is well-conditioned 
by taking $m\geq n\log n$, as opposed to the $O(n^2)$ samples required with equispaced sampling. In this way, we can increase the degree of polynomials adaptively as we sample more. 
To solve the least-squares problems, we apply 
a Krylov subspace method for the normal equation: 
this converges in a constant number of steps since the Vandermonde system is well-conditioned. 
To esimate the variance, 
}

\subsubsection{Variance}\label{optscale:adaptvariance}
We examine the variance and convergence of the estimator $\hat I_N=c_0$ in MCLS with the weighted sampling~\eqref{eq:MCLSweight}. 
\begin{theorem}\label{thm:lsmcoptvariance}
Fix $n$ and the basis functions $\{\phi_j\}_{j=0}^n$. 
Then with the weighted sampling $\sim \frac{1}{w}$, 
the MCLS estimator $\hat I_N$ in~\eqref{eq:lsmcval}
where $\c$ is the solution of~\eqref{eq:MCLSweight} satisfies 
 \begin{equation}
   \label{eq:ISvariance}
\sqrt{N}(\hat I_N-I)\xrightarrow{d} \mathcal{N}(0,    \min_\c\|\sqrt{w}(f-\sum_{j=0}^nc_j\phi_j)\|_2^2)   .
 \end{equation}  
\end{theorem}
{Proof.} 
We argue as in Theorem~\ref{eq:MCLSestimate}. Again we write 
$f = \sum_{j=0}^nc_j^\ast\phi_j+g=:f_1+g$. 
Then $\hat I_n-I = [1,0,\ldots,0](\tilde\V^T\tilde\V)^{-1}\tilde\V^T\tilde\g$, where $\tilde\g=[\sqrt{w(\tilde\x_1)}g(\tilde\x_1),\ldots,\sqrt{w(\tilde\x_N)}g(\tilde\x_N)]^T$. 
We have 
\begin{align*}
\frac{1}{N}(\tilde \V^T\tilde \V)_{i+1,j+1} &= \frac{1}{N}\sum_{\ell=1}^Nw(\tilde\x_\ell)\phi_i(\tilde\x_\ell)\phi_j(\tilde\x_\ell)\\
&\rightarrow 
\int_\Omega w(\tilde\x)\phi_i(\tilde\x)\phi_j(\tilde\x)\frac{d\tilde\x}{w(\tilde\x)}=\int_\Omega \phi_i(\tilde\x)\phi_j(\tilde\x)d\tilde\x=\delta_{ij}  
\end{align*}
as $N\rightarrow \infty$ by the law of large numbers, so $\frac{1}{N}(\tilde \V^T\tilde \V)\rightarrow \I_{n+1}$. 
Hence
\begin{align}
   \label{eq:cg0til}  
\sqrt{N}(\hat I_N-I)&=
\sqrt{N}[1,0,\ldots,0]^T
\frac{1}{N}(\frac{1}{N}\tilde\V^T\tilde\V)^{-1}\tilde\V^T\tilde\g\\
&\rightarrow \sqrt{N}\left(\frac{1}{N}\sum_{i=1}^Nw(\tilde\x_i)g(\tilde\x_i)\right)
\xrightarrow{d} \mathcal{N}(0,\|\sqrt{w}g\|^2_2)   
\end{align}
by the central limit theorem, where we used the fact $\int_\Omega g(\x)d\x=0$ for the mean and $\int_\Omega(w(\tilde\x)g(\tilde\x))^2\frac{d\tilde\x}{w(\tilde \x)}=\|\sqrt{w}g\|_2^2 $ for the variance. 
\hfill$\square$

Using the samples $\{\tilde \x_i\}_{i=1}^N$ and 
solution $\c=[c_0,\ldots,c_n]^T$ for~\eqref{eq:MCLSweight}, 
the variance $\|\sqrt{w}g\|_2$ above can be estimated via 
\begin{equation}  \label{eq:wg}
\|\sqrt{w}g\|_2^2\approx \frac{1}{N-n-1}\sum_{i=1}^N (w(\tilde\x_i))^2(f(\tilde\x_i)-
\sum_{j=0}^nc_j\phi_j(\tilde\x_j)). 
\end{equation}
 Note the power 2 in $(w(\tilde\x_i))^2$, due to the nonuniform sampling $\sim \frac{1}{w}$. 

\subsubsection{Relation to MC with importance sampling}

As briefly mentioned in~\cite{cohenoptimal},  
sampling from a different probability distribution 
 is a common technique in Monte Carlo known as \emph{importance sampling} (IS)~\cite[Ch~.9]{owenmcbook}. 
However, IS is fundamentally different from the above strategy: 
 The goal above is to improve the conditioning $\kappa_2(\tilde\V)$, not to reduce the MC variance as in IS (we derive the MCLS variance shortly).

The variance $\|\sqrt{w}g\|_2^2$ in~\eqref{eq:ISvariance} is different from the constant $\|g\|_2^2$ in 
MCLS~\eqref{eq:MCLSestimate}, and either can be smaller. 
This effect is similar to MC with IS, although the variance with IS is 
 $\|\frac{1}{\sqrt{w}}(wf-c_0^*)\|^2$~\cite[Ch.~9]{owenmcbook}, which is again different. 
Our viewpoint provides a fresh way of understanding IS: it samples from a probability distribution $p$ so that $f/p$ can be approximated well by a constant function. 

\subsection{MCLSA with polynomials: adaptively chosen degree}\label{optscale:adapt}
For a fixed degree $k$ (and hence $n$), 
the Vandermonde matrix $\tilde\V$ becomes well-conditioned if we sample at $N=O(n\log n)$ or more points according to the probability measure~\eqref{eq:optsample}, as we outlined above. 
On the other hand, generally speaking, the larger the degree $k$, the better the approximation quality $\min_{\c}\|f-\sum_{j=0}^nc_j\phi_j\|_2$ becomes, and therefore the smaller the variance; recall Table~\ref{tab:compare}. 

This motivates an MCLSA method where we increase the degree $k$ (and hence $n$) with $N$. 
Specifically, we choose $k$ to be the largest integer for which 
the resulting $\tilde\V$ is well-conditioned with high probability, that is, 
$n=_{d+k}C_d\lesssim N/\log N$. In the experiments below, we use the criterion $n=_{d+k}C_d\leq N/10$, which always gave $\kappa_2(\tilde\V)\leq 3$. 

We remark that while increasing $k$ is always beneficial for improving the approximation $\|f-p\|_2$ (assuming $\kappa_2(\tilde\V)=O(1)$), it comes with the obvious computational cost of (i) building the matrix $\V$, and (ii) solving the least-squares problem~\eqref{eq:ls}. (i) requires $O(Nnk)$ operations, and (ii) requires $O(Nn)$ operations using the conjugate gradient algorithm. In our setting this translates to $O(N^2)$ operations, which will at some $N$ exceed the sampling cost $C_fN$ ($C_f$ is the cost of sampling $f$). A practical strategy would be to set an upper bound on $k$ (hence $n$) so that $n\lesssim C_f$. 

While we focus on the use of polynomials with total degree, 
many other choices are possible. In particular, growing the basis adaptively 
 (depending on $f$) is 
a commonly employed technique~\cite{cohen2015approximation,gerstner2003dimension}, 
which is also a successful strategy in sparse grids~\cite{bungartz2004sparse}. 

\subsection{Convergence of MCLSA}\label{sec:converge}
Here we briefly discuss the convergence of MCLSA with polynomials. 
We show that when $f$ is analytic in $\Omega$ the convergence is superalgebraic, satisfying the last desideratum in Section~\ref{sec:desire}. 
\begin{proposition}\label{prop:analyconv}
Let $f:[0,1]^d\rightarrow \mathbb{R}$ be analytic in an open set containing $[0,1]^d$. Then 
MCLSA (in which \eqref{eq:MCLSweight} is solved with $n=O(N/\log N)$ and $\{\phi_j\}_{j=0}^n$ is a basis for polynomials of bounded total degree)
has error 
$\E[|\hat I_N-I|]=O(\frac{\exp(-cN^{-1/d}/\sqrt{d})}{\sqrt{N}})$ for some constant $c>0$.
\end{proposition}
{\sc proof.}
The analyticity assumption implies 
 that $f$ can be approximated by a $d$-dimensional 
 polynomial of total 
 degree $k$ as~\cite{trefethen2017multivariate}
\begin{equation}  \label{eq:f-p}
\inf_{\mbox{deg}(p)\leq k}\|f-p\|_{[0,1]^d} = O(\exp(-ck/\sqrt{d})),
\end{equation}
where 
$\|\cdot\|_{[0,1]^d}$ denotes the supremum norm and 
$c$ is a constant depending on the location of the singularity (if any) of $f$ nearest $[0,1]^d$ with respect to the radius of the Bernstein ellipse. 

\ignore{
Then with $N$ sample points, we can take $k\approx N^{-1/d}$, so 
the approximation error decays like $\exp(-cN^{-1/d})$; resulting in 
the MCLS convergence $O(\exp(-cN^{-1/d})/\sqrt{N})$.
We note that asymptotically as $N\rightarrow \infty$, the $\exp(-cN^{-1/d})$ convergence is faster than any algebraic convergence $N^{-c}$, suggesting superalgebraic convergence (as desired in Section~\ref{sec:desire});
however, clearly the $-1/d$ power in $N^{-1/d}$ means an extremely large $N$ is necessary to obtain good convergence. 


The preceding argument assumes the use of maximum degree; in practice a more common choice is the total degree, as we just did above. 

In this case, } 

The number of basis functions for $d$-dimensional polynomials of degree $k$ 
is $\big(
\begin{smallmatrix}
d+k\\k  
\end{smallmatrix}
\big)$, which 
is $N=O(k^d)$ in the limit $N\rightarrow \infty$ (hence $k\gg d$), so 
$k=O(N^{-1/d})$. 
Together with~\eqref{eq:f-p}, we conclude that with $N$ points, the approximation quality $\|f-p\|_{[0,1]^d}$ is 
\begin{equation}
  \label{eq:inffpdegpk}
\inf_{\mbox{deg}(p)\leq k}\|f-p\|_{[0,1]^d} = O(\exp(-cN^{-1/d}/\sqrt{d})).   
\end{equation}
From this it immediately follows that 
$\|f-p\|_2 = O(\exp(-cN^{-1/d}/\sqrt{d}))$, and the result follows from 
the fact that the MCLS convergence is $O(\frac{\|f-p\|_2}{\sqrt{N}})$. 
\hfill $\square$ 


We note that when $d\gg 1$, 
$N=O(k^d)$ with $k\gg d$ means $N$ is astronomically large. 
Thus it may be more reasonable to consider 
the regime $d\gg k$, 
in which case $N=O(d^k)$, 
hence $k=\log N^{1/d}$. Therefore the function approximation quality is $O(\exp(-ck/\sqrt{d}))=O(\exp(-c\log N^{1/d}/\sqrt{d}))=O(N^{-c/d\sqrt{d}})$, suggesting an algebraic convergence. 

We also note that Trefethen~\cite{trefethen2017multivariate} argues that for a certain natural class of functions analytic in the hypercube (with singularities outside), neither the total nor the maximum degree would be the optimal choice. 
In particular, employing the \emph{Euclidean} degree minimizes the approximation error for a fixed $N$. 
In terms of the function approximation quality~\eqref{eq:f-p}, the use of Euclidean degree removes the $1/\sqrt{d}$ factor in the exponent in~\eqref{eq:inffpdegpk}. 
However, our experiments suggest that with MCLSA, the total degree and Euclidean degree perform almost equally well for the class of functions considered in~\cite{trefethen2017multivariate}; we suspect that the $1/\sqrt{N}$ term in the convergence $\|f-p\|_2/\sqrt{N}$ is playing a significant role. For less smooth functions, the total degree appears to have better accuracy. 
They both perform significantly better than the maximum degree. For these reasons, in our experiments below we employ the total degree unless otherwise specified. 

\subsubsection{Exact degree}\label{sec:exact}
Recalling the discussion in Section~\ref{sec:careful}, MCLSA 
integrates $f$ exactly if it
is a polynomial of degree at most $k$. 
Thus as long as 
$\V$ has full column rank, 
we obtain $p=f$, and hence the integral is exact. 
Such (randomized) integration formulae that provide exact results for polynomials of bounded degree were investigated by Haber~\cite{haber1968combination,haber1969stochastic}, and are called stochastic quadrature formulae. 

In the design of algorithms for cubature~\cite{cools1997constructing}, 
significant effort is made to maximize the polynomial degree of exactness for a specified number of (deterministic) sample points. 
In such algorithms, the exactness holds for a linear space whose dimension is higher than the number of sample points (Gauss quadrature is a typical example with double-degree exactness). 
By contrast, in MCLS the sample points are taken randomly, and so the estimate is exact only for $f\in\mbox{span}\{\phi_j\}_{j=0}^{n}$. 
Therefore optimality in terms of the polynomial degree exactness 
is not attained. The gain in MCLS (relative to cubature) is rather in the $O(1/\sqrt{N})$ MC convergence. 




\subsection{Variance estimate when $\kappa_2(\tilde\V)\not\approx 1$}\label{sec:kappa}

In MCLSA we allow $n$ to grow with $N$, so that $n=O(N)$ up to logarithmic factors. 
In this case we cannot use the limit 
$\frac{1}{N}(\tilde \V^T\tilde \V)\rightarrow \I_{n+1}$ as 
$N\rightarrow\infty$, which requires $n$ to be fixed as was in Theorem~\ref{thm:lsmcoptvariance}. Hence the variance estimate via~\eqref{eq:ISvariance} is not directly applicable. 

Here we 
examine the variance of the MCLS estimator $c_0$ in~\eqref{eq:MCLSweight} 
when $n$ is not negligible compared with $N$, and therefore $\kappa_2(\tilde\V)$ cannot be taken to be 1. Note that here $\kappa_2(\tilde\V)$ does depend on the specific choice of basis $\{\phi_j\}_{j=0}^n$ and not just its span. 
To exclude ill-conditioning caused by a poor choice of $\{\phi_j\}_{j=0}^n$ (rather than poor sample points $\{\x_i\}_{i=1}^N$ or too small $N$), 
we continue to assume that the basis functions are orthonormal in the continuous setting $\int_\Omega \phi_i(\x)\phi_j(\x)d\x = \delta_{ij}$. 

Our goal is to estimate $\EE[\Delta c_0^2\big|\kappa_2(\tilde\V)\leq K]$, where 
$\Delta c_0=c_0-c_0^*$ is the error in $c_0$. 
The expectation $\EE$ is taken over the random samples $\{\tilde\x_i\}_{i=1}^N\sim \frac{1}{w}$, conditional on $\kappa_2(\tilde\V)\leq K$ for some $K>1$. In practice, we will take $K$ to be the exact or estimated value of the particular $\kappa_2(\tilde\V)$. 
Better yet would be to condition on $\kappa_2(\tilde\V)= K$ and estimate $\EE[\Delta c_0^2\big|\kappa_2(\tilde\V)= K]$, but analyzing this appears to be difficult. The next result bounds the variance of the whole vector $\Delta \c=\c-\c^*$. 

\begin{proposition}In the above notation, 
\begin{equation}  \label{eq:cnorm}
\EE[\|\Delta \c\|_2^2\big|\kappa_2(\tilde\V)\leq K]\leq  
\frac{n+1}{N}\frac{K^4}{1-\epsilon_K}
\min_\c\|f-\sum_{j=0}^nc_j\phi_j\|_2^2, 
\end{equation}  
where $\epsilon_K=2(n+1)\exp(-\frac{c_{\delta_K} N}{n+1})$, 
$c_{\delta_K}=\delta_K+(1-\delta_K)\ln (1-\delta_K))>0$,
in which $\delta_K$ is defined via 
$K= \sqrt{\frac{1+\delta_K}{1-\delta_K}}$. 
\end{proposition}
{\sc proof}.
We first investigate $\P[\kappa_2(\tilde\V)\leq K]$ and show that the probability 
is at least $1-\epsilon_K$ (which is $\approx 1$ when $K$ is moderately larger than 1 and $N\geq n\log n$). 
Noting that 
$\frac{1}{N}\tilde\V^T\tilde\V
=\frac{1}{N}\sum_{i=1}^Nw(\tilde\x_i)[1,\phi_1(\tilde\x_i),\ldots,\phi_n(\tilde\x_i)]^T
[1,\phi_1(\tilde\x_i),\ldots,\phi_n(\tilde\x_i)] = :\frac{1}{N}\sum_{i=1}^N\tilde\X_i$
 can be regarded as a sum of $N$ independent positive semidefinite matrices $\tilde\X_i\succeq 0$, the matrix Chernoff inequality~\cite{tropp2012user} implies that for $0<\delta<1$, 
\begin{equation}  \label{eq:Pchernoff}
\P[\|\tilde\V^T\tilde\V-\I_{n+1}\|_2>\delta]\leq 2(n+1)\exp(-c_\delta/R),   
\end{equation}
where $c_\delta = \delta+(1-\delta)\ln (1-\delta))>0$ and 
$R$ is a bound such that $R\geq \|\X_i\|_2$  almost surely; here we can take 
$R = \frac{n+1}{N}$, by~\eqref{eq:optsample}. When $N\approx n\log n$, the right-hand side of~\eqref{eq:Pchernoff} decays algebraically with $N$. Note further that the complement 
$\P[\|\tilde\V^T\tilde\V-\I_{n+1}\|_2\leq \delta]$
implies $\sqrt{1-\delta}\leq \sigma_i(\tilde\V)\leq \sqrt{1+\delta}$, so 
$\kappa_2(\tilde\V)\leq \sqrt{\frac{1+\delta}{1-\delta}}$. Thus 
defining $\delta_K$ via $K= \sqrt{\frac{1+\delta_K}{1-\delta_K}}$ and 
setting $c_{\delta_K}=\delta_K+(1-\delta_K)\ln (1-\delta_K))>0$, 
from~\eqref{eq:Pchernoff} we obtain 
\begin{equation}  \label{eq:PchernoffK}
\P[\kappa_2(\tilde\V)\leq K]\geq 1-2(n+1)\exp(-c_{\delta_K} N/n)=:
1-\epsilon_K. 
\end{equation}

\ignore{
In view of~\eqref{eq:MCLSweight}, 
define $\tilde \V=\W^{1/2}\V$, $\tilde \f=\W^{1/2}\f$, and 
$\tilde \g=[\sqrt{w(\x_1)}g(\x_1),\ldots,\sqrt{w(\x_N)}g(\x_N)]^T$. 
The solution $\c$ in \eqref{eq:MCLSweight} is 
$\c=(\tilde\V^T\tilde\V)^{-1}\tilde\V^T\tilde\f$, and the error in $\c$ 
is 
\begin{equation}  
\label{eq:MCLSerrcohen}
\Delta \c=(\tilde\V^T\tilde\V)^{-1}\tilde\V^T\tilde \g. 
 \end{equation}
} 

We now turn to bounding $\Delta\c$. 
First decompose $f=f_1+g$ as in the proof of Theorem~\ref{thm:mainvar}
(the nonuniform sampling has no effect on the decomposition). 
We have $\Delta \c=(\tilde\V^T\tilde\V)^{-1}\tilde\V^T\tilde \g$, and 
so
$\|\Delta \c\|_2\leq \frac{1}{(\sigma_{\min}(\tilde\V))^2}\|\tilde\V^T\tilde \g\|_2$. 
Now since $\tilde\V$ has rows of norm $\sqrt{n+1}$, it follows that 
$\|\tilde\V\|_F=\sqrt{N(n+1)}$, and hence 
$\|\tilde\V\|_2\geq \sqrt{N}$ (this bound is tight when $\tilde\V$ is close to having orthonormal columns). 
We thus obtain 
$\frac{1}{\sigma_{\min}(\tilde\V)}\leq \frac{\kappa_{2}(\tilde\V)}{\sqrt{N}}$, so 
$\|\Delta \c\|_2\leq  \frac{1}{N}(\kappa_2(\tilde\V))^2\|\tilde\V^T\tilde \g\|$. 
It follows that 
\begin{equation}  \label{eq:EEc2}
\EE[\|\Delta \c\|_2^2\big|\kappa_2(\tilde\V)\leq K]\leq  \frac{1}{N^2}K^4
\EE[\|\tilde\V^T\tilde \g\|_2^2\big|\kappa_2(\tilde\V)\leq K]. 
\end{equation}
To bound the right-hand side, we next examine $\EE[\frac{1}{N}\|\tilde\V^T\tilde \g\|_2^2\big|\kappa_2(\tilde\V)\leq K]$. 
Using~\eqref{eq:PchernoffK} and Markov's inequality we obtain 
\begin{equation}  \label{eq:markov}
\EE[\frac{1}{N}\|\tilde\V^T\tilde \g\|_2^2\big|\kappa_2(\tilde\V)\leq K]
\leq \frac{\EE[\frac{1}{N}\|\tilde\V^T\tilde \g\|_2^2]}{1-\epsilon_K}  .
\end{equation}
The remaining task is to bound $\EE[\|\tilde\V^T\tilde \g\|_2^2]$. 
For this, the $k$th element 
$\EE[\frac{1}{N}(\tilde\V^T\tilde \g)_k^2]$
has been worked out in~\cite[\S~3]{cohenoptimal}, which we repeat here:
\begin{align*}
\mathbb{E}_{\frac{1}{w}}[\frac{1}{N}(\tilde\V^T\tilde \g)_k^2]&=
\mathbb{E}_{\frac{1}{w}}[
\frac{1}{N}\sum_{i=1}^N\sum_{j=1}^Nw(\tilde\x_i)w(\tilde\x_j)g(\tilde\x_i)\phi_k(\tilde\x_i)g(\tilde\x_j)\phi_k(\tilde\x_j)]  \\
&=
\frac{1}{N}\mathbb{E}_{\frac{1}{w}}[\sum_{i=1}^N\left(w(\tilde\x_i)g(\tilde\x_i)\phi_k(\tilde\x_i)\right)^2]=\int_\Omega (w(\tilde\x)g(\tilde\x)\phi(\tilde\x))^2\frac{1}{w(\tilde\x)}d\tilde\x\\
&=\int_\Omega w(\x)(g(\x)\phi(\x))^2d\x. 
\end{align*}
The $i\neq j$ terms disappear in the second expression because their expectation are a product of two terms both equal to
$\int_\Omega w(\tilde\x)g(\tilde\x)\phi(\tilde\x)\frac{1}{w(\tilde\x)}d\tilde\x=
\int_\Omega g(\x)\phi(\x)d\x=0$. 
Therefore 
\begin{equation}  \label{eq:cvariance}
\mathbb{E}_{\frac{1}{w}}[\frac{1}{N}\|\tilde\V^T\tilde \g\|_2^2]
=\int_\Omega w(\x)\left(\sum_{k=0}^n(\phi_k(\x))^2(g(\x))^2\right)d\x
=(n+1)\int_\Omega (g(\x))^2d\x=(n+1)\|g\|_2^2, 
\end{equation}
where we used the definition~\eqref{eq:optsample} of $w$. 
Putting together~\eqref{eq:EEc2}, \eqref{eq:markov} and~\eqref{eq:cvariance} 
completes the proof. 
\hfill$\square$

We now turn to estimating $\EE[|\hat I_N-I|^2]=\EE[\Delta c_0^2]$. 
The left-hand side in~\eqref{eq:cnorm} is 
$\EE[\Delta c_0^2]+\EE[\Delta c_1^2]+\cdots+\EE[\Delta c_n^2]$. This is a sum of $n+1$ terms, and 
\ignore{
Thus 
$\frac{n+1}{N}\|g\|_2^2=\mathbb{E}_{\frac{1}{w}}[\frac{1}{N}(\tilde\V^T\tilde \g)^2]=\sum_{k=0}^{n}\mathbb{E}_{\frac{1}{w}}[\frac{1}{N}(\tilde\V^T\tilde \g)_k^2]
$, which is a sum of $n+1$ positive terms. 
~\eqref{eq:markov}
}
we expect no term to dominate the others, suggesting $\EE[\Delta c_0^2\big|\kappa_2(\tilde\V)\leq K]\lessapprox K^4\frac{\|g\|_2^2}{N}$. 
Together with the fact 
$N\EE [\Delta c_0^2\big|\kappa_2(\tilde\V)\leq K]\rightarrow \|\sqrt{w}g\|^2_2$ as $N\rightarrow \infty$ 
(and hence $\kappa_2(\tilde \V)\rightarrow 1$)
as shown in Theorem~\ref{thm:lsmcoptvariance}, this suggests the estimate 
$\EE [\Delta c_0^2]\lessapprox \frac{1}{N}K^4\|\sqrt{w}g\|^2_2$, leading to the error estimate for $\hat I_N-I=\Delta c_0$
\begin{equation}
  \label{eq:c0estimate}
\EE [|\Delta c_0|]\lessapprox \frac{(\kappa_2(\tilde\V))^2}{\sqrt{N}}\|\sqrt{w}g\|_2.
\end{equation}
 The discussion in this paragraph is clearly heuristic; 
rigorously and tightly bounding $\EE [\Delta c_0^2]$ is left as an open problem. 

When the standard sampling strategy is employed (i.e., $w=1$), essentially the whole argument carries over: for example, $\|\V\|_2\geq \sqrt{N}$ continues to hold since the first column of $\V$ is a $N$-vector of 1's. 
However, the value of $R$ becomes larger and $\epsilon_K$ decreases slower as $N$ grows, with $\epsilon_K\approx 1$ if $N\ll n^2$.

\begin{figure}[htbp]
  \centering
      \includegraphics[width=0.45\textwidth]{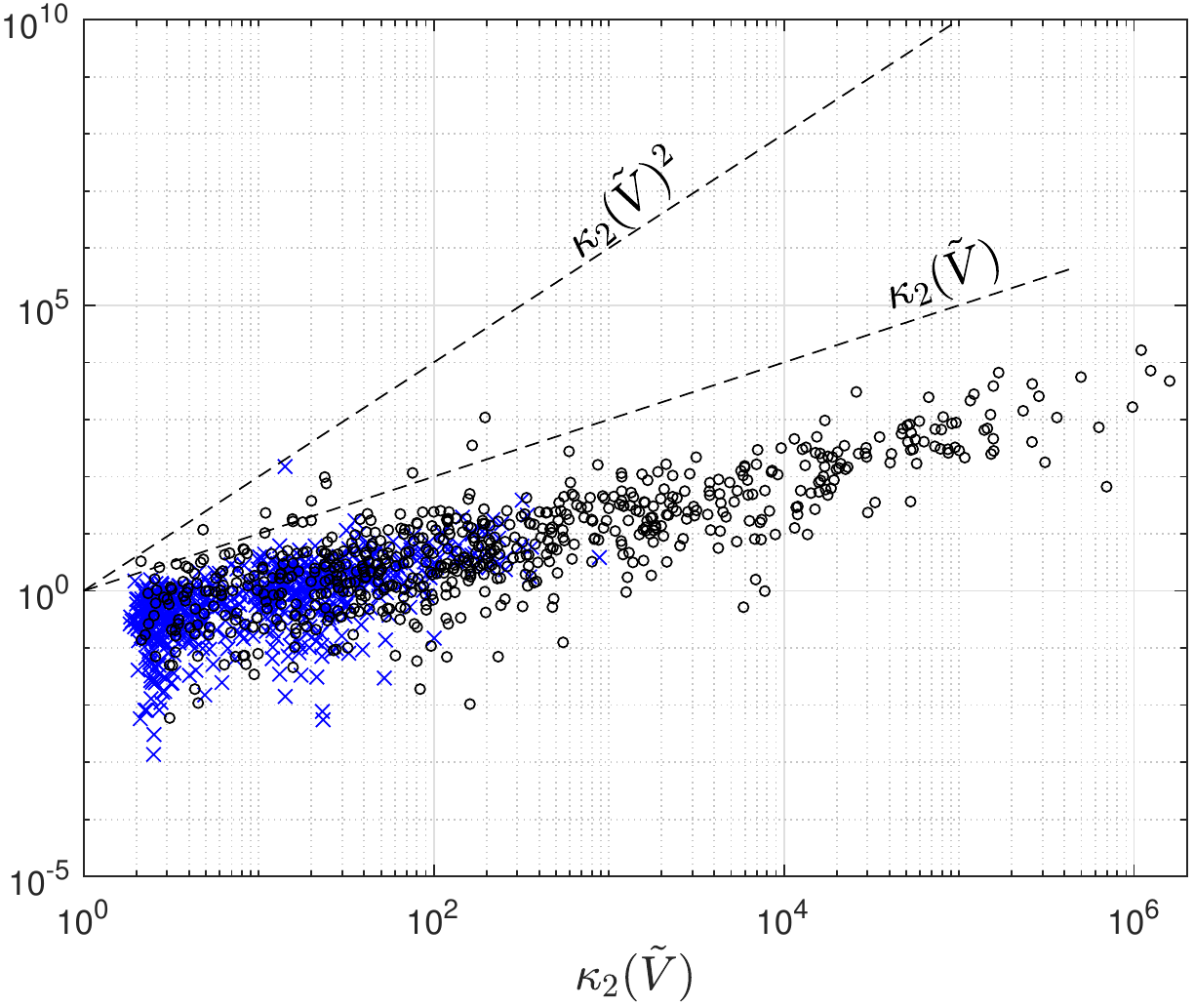}
  \caption{Scatterplots of $|c_0-\bar{c_{0}}|/\frac{\|\sqrt{w}\hat f\|_2}{\sqrt{N}}$ (error divided by conventional estimate~\eqref{eq:ISvariance}) against $\kappa_2(\tilde \V)$, with MCLSA (blue cross) and MCLS (black cirles). 
}
  \label{fig:condVdependence}
\end{figure}

To illustrate the estimate~\eqref{eq:c0estimate}, 
we perform the following experiment:
consider computing $\int_{[0,1]^3}x_1^{10}x_2^5x_3^7 dx$.
We take $320$ instances of MCLS and MCLSA, using polynomials of total degree between $5$ and $20$, and varying the aspect ratio $N/n$ in $\{1.1,1.2,\ldots,2\}\cup \{3,4,\ldots,10\}$ (the larger the aspect ratio, the better-conditioned $\tilde \V$ tends to be). 
Figure~\ref{fig:condVdependence} illustrates the result, which is a scatterplot of the observed values $\kappa_2(\tilde \V)$ versus $|c_0-\bar{c_{0}}|/\frac{\|\sqrt{w}g\|_2}{\sqrt{N}}$, the error 
divided by the ``conventional'' convergence analysis~\eqref{eq:ISvariance} 
where $\kappa_2(\tilde \V)\rightarrow 1$ is assumed. 
Observe that $\kappa_2(\tilde\V)\ll \kappa_2(\V)$ 
and in particular $\kappa_2(\tilde\V)=O(1)$ in most cases (even with the small aspect ratio 1.1), reflecting the optimal sampling used in MCLSA. 


While the estimate~\eqref{eq:c0estimate} suggests the error would scale like $\kappa_2(\tilde \V)^2$, our experiments indicate linear (or less) dependence, for both MCLS and MCLSA. 
We have also verified that the errors are smaller than 
$2\kappa_2(\tilde \V)\frac{\|\sqrt{w}g\|_2}{\sqrt{N}}$ in at least $95\%$ of the instances (the plot looks much the same under different settings e.g. different $f$ or $d$). 
In view of these, in what follows we estimate the confidence intervals of MCLSA via 
\begin{equation}  \label{eq:useCI}
\hat I_N\pm 2\kappa_2(\tilde \V)\frac{\|\sqrt{w}g\|_2}{\sqrt{N}}  . 
\end{equation}
In practice,  $\|\sqrt{w}g\|_2$ is estimated by~\eqref{eq:wg} and 
$\kappa_2(\tilde \V)$ can be estimated by a condition number estimator (e.g. MATLAB's {\tt condest}; here we computed it to full precision using {\tt cond}). 
We have also run MCLSA multiple times and confirmed that the actual errors lay within the confidence interval at least $95\%$ of the time. 
Note that~\eqref{eq:useCI} reduces to the conventional confidence interval in the standard Monte Carlo setting where $\kappa_2(\tilde \V)=1$ and $w=1$. 

The crux of this subsection  is that it is important to ensure 
$\tilde\V$
 is well conditioned in MCLS(A), and thus so is the use of the optimal sampling in Section~\ref{sec:leghigh}. 

Let us remark that in numerical analysis (in particular numerical linear algebra), matrix ill-conditioning is often regarded as a problem that manifests itself due to roundoff errors in finite precision arithmetic (i.e., ill-conditioning magnifies the effect of roundoff errors): in particular, 
issues associated with ill-conditioning usually disappear if high-precision arithmetic is used. 
Here ill-conditioning is playing a different (one might say more fundamental) role, and its effect persists even if one used exact arithmetic for the least-squares problem. 


\subsection{Numerical examples}\label{optscale:ex}
We present examples to compare 
 Monte Carlo (MC), MCLS with degree $k=5$ (MCLS-deg5), and MCLSA 
(the $j$th data point uses a degree $j$ polynomial in the figures). 
We take the dimension $d=6$, and integrate three functions of varying smoothness: 
\begin{itemize}
\item $f(\x) = \sin(\sum_{i=1}^d x_i)$. An analytic function; an example of Genz' first problem~\cite{bigoni2016spectral,genz1984testing}.
\item $f(\x) = \sum_{i=1}^d \exp(-|x_i-1/2|)$. Genz' fifth problem; function with  singularity of absolute value type, aligned with the grids.
\item A basket option arising in finance~\cite{glasserman2013monte}
  \begin{equation}
    \label{eq:basketf}
f(\x) = \exp(-rT)\max(0,S_0\exp((r-\frac{1}{2}\sigma^2)T+\sigma\sqrt{T}{\bf L}\Phi^{-1}(\x))-K).
  \end{equation}
Here we took $r=0.05,T=1,\sigma=0.2,K=S_0=10$, and ${\bf L}$ is the covariance matrix with $0.1$ in all the off-diagonals. $\Phi^{-1}$ is the inverse map of the cumulative distribution function of the multivariate normal distribution~(e.g. \cite[\S~2.3]{glasserman2013monte}). This function also has singularities of $|x|$ type, but now the singularity is not aligned with the grid. 
\end{itemize}

\begin{figure}[htbp]
  \begin{minipage}[t]{0.325\hsize}
  \centering
      \includegraphics[width=1\textwidth]{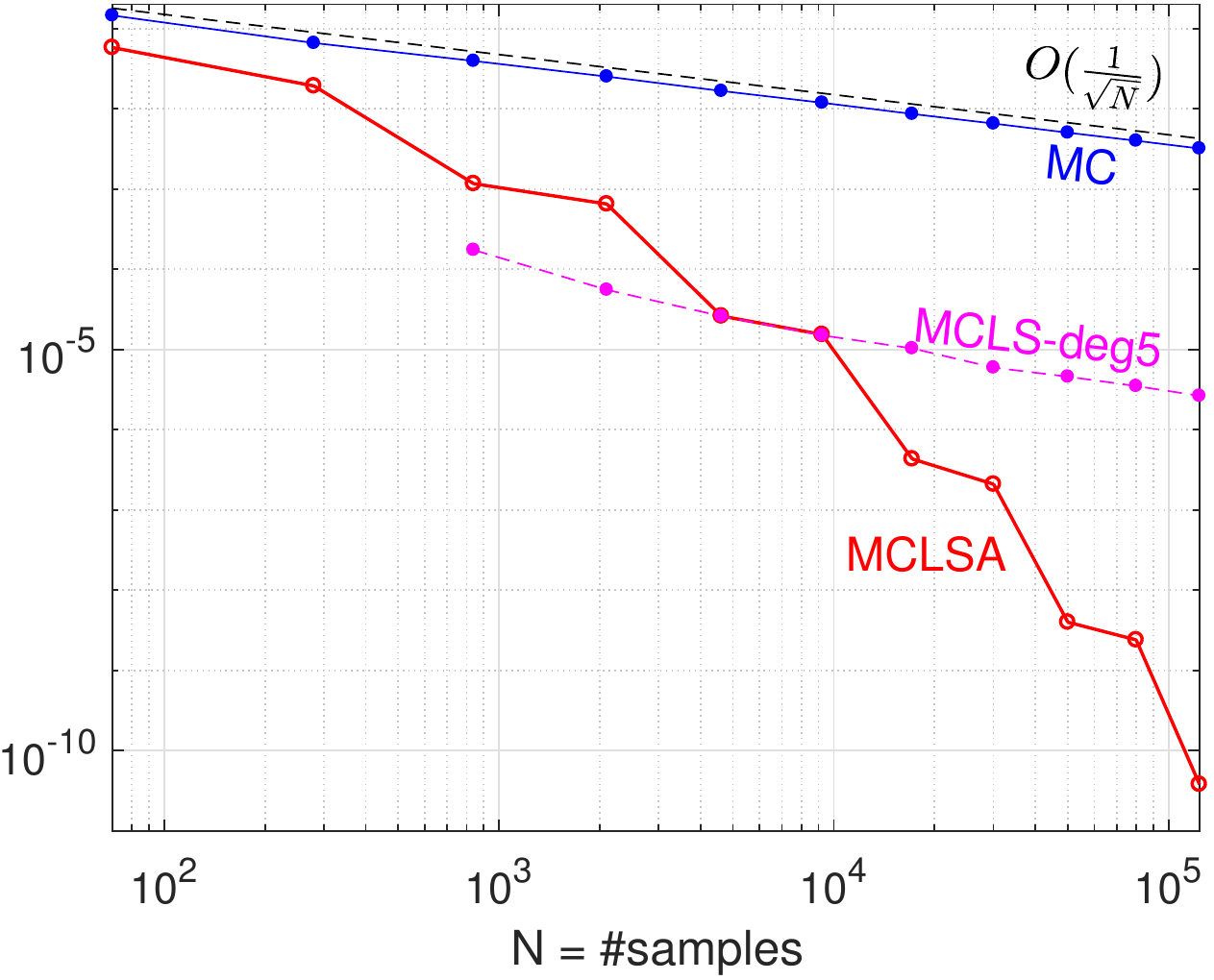}
  \end{minipage}
  \begin{minipage}[t]{0.33\hsize}
  \centering
      \includegraphics[width=.97\textwidth]{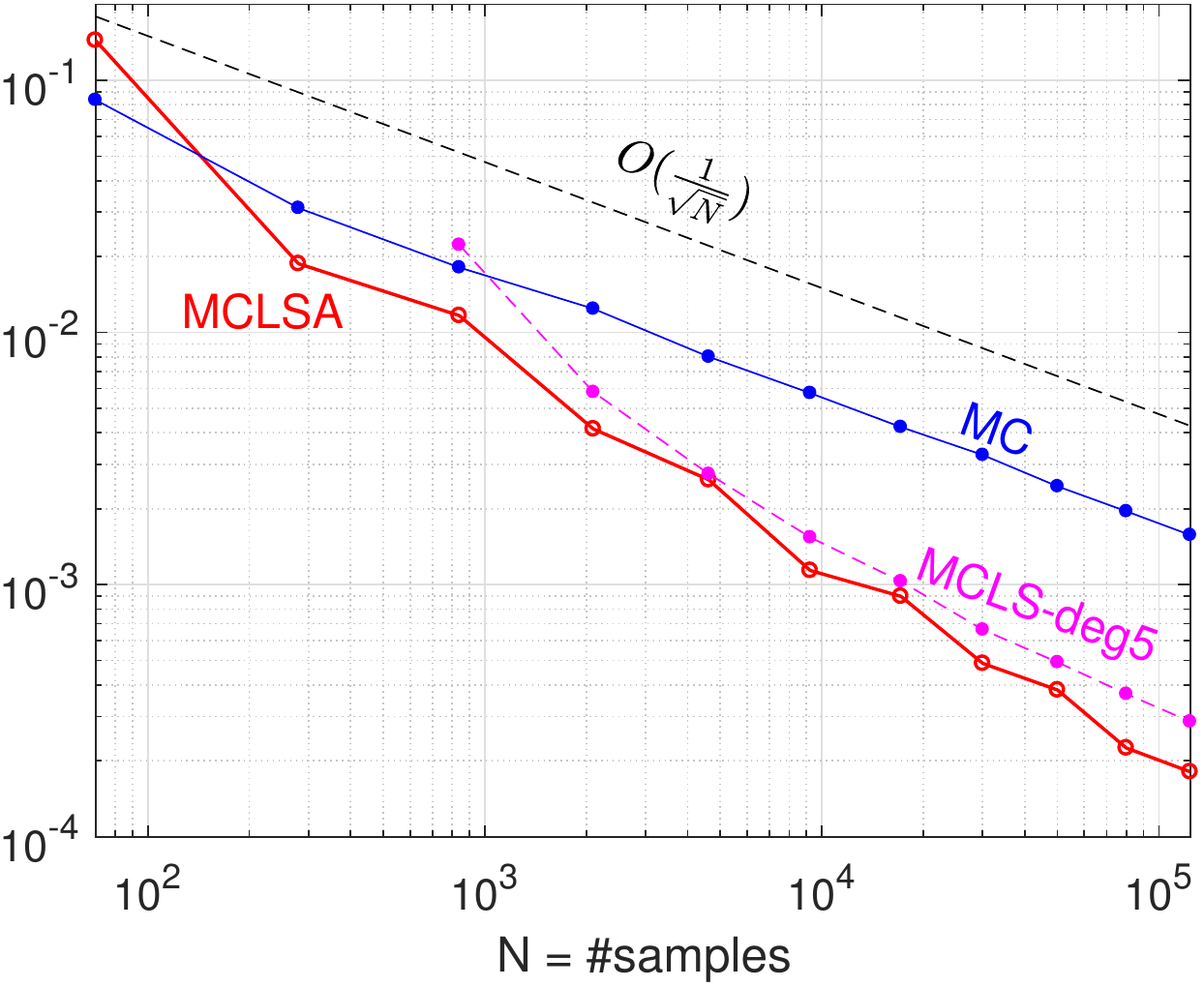}
  \end{minipage}
  \begin{minipage}[t]{0.325\hsize}
  \centering
      \includegraphics[width=.985\textwidth]{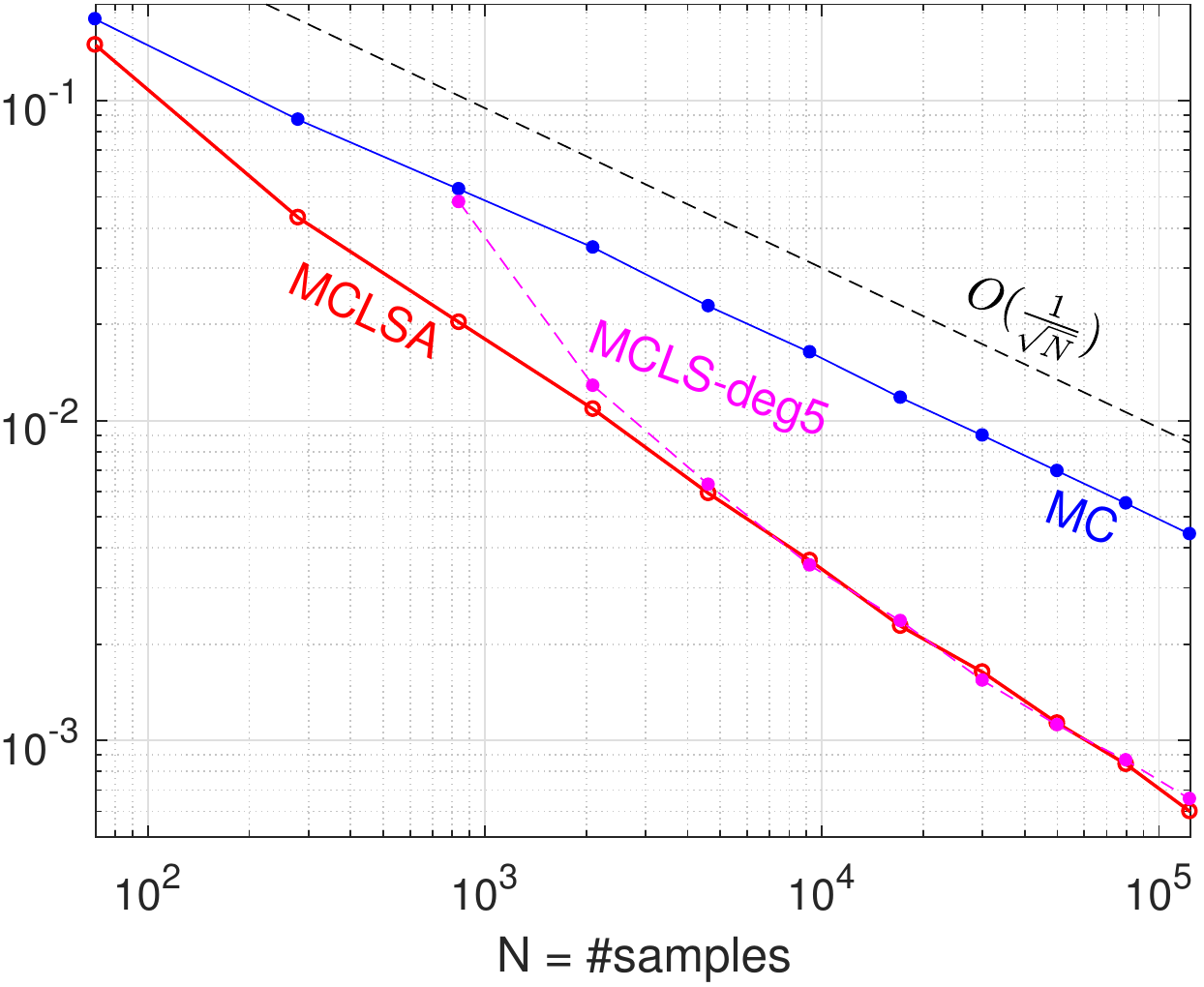}
  \end{minipage}
  \caption{Convergence (confidence interval widths) of MCLS with polynomials 
(fixed degree and adaptive degree with optimal sampling) for integration over $[0,1]^6$. Left: $f(\x) = \sin(\sum_{i=1}^d x_i)$ (smooth/analytic), center: $f(\x) = \sum_{i=1}^d \exp(-|x_i-1/2|)$ (absolute value singularity), right: basket option. 
}
  \label{fig:Legpolyinterp}
\end{figure}

The results are shown in Figure~\ref{fig:Legpolyinterp}. 
As before,  MCLS has the same $O(1/\sqrt{N})$ asymptotic convergence as MC (the apparently faster convergence in the pre-asymptotic stage is due to $\kappa_2(\tilde\V)$ decreasing, as in Figure~\ref{fig:1d}), with 
 reduced variance due to the improved underlying approximation, as indicated in Table~\ref{tab:compare}. Moreover, MCLSA exhibits dramatically improved convergence, combining the superalgebraic convergence of the approximant as in Proposition~\ref{prop:analyconv}
and (to a lesser extent) the $O(1/\sqrt{N})$ convergence of MC. 
There is no data for MCLS-deg5 with small values of $N$ because we require $N>n$ in MCLS. 

With a standard uniform (non-optimal) random sampling (not shown), the asymptotics of MCLS-deg5 look much the same but the initial data points will be significantly higher, caused by $\kappa_2(\V)\gg 1$. More importantly, the degree in MCLSA will need to be much lower for the same number $N$ (hence larger error) to ensure $\kappa_2(\V)=O(1)$. 

\section{MCLSA with sparse grids}\label{sec:spgridmc}
In this section we combine MCLSA with another 
popular class of technique in function approximation: sparse grids~\cite{bungartz2004sparse,smolyak1963quadrature}. 
Specifically, we simply take a single basis function $\phi_1=p_s$, where $p_s$ is a sparse grid approximant to $f$, of varying levels. Namely, the algorithm proceeds as follows. 
  \begin{enumerate}
  \item Obtain $p_s\approx f$ using sparse grids (using $N_s$ samples, depending on level $s$ and underlying 1-dimensional quadrature rule)
  \item Obtain random samples $\{\x_i\}_{i=1}^N\in[0,1]^d$, and 
solve the least-squares problem
    \begin{equation}
      \label{eq:spgridMCLS}
\min_{c\in\mathbb{R}^2}
{\left\|
  \begin{bmatrix}
1&p_s(\x_1)\\
1&p_s(\x_2)\\
\vdots&\vdots\\
1&p_s(\x_N)\\
  \end{bmatrix}
  \begin{bmatrix}
c_0\\c_1
  \end{bmatrix}
-  \begin{bmatrix}
    f(\x_1)\\f(\x_2)\\\vdots\\ f(\x_N)
  \end{bmatrix}
\right\|_2}      ,
    \end{equation}
Finally, take   $\hat I_{N_s+N} = c_0+c_1\int_{[0,1]^d} p_s(\x)d\x$. 
  \end{enumerate}

Here the number of sample points is the sum $N_s+N$, where $N_s$ (deterministic) samples are used to obtain the sparse grid approximant, and $N$ random samples are taken in the least-squares problem~\eqref{eq:spgridMCLS}. Note that the samples on the sparse grid are not included in \eqref{eq:spgridMCLS}, because on the sparse grid we have $p_s=f$ by construction. 

This integrator is another instance of MCLSA, where (unlike in Section~\ref{optscale:adapt}, where $n$ grows with $N$) $n$ is fixed to 1 but $\phi_1=p_s$ is chosen adaptively. Namely, given a computational budget $N_{total}$, one would need to choose $s$, which then determines $N_s$ and $N=N_{total}-N_s$. 

\subsection{Convergence}\label{ex:spconv}
By  Theorem~\ref{thm:mainvar}, the MCLS error estimate is $\frac{\min_\c\|f-\sum_{j=0}^nc_j\phi_j\|_2}{\sqrt{N}}$, which with~\eqref{eq:spgridMCLS} 
is $\frac{\|f-p_s\|_2}{\sqrt{N}}$. 
The convergence of $\|f-p_s\|_2$ is a central subject in sparse grids. 
Roughly speaking, for $k(\geq 2)$-times differentiable functions, we have 
$\|f-p_s\|_2=O(N_s^{-2})$ with a piecewise linear approximation and 
$\|f-p_s\|_2=O(N_s^{-k-1})$ with global polynomial interpolation (e.g. at Chebyshev points), up to factors $O(\log N_s^{kd})$, which often plays a significant role. 
 We refer to~\cite[Lem.~3.13, Lem.~4.9]{bungartz2004sparse} for details. 
The convergence of MCLS is thus $O(N_s^{-2}N^{-1/2})$ and $O(N_s^{-k-1}N^{-1/2})$ respectively, again up to factors $O(\log N_s^{kd})$. 

\subsection{Numerical examples}\label{ex:spgrid}
\begin{figure}[htbp]
  \begin{minipage}[t]{0.5\hsize}
      \includegraphics[width=0.95\textwidth]{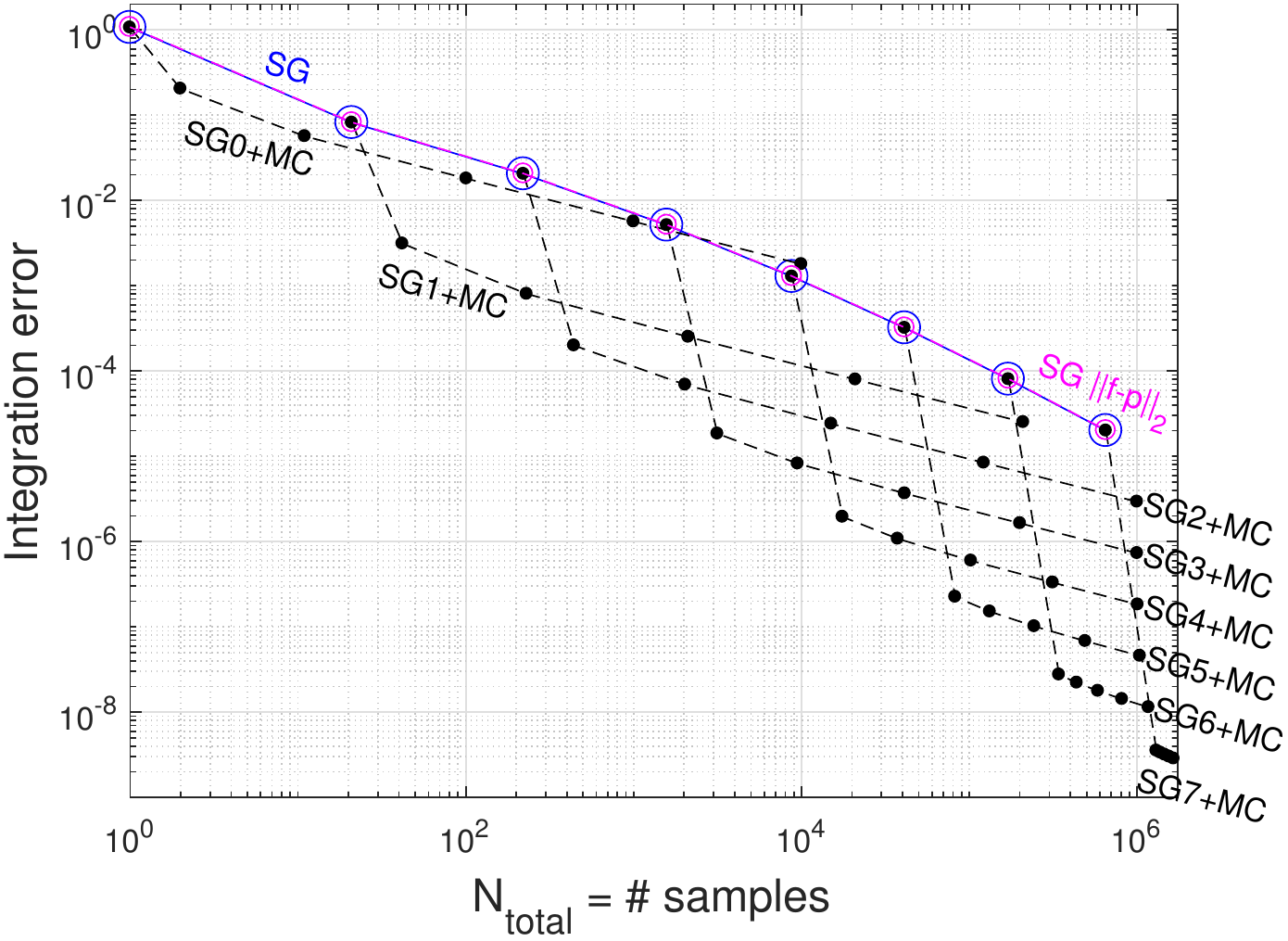}
  \end{minipage}   
  \begin{minipage}[t]{0.5\hsize}
      \includegraphics[width=0.95\textwidth]{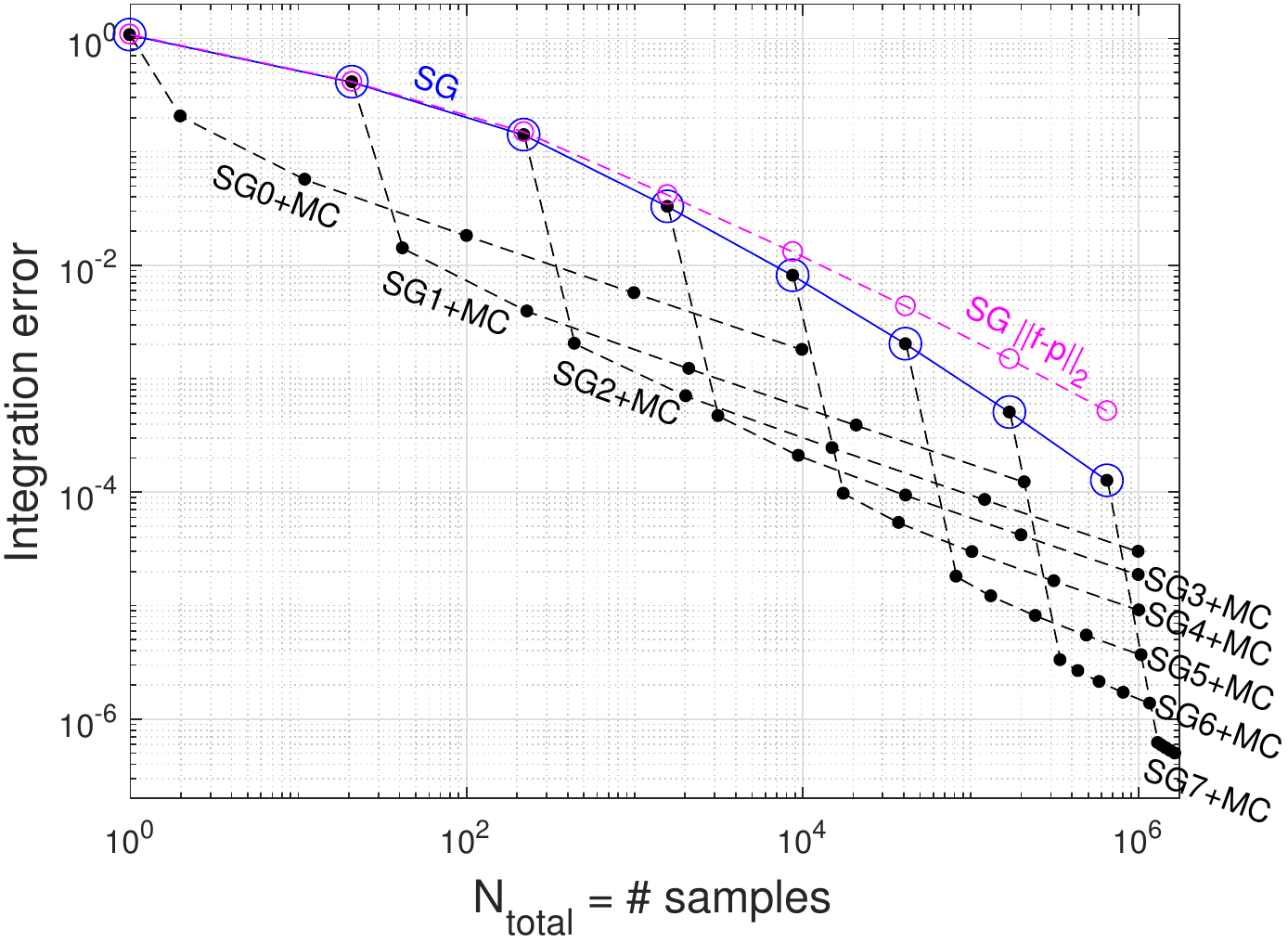}

  \end{minipage}
  \caption{
Integration error of SG quadrature and SG combined with MCLS (for which the $95\%$ confidence intervals are shown). $d=10$, $f(\x) = \sum_{i=1}^d \exp(-|x_i-1/2|)$. 
Left: SG using piecewise linear approximation. 
Right: SG 
using Clenshaw-Curtis. 
}
  \label{fig:spgridplusMC}
\end{figure}

\begin{figure}[htbp]
  \begin{minipage}[t]{0.5\hsize}
     \includegraphics[width=0.95\textwidth]{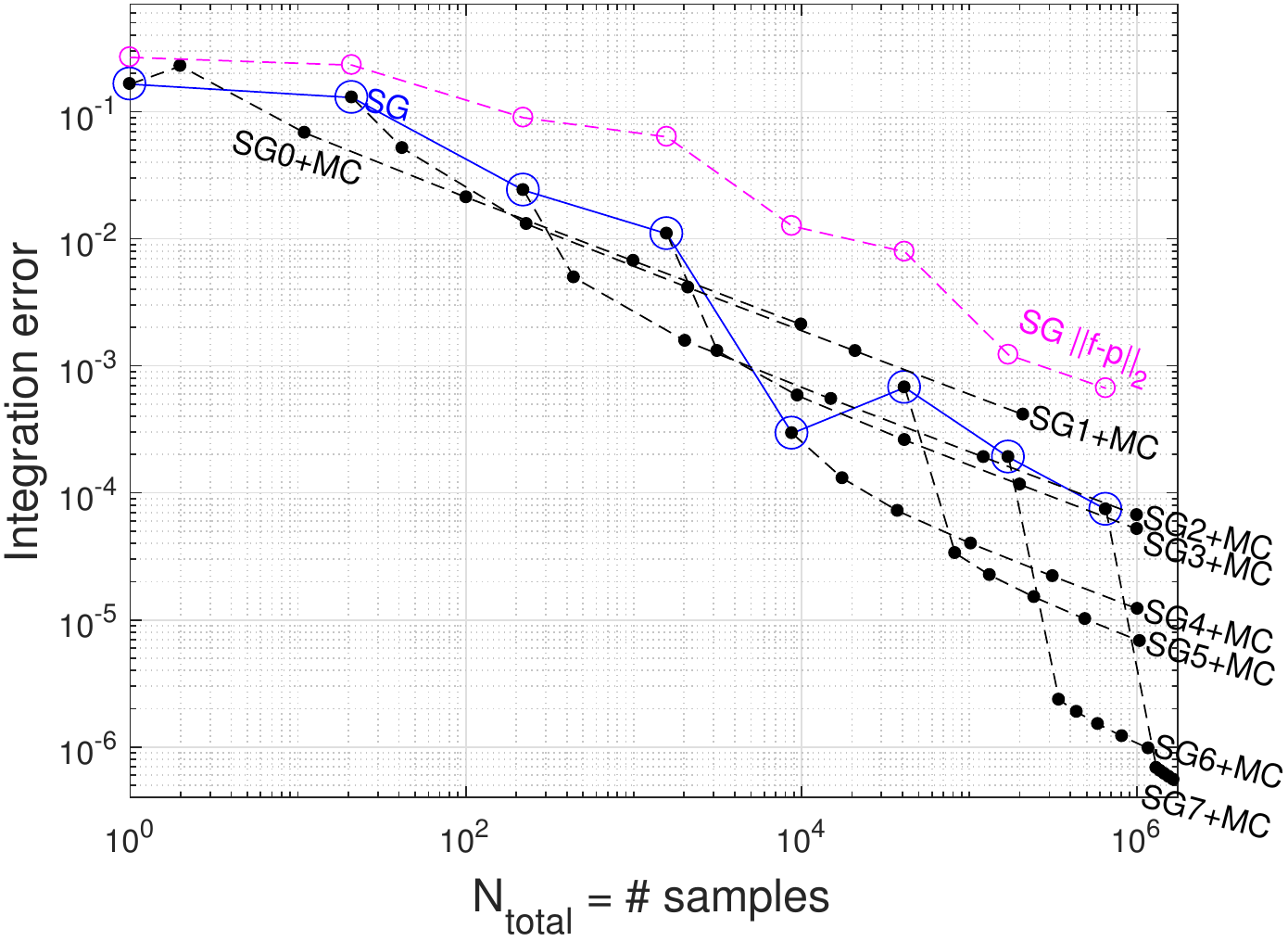}
  \end{minipage}   
  \begin{minipage}[t]{0.5\hsize}
      \includegraphics[width=0.95\textwidth]{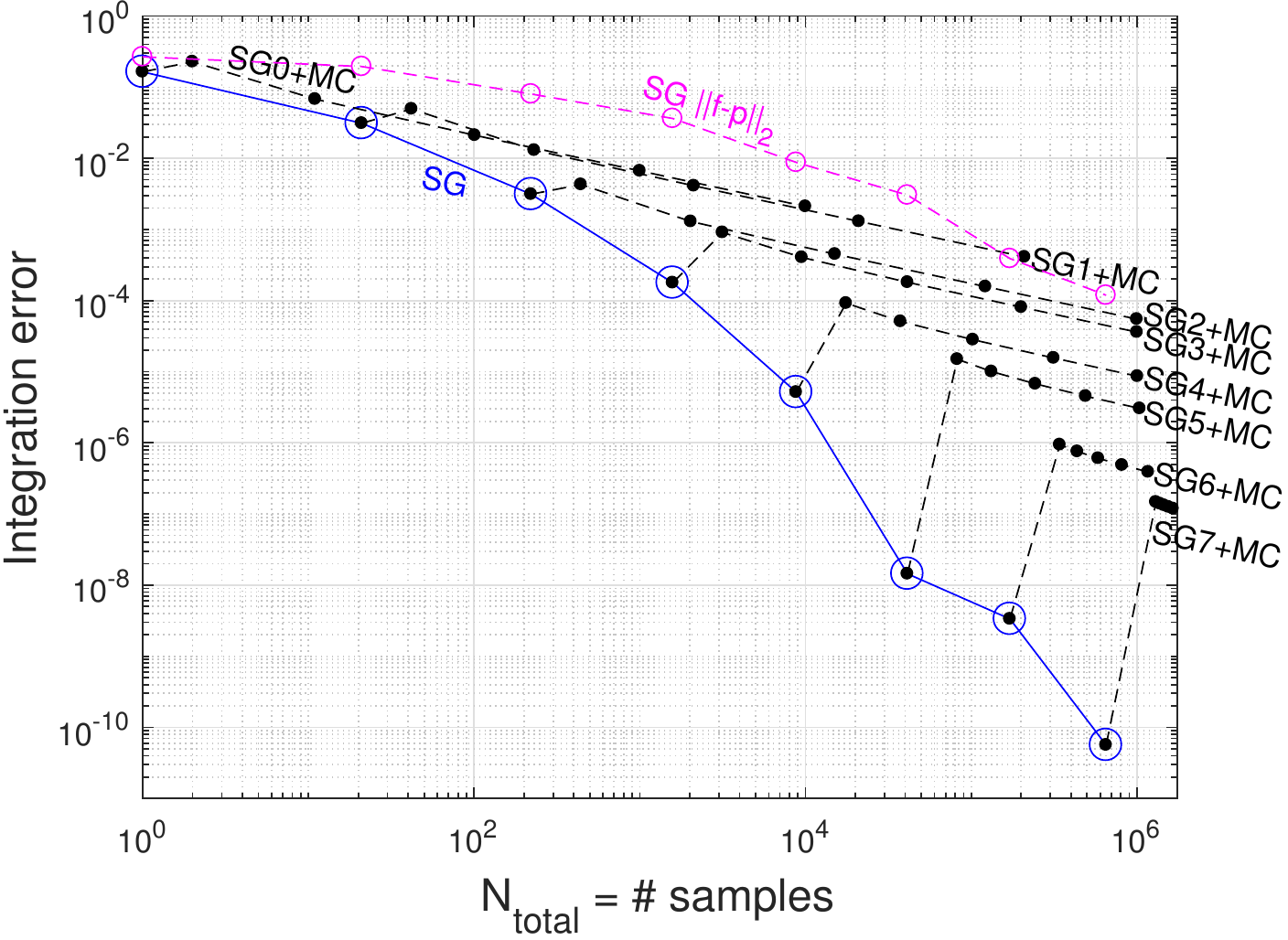}
  \end{minipage}
  \caption{
Repeating Figure~\ref{fig:spgridplusMC}, with $f =  \sin(\sum_{i=1}^d x_i)$, $d=10$. 
}
  \label{fig:spgridplusMC2}
\end{figure}

We used the Sparse Grid Interpolation Toolbox~\cite{Klimke_spalg,spgridtoolboxklimke2007} 
to find sparse grid (SG) interpolants to $f$. 
For the underlying one-dimensional quadrature rule we used two types of integrators: 
piecewise linear approximants and the Clenshaw-Curtis rule~\cite{novak1996high}. 
The results are shown in  Figures~\ref{fig:spgridplusMC} and \ref{fig:spgridplusMC2}. 

The solid curves show the convergence of sparse grid (SG) integration (i.e., $I-\int_\Omega p_s(\x)d\x$), as the level $s$ increases (the xth SG point takes $s=$x$-1$). Emanating from each SG point is the MCLS approximant (shown as SGx+MC) using that SG approximant as $\phi_1$. 

Note in Figure~\ref{fig:spgridplusMC} the sudden drop in error between each SG point and the first SGx+MC point. 
This is among the highlights of this paper, and warrants further explanation. 
Recall that the MCLS error estimate is  $\frac{\|f-p_s\|_2}{\sqrt{N}}$. 
 With SG, the integration errors here are roughly $\|f-p_s\|_2$ the upper bound in~\eqref{eq:interr}. By investing $N$ samples in addition to $N_s$ for SG, the error of MCLS improves (suddenly) from $\|f-p_s\|_2$ to $\frac{\|f-p_s\|_2}{\sqrt{N}}$. For example, SG7+MC uses $N_s\approx 10^6$ points for SG, so by taking $N=10^6$ more samples to double the overall samples (which corresponds to only a slight move in the $x$-axis in the graph), we improve the accuracy by a factor $\approx \sqrt{N}=10^3$. 

In addition to improving the accuracy, another  benefit of combining SG with MCLS is that it comes with a confidence interval, which the sparse grid alone does not provide. 

Each SGx+MC is an MCLS method. We see from the figures that the best accuracy for a given  $N_{total} = N_s+N$ is obtained by adaptively choosing x, growing it with $N_{total}$. This is MCLSA, which appears to converge faster than $O(1/N)$ in Figure~\ref{fig:spgridplusMC}. 
The optimal choice of $N_s$ and $N$ is a nontrivial matter and left for future work. In this example, taking $N_s\approx N$ appears to be a reasonable choice. Conceptually, we are spending $N_s$ work to approximate the function, and $N$ work to integrate. 

It is worth noting, however, that for smooth integrands, sparse grid quadrature based on Clenshaw-Curtis 
 can give integration accuracy much better than the function approximation, just like in the 1-dimensional case as we mentioned in Section~\ref{sec:mcasquad}. 
In such cases, using Monte Carlo does not improve the integration accuracy. 
This is illustrated in Figure~\ref{fig:spgridplusMC2} (right), where the integrand is analytic and the SG quadrature gives much higher accuracy than $\|f-p\|_2$, in fact higher than when combined with MCLS, unless a huge number of Monte Carlo samples are taken. 
We also note that while piecewise linear approximation gave a better approximate integral in Figure~\ref{fig:spgridplusMC}, the situation is opposite in Figure~\ref{fig:spgridplusMC2}; this reflects a well-known fact that for smooth (indeed analytic) functions, global approximation by high-degree polynomial approximation outperforms piecewise low-degree approximation. 
In any event, in all cases MCLS provides a confidence interval, and 
the error of MCLS is accurately estimated by $\frac{1}{\sqrt{N}}\|f-p_s\|_2$. 


While we have focused on the basic versions of SG approximants, other variants have been proposed. 
Most notably, the degree-adaptive sparse grids~\cite{gerstner2003dimension} is often able to detect low-dimensional structure in $f$ (if present) to significantly improve the approximation quality $\|f-p_s\|$, especially in higher dimensions $d$. 
These SG variants can be combined with MCLS in the same way. 

Combining SG with Monte Carlo has been considered for example in~\cite{nobile2016adaptive,nobile2015multi,Teseithesis}, where SG is used to estimate the integral of a control variate, for solving stochastic PDEs. Here SG is used as a control variate itself, for the more classical problem of integration. We believe our function approximation viewpoint makes the convergence analysis more transparent. 


We note that all the algorithms here 
and in Section~\ref{sec:leg}
are linear integrators, that is, 
$\hat I_N(f_1+f_2)=\hat I_N(f_1)+\hat I_N(f_2)$. Thus for example if $f$ is a sum of an analytic function and noise, then since the analytic part converges faster than $O(1/\sqrt{N})$, we expect MCLSA
to have the asymptotic convergence $O(1/\sqrt{N})$ with 
constant being the noise level.

\section{MCLS with other high-dimensional function approximation methods}\label{sec:mainalg}

Approximating a high-dimensional function $f$ is a very active area of research.
In addition to the polynomial least-squares and sparse grid approximation covered above, notable directions in this area include the use of separable (low-rank) functions~\cite{beylkin2009multivariate}, 
polynomial approximation combined with judicious choice of 
basis~\cite{cohen2015approximation} or compressed sensing~\cite{adcock1703.06987}, 
the functional analogue~\cite{bigoni2016spectral} of the tensor train decompositions~\cite{oseledets2011tensor}, 
radial basis functions~\cite{wendland2004scattered}, 
and approximation by Ridge functions~\cite{petrushev1998approximation}. 
It is premature to make a call on which method is the best, and the method of choice would most likely be problem-dependent. Here we briefly describe some of these possibilities that can easily be combined with MCLS. 

\subsection{Low-rank approximation}\label{sec:lowrank}
$p:\Omega\rightarrow \mathbb{R}$ is said to be of rank $r$ if 
\begin{equation}  \label{eq:beylkineq}
p(\x) = \sum_{j=1}^r g_{1j}(x_1)g_{2j}(x_2)\cdots g_{dj}(x_d), 
\end{equation}
for univariate functions $g_{ij}$. 
Functions admitting such representation with small $r$ are called 
low-rank or separable. 
Note that such functions are straightforward to integrate using univariate quadrature rules. 
Thus it is of interest to approximate $f$ by low-rank functions in the context of integration. 

To obtain such approximant $p\approx f$, 
a common approach is alternating least-squares, in which one modifies a single coordinate $\{g_{ij}(x_i)\}_{j=1}^r$ at a time. 
The algorithm in~\cite{beylkin2009multivariate} is an example, and 
it is straightforward to incorporate into MCLS as 
it  allows the sample points to be arbitrary. 



Low-rank approximation can be extremely powerful, for example for most of Genz' functions~\cite{genz1984testing}, which are rank-1. 
However, unlike the methods presented in the previous sections, the operation is not linear in $f$, that is, the output of  $f=f_1+f_2$ is not always the sum of the outputs of $f_1$ and $f_2$. The effectiveness of low-rank approximation usually deteriorates as the rank of $f$ increases. 


We also mention a recent work \cite{chevreuil2015least} that presents algorithms for approximating a function via a low-rank and sparse function, 
employing ideas in compressed sensing. As always, once a good approximant is obtained, we can combine it with MCLS to obtain accurate approximate integrals.

\subsection{Sparse approximation }\label{sec:compressedsensing}
Even when the problem lies in a high-dimensional space, it is sometimes possible to represent the function with compact storage (e.g. polynomials with sparse coefficients), largely independent of the ambient dimension. Such functions are identified in a number of studies, see e.g. \cite{CHKIFA2015400,cohen2015approximation} and reference therein. 
These include solutions of certain types of PDEs, which can be provenly approximated by a tensor product of Legendre polynomials, with error decaying algebraically with the number of nonzero coefficients. 
This fact is exploited in \cite{adcock1703.06987} to devise algorithms based on $\ell_1$ minimization to find polynomial approximation to high-dimensional problems. Incorporating these into MCLS(A) would be an interesting topic for the future. 

\ignore{
One approach would be to use the $\ell_1$-minimization algorithm in \cite{adcock1703.06987} to approximate the function, then integrate (either that approximant, or by Monte Carlo applied to $f-\hat f$). 
Here we take a simpler approach, wherein we identify the approximate sparsity more inexpensively, using the so-called basic thresholding~\cite[\S~3.3]{foucart2013mathematical} algorithm in compressed sensing. 

One simple way to estimate the coefficient $c_{i_1,i_2,\ldots ,i_d}$ in the Legendre expansion 
\begin{equation}  \label{eq:legexp}
f =  \sum_{i_1,\ldots,i_d} c_{i_1,\ldots,i_d}\prod_{j=1}^dP_{i_j}(x_j), 
\end{equation}
is to use Monte Carlo: this gives 
$\hat c_{i_1,\ldots,i_d} = \frac{1}{N}f(x)\prod_{j=1}^dP_{i_j}(x_j)$, and the variance is 
\[
\|f(x)\prod_{j=1}^dP_{i_j}(x_j)-c_{i_1,\ldots,i_d}\|^2\leq 
\|f(x)\prod_{j=1}^dP_{i_j}(x_j)\|^2.
\]
This motivates us to bound the term $\|P_i(x)P_j(x)\|^2$ from above. 
To do this we can use Adams's result~\cite{adams1878expression} that 
finds the coefficients for $\tilde P_i(x)\tilde P_j(x)=\sum_{k}c_k\tilde P_k(x)$, as then 
$\|P_i(x)P_j(x)\|^2=\sum_kc_k^2\|\tilde P_k(x)\|^2_2/\|\tilde P_i(x)\|_2\|\tilde P_j(x)\|_2$. (NO wrong! they don't go to zero as absolute values are taken..)

The key here is that it grow very moderately ($\leq 3$ for $n\leq 10^4$), and usually only few terms $i_j$ are nonzero in our context. 
} 



\section{Relations to classical Monte Carlo techniques}\label{sec:discuss}
Thus far in this paper we have mainly developed MCLS from the approximation theory viewpoint. Here we discuss various aspects of MCLS in terms of its relation to classical methods and techniques in Monte Carlo integration. 
\subsection{MCLS with quasi-Monte Carlo}\label{sec:lsqmc}
Quasi-Monte Carlo (QMC)~\cite{dick2013high}, \cite[Ch.~5]{glasserman2013monte} 
is a widely used method to obtain improved MC-type estimates, 
with convergence typically improved to close to $O(1/N)$ for ``sufficiently nice'' functions. 
The idea is to choose the sample points to lie more evenly in $[0,1]^d$ to reduce the so-called discrepancy. QMC is usually tied to the hypercube more than MC. 
Given the success of QMC in a number of applications, it is of interest to compare QMC with MCLS, and moreover to combine the two: take QMC sample points in MCLS. We refer to this algorithm as QMCLS. 

\ignore{
\begin{figure}[htbp]
  \begin{minipage}[t]{0.325\hsize}
  \centering
      \includegraphics[width=1\textwidth]{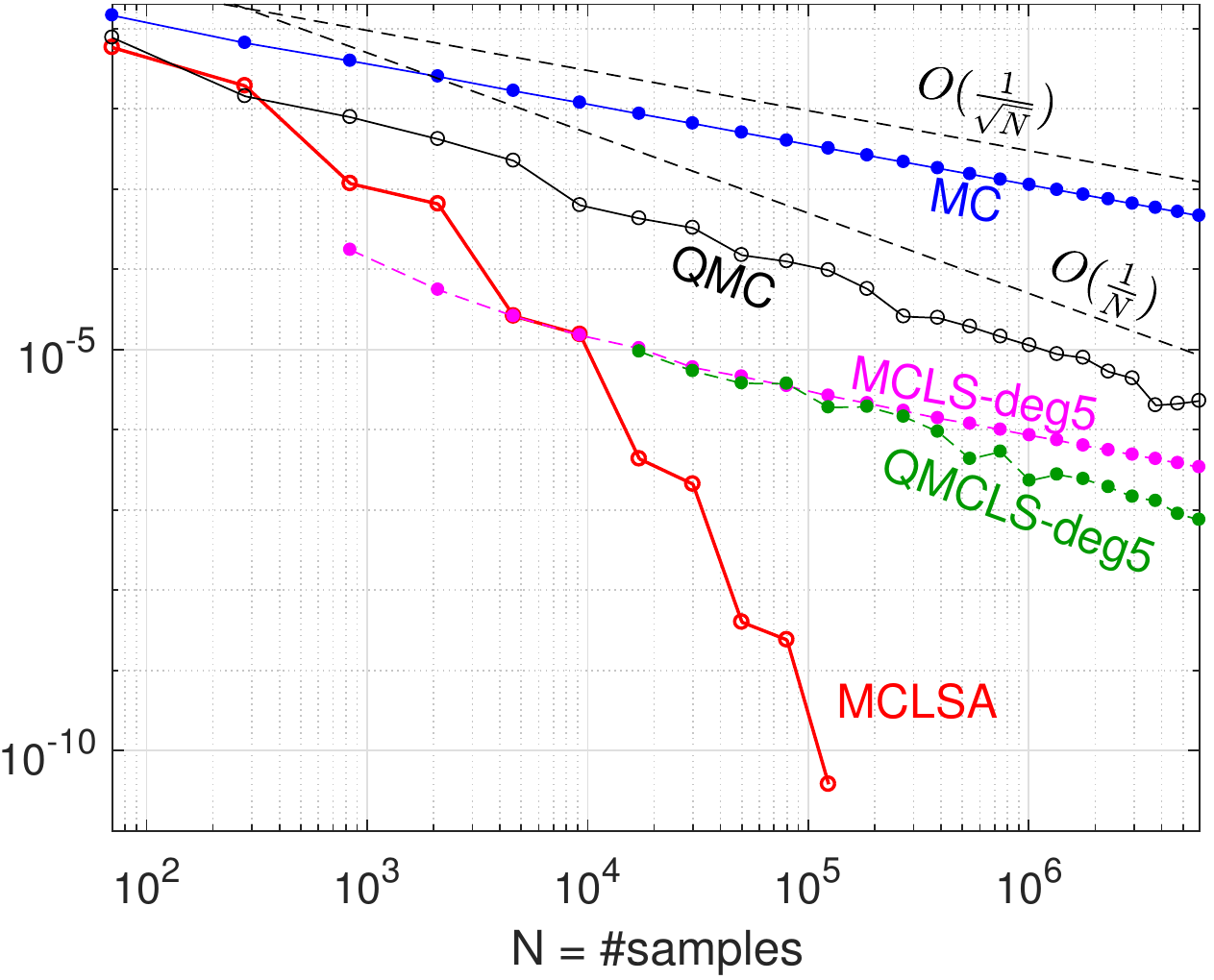}
  \end{minipage}
  \begin{minipage}[t]{0.33\hsize}
  \centering
      \includegraphics[width=.97\textwidth]{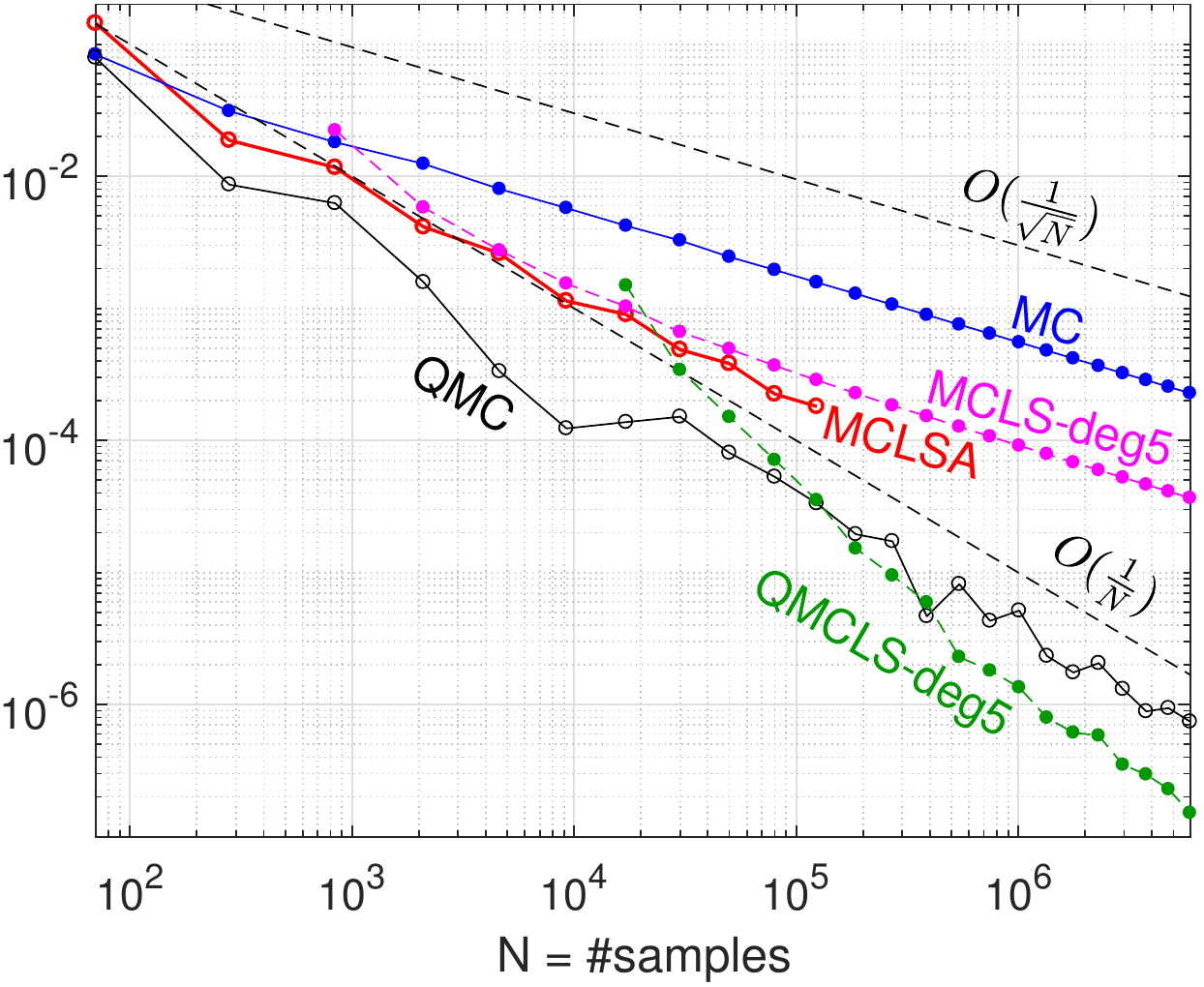}
  \end{minipage}
  \begin{minipage}[t]{0.325\hsize}
  \centering
      \includegraphics[width=.97\textwidth]{figs/norm1vs2CIRRRbasketd6Novern10k24donorm1fixd5}
  \end{minipage}
  \caption{
Same plots as Figure~\ref{fig:Legpolyinterp}, but with QMC versions included. 
Left: $f(\x) = \sin(\sum_{i=1}^d x_i)$, center: $f(\x) = \sum_{i=1}^d \exp(-|x_i-1/2|)$, right: basket option. 
}
  \label{fig:Legpolyinterp2old}
\end{figure}
}

\begin{figure}[htbp]
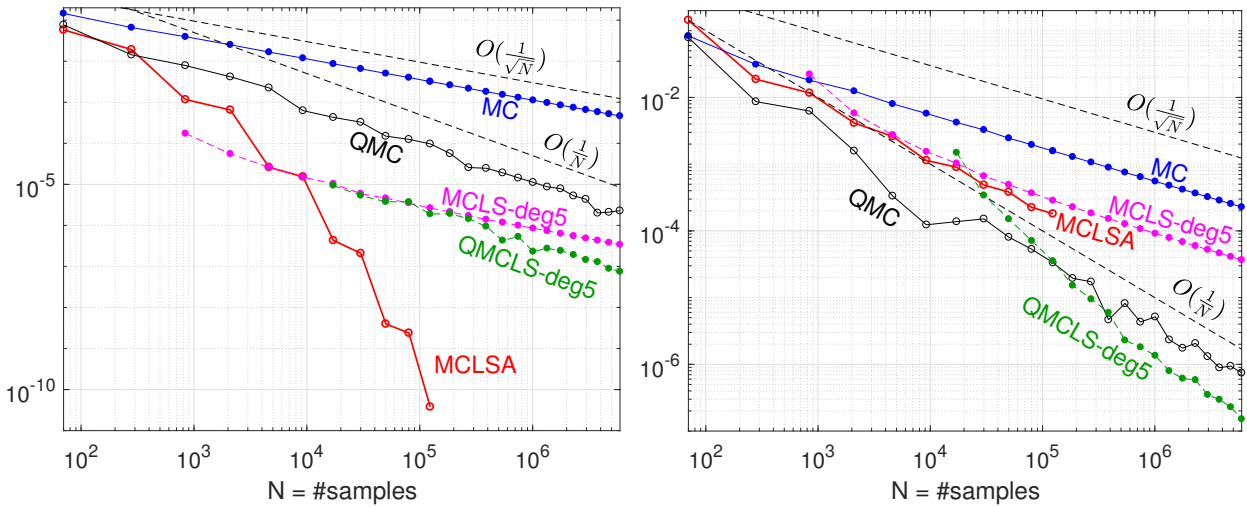

  \begin{minipage}[t]{0.5\hsize}
  \centering
      \includegraphics[width=1\textwidth]{figs/norm1vs2CIRRRsinsum6d6Novern10k24donorm1fixd5}
  \end{minipage}
  \begin{minipage}[t]{0.5\hsize}
  \centering
      \includegraphics[width=.98\textwidth]{figs/norm1vs2CIRRRgenzsumexpabsd6d6Novern10k24donorm1fixd5}
  \end{minipage}
  \caption{
Same plots as Figure~\ref{fig:Legpolyinterp} (left and center), but with QMC versions included. 
Left: $f(\x) = \sin(\sum_{i=1}^d x_i)$, Right: $f(\x) = \sum_{i=1}^d \exp(-|x_i-1/2|)$. 
}
  \label{fig:Legpolyinterp2}
\end{figure}

Figure~\ref{fig:Legpolyinterp2} shows the results of applying QMC and QMCLS to the same functions as in Figure~\ref{fig:Legpolyinterp} (we omit the basket option to make the figures large enough for visibility; 
the qualitative behavior is similar to the right plot). 
A few remarks are in order. 
First, for very smooth or analytic functions, MCLSA outperforms all other methods, reflecting the superalgebraic convergence. 
We note that an QMCLSA algorithm (integrating QMC with adaptively chosen degree) appears to be difficult as the optimal sampling described in Section~\ref{sec:leghigh} is not uniform in $[0,1]^d$ (a possibility would be to take QMC sample points and choose the degree s.t. $N=O(n^2)$). 
Second, QMCLS with fixed degree appears to have the same asymptotic convergence as QMC, close to $O(1/N)$.
The constant is smaller with QMCLS than QMC, but by how much depends on the function: the smoother, the wider the gap appears to become. 
Moreover, even for smooth functions (left plot), the improvement gained by QMCLS compared with QMC appears to be smaller than the difference between MCLS and MC.

\begin{figure}[htbp]
  \begin{minipage}[t]{0.5\hsize}
      \includegraphics[width=0.9\textwidth]{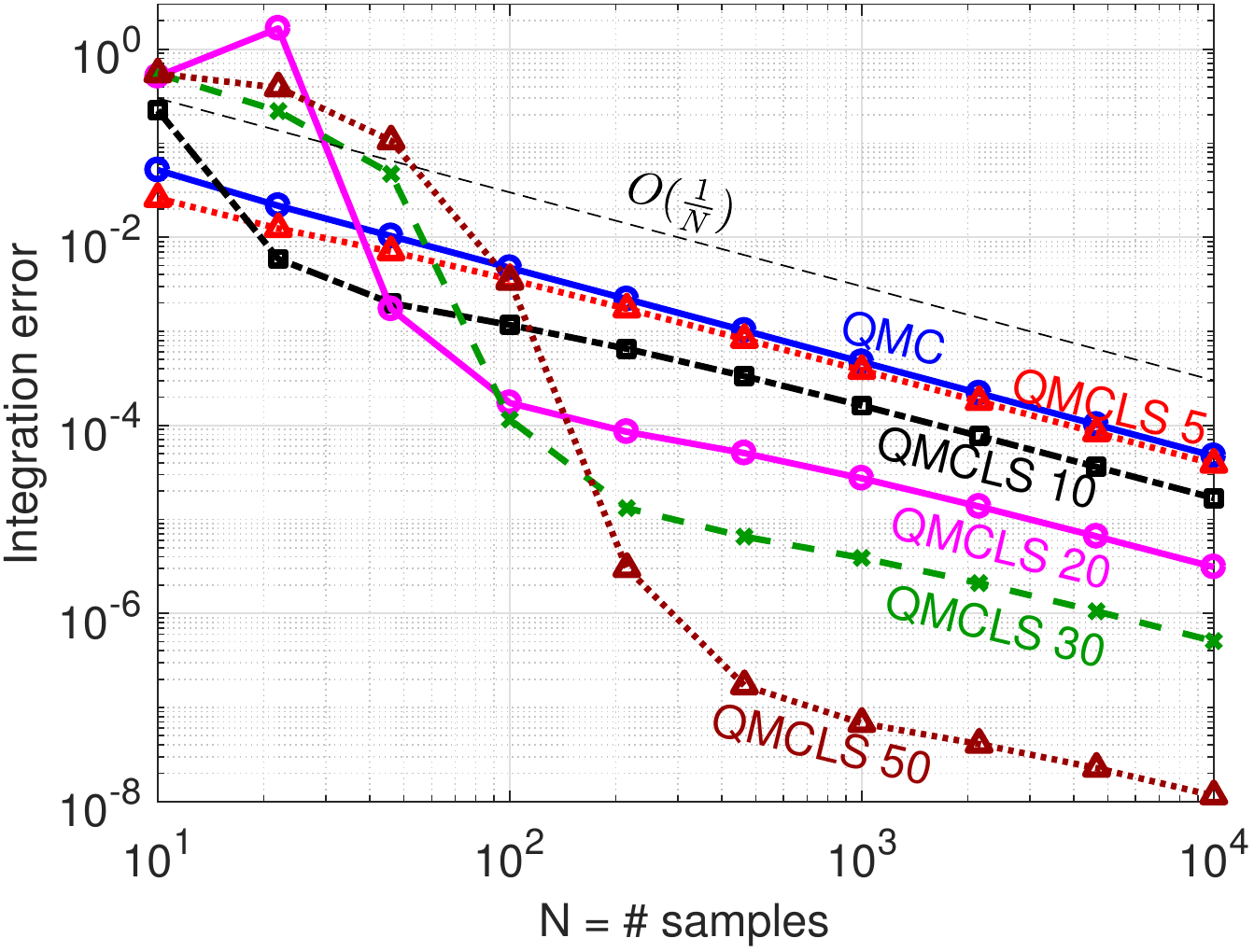}    
  \end{minipage}   
  \begin{minipage}[t]{0.5\hsize}
      \includegraphics[width=0.9\textwidth]{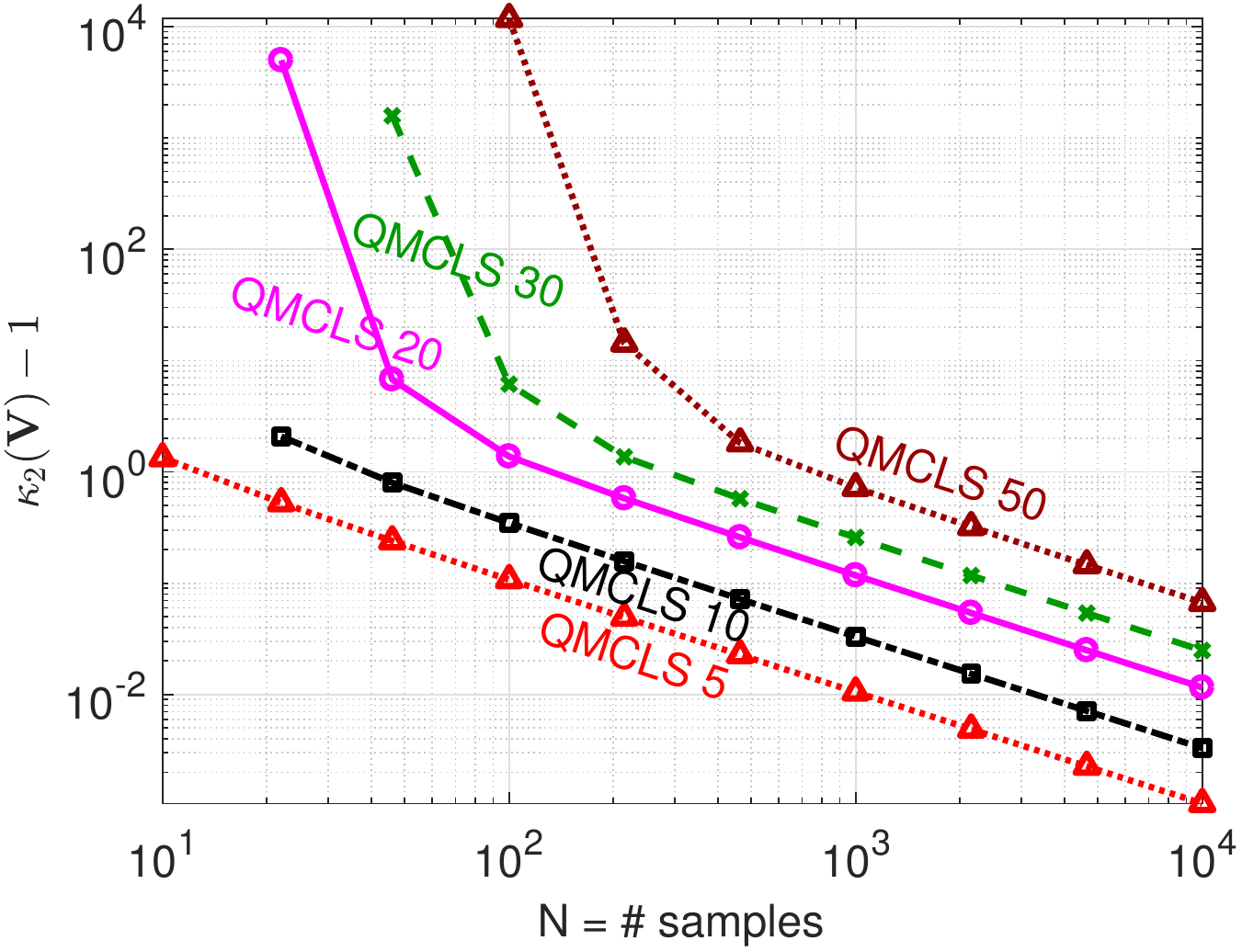}    
  \end{minipage}
  \caption{Same as Figure~\ref{fig:1d}, but using equispaced sample points on $[-1,1]$. 
Left: QMC vs. QMCLS errors for approximating $\int_{-1}^1\frac{1}{25x^2+1}$.  
Right: Conditioning $\kappa_2(\V)-1$. 
}
  \label{fig:1dqmc}
\end{figure}

To gain more insight, we repeat Figure~\ref{fig:1d}, a one-dimensional integration $\int_{-1}^1\frac{1}{25x^2+1}dx$, now using 
equispaced sample points on $[-1,1]$.  The result is shown in Figure~\ref{fig:1dqmc}. 
We see again that the gain provided by QMCLS is smaller if the polynomial degree is low. With sufficiently high degree, QMCLS becomes significantly better than QMC, with improvement factor similar to that of MCLS relative to MC. The asymptotic convergence of QMCLS appears to be the same as QMC, which here is $O(1/N)$. 
Making these observations precise is left for future work. 


\subsection{Relation to Monte Carlo variance reduction methods}
The function approximation viewpoint often gives us  fresh understanding of classical techniques for variance reduction in Monte Carlo. 
In Section~\ref{sec:leghigh} we mentioned how the nonuniform sampling is related to MC with importance sampling. 
In this subsection we describe other connections. 
Comprehensive treatments of variance reduction techniques in MC include 
Lemieux~\cite{LemieuxMC}, Owen~\cite{owenmcbook}, and Rubinstein and Kroese~\cite{rubinstein2016simulation}. 

\subsubsection{Multiple control variates} \label{sec:controlvariates}
We start with control variates, which has the strongest connection to MCLS. 
In the most basic form, one applies MC to $f-g$ instead of $f$, where 
$g:\Omega\rightarrow \mathbb{R}$ such that 
$\int_\Omega g(\x)d\x$ is known (and often assumed to be 0), and later added to the MC estimate to obtain 
$I\approx 
 \frac{1}{N}\sum_{i=1}^N(f(\x_i)-g(\x_i)) + \int_\Omega g(\x)d\x$. 
The analysis is usually presented in terms of the correlation between the integrand $f$ and the control variate $g$, in view of 
\[
\mbox{Var}(f-g) = \mbox{Var}(f)+\mbox{Var}(g)-2\mbox{Cov}(f,g). 
\]
This shows that a good control variate $g$ is one that ``correlates well'' with $f$. 

While the above description may bear little resemblance to MCLS, with a few more steps we can essentially obtain MCLS. First, we take multiple control variates $g_1,\ldots,g_n$. Second, we  use the estimator $I\approx \frac{1}{N}\sum_{i=1}^N(f(\x_i)-\sum_{j=1}^n c_j g(\x_i)) + \sum_{j=1}^n\int_\Omega c_jg_j(\x)d\x$ for some scalars $c_1,\ldots,c_n$. Third, we choose the scalars $c_j$ via regression, in which a least-squares problem~\eqref{eq:lsdisplay} is solved with $\phi_j=g_j$. We thus arrive at MCLS. We note that in MC with multiple control variates, $c_j$ are usually determined in a separate computation, as described at the end of Section~\ref{sec:bias} (which makes the estimator unbiased, but one might say it ``wastes'' the pilot samples). 

For the natural choice where $\phi_j$ is taken to be low-degree polynomials, 
this method is mentioned briefly in~\cite[Sec.~8.9]{owenmcbook}, but appears to not have been studied intensively. 

Despite the strong connection between MCLS and MC with control variates, 
the derivations are entirely different. 
Moreover, while ``approximating $f\approx p$'' and ``correlating $f,g$'' are closely related notions, conceptually we 
believe the approximation viewpoint is more transparent, which reveals the direct link between the MCLS variance and the approximant quality $\|f-p\|_2$, and leads to extensions. 
Indeed starting from MCLS with polynomial basis functions, we have introduced an MCLSA method that converges faster than $O(1/\sqrt{N})$ for smooth functions. Such extensions are unnatural and difficult from a control variate viewpoint.



\subsubsection{Stratification}\label{sec:stratify}

Stratification in MC is a technique to promote the samples to be more uniformly distributed, somewhat similar to QMC. 
Specifically, the domain is split into a certain number of subdomains, 
the number of random samples taken in each subdomain is specified. This leads to reduced variance~\cite[Sec.~8.4]{owenmcbook}. Using Lagrange multipliers, we can derive 
 the optimal budget allocation per stratum based on the variance of $f$ in each stratum. 

The function approximation viewpoint gives us a simple interpretation of stratification: the underlying function approximating $f$ is a \emph{piecewise constant} function. This is straightforward to see from the fact that each stratum is integrated as in standard MC. 

This observation leads immediately to an MCLS variant where we employ piecewise \emph{polynomial} approximants (which are not required to match on the boundaries). 
In numerical analysis, this process is called domain subdivision, and used effectively for example in rootfinding~\cite{boyd2002computing,nakatsukasa2013computing}. 
In the context of MC, this is 
a combination of stratification and control variates. 
Again, one can work out the optimal budget allocation per stratum. 
Importantly, here we allow the coefficients $\c$ to differ between strata. 
In MC, it appears to be more common to use a global $\c$ for all strata, which, from our viewpoint, clearly gives poorer approximation. 
Use of different coefficients is mentioned in the appendix of~\cite{hickernell2005control} as a comment by L'Ecuyer, but is apparently not widely used. 
We can even allow the basis functions (control variates) themselves to differ between strata. 

\begin{figure}[htbp]
  \begin{minipage}[t]{0.5\hsize}
  \centering
\includegraphics[width=.7\textwidth]{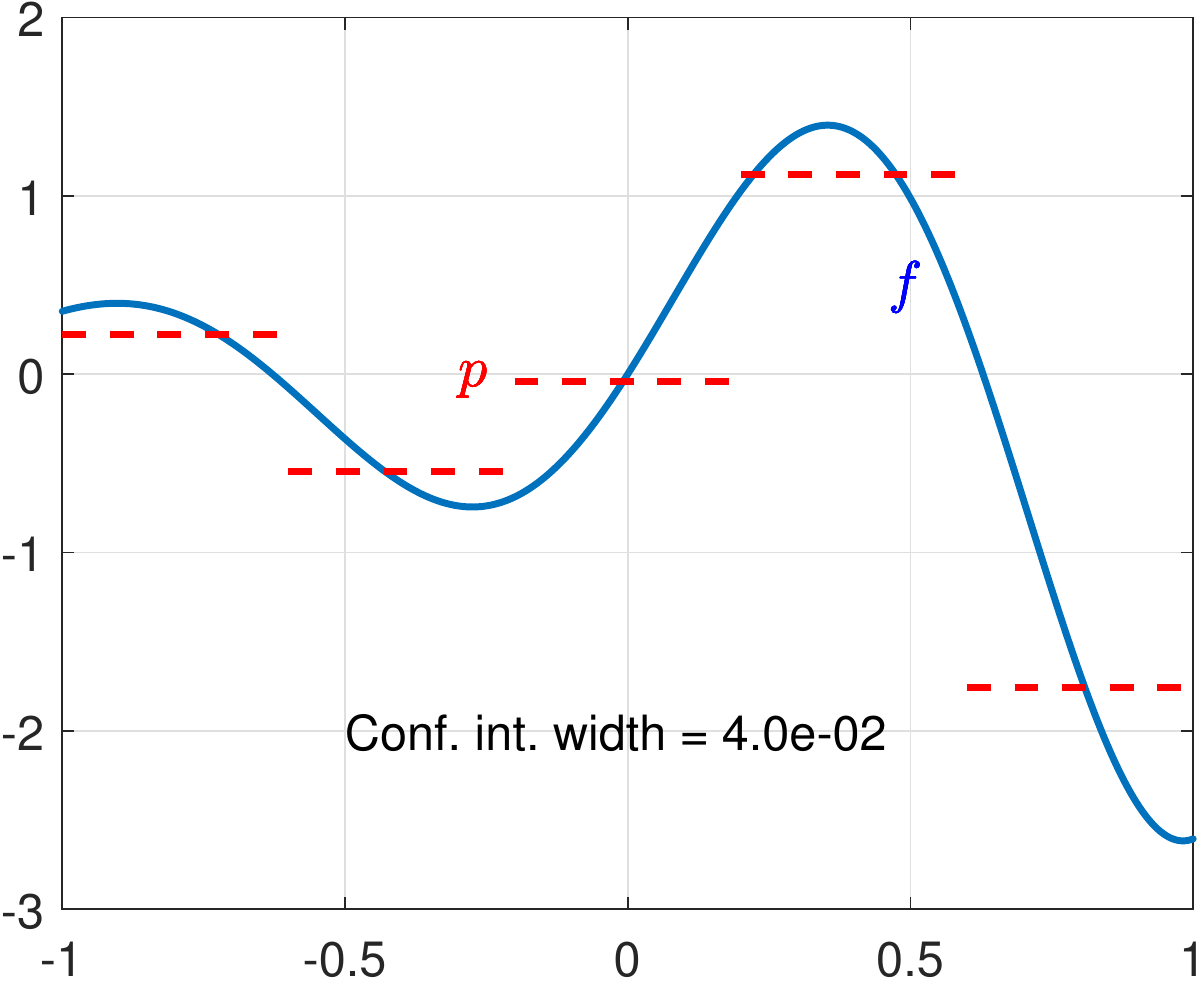}      

  \end{minipage}
  \begin{minipage}[t]{0.5\hsize}
  \centering
\includegraphics[width=.7\textwidth]{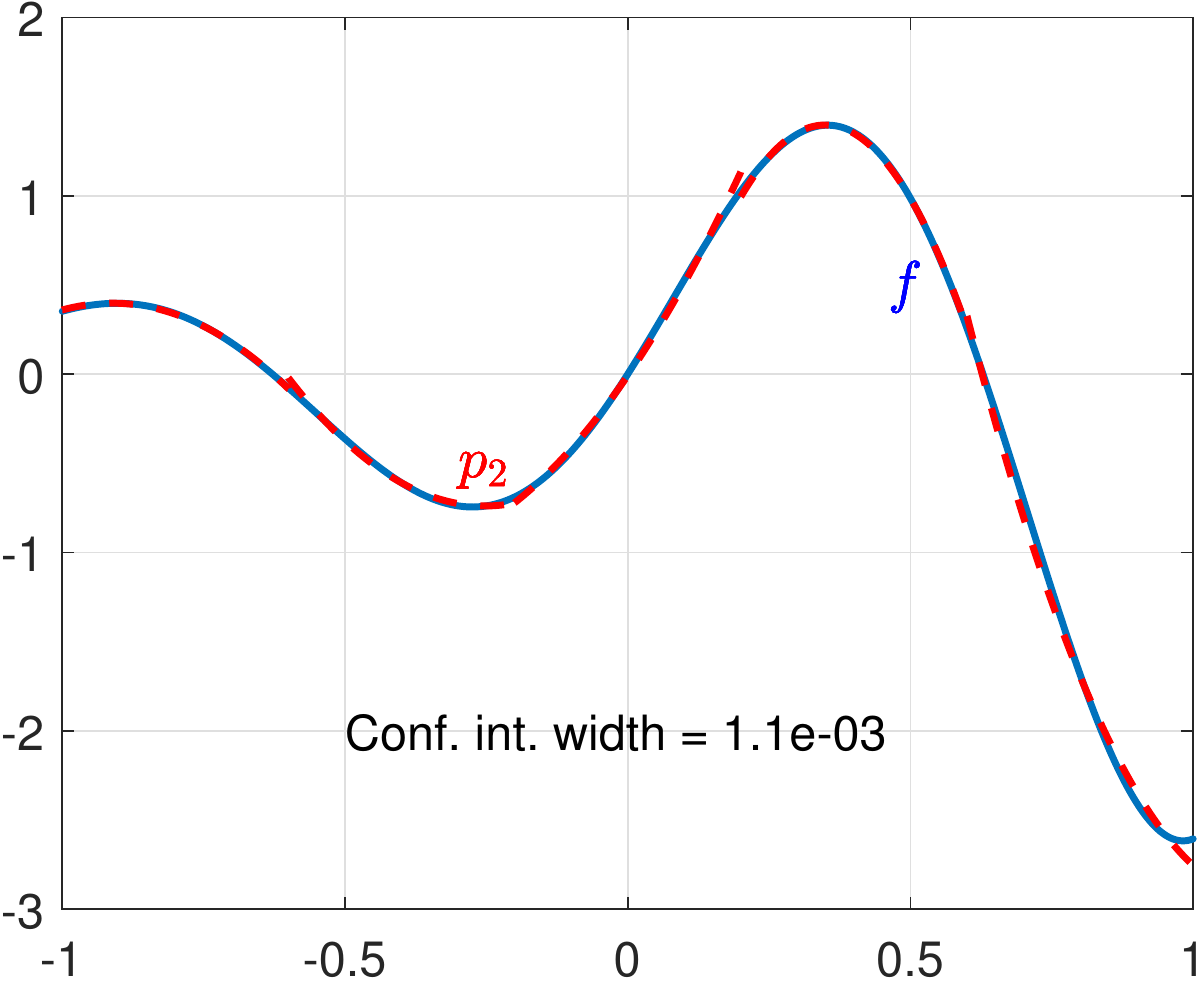}      
  \end{minipage}
  \caption{Underlying approximants $p$ (piecewise polynomials) for stratified MC and MCLS. Left: 
stratified MC uses piecewise constants. 
Right: stratified MCLS uses piecewise polynomials (here degree 2). 
Used $100$ total sample points. 
Compare with Figure~\ref{fig:plotkeyidea}, where MCLS used a degree 5 polynomial; if we do so here, the MCLS confidence interval width becomes $1.8\times 10^{-6}$. 
}
  \label{fig:strat}
\end{figure}

To illustrate the situation, Figure~\ref{fig:strat} shows the underlying approximation in MC with stratification and MCLS with piecewise polynomial approximants. It essentially repeats Figure~\ref{fig:plotkeyidea} employing stratification; observe how stratification improves the approximation quality $\|f-p\|_2$  for both MC and MCLS. Again, this directly means the variance is reduced. 

While stratification effectively reduces the variance, in high dimensions $d\gg 1$ one faces the obvious challenge that stratifying in each dimension leads to at least $2^d$ strata, making it impractical, as at least a few samples per stratum are needed for a confidence interval. Perhaps a good strategy would be to employ a judicious choice of stratification, where we split only in directions that matter. 


\subsubsection{Antithetics}\label{sec:anti}
Antithetics is a method where sample points are taken symmetrically. 
To simplify the argument, here we assume $\Omega=[-1,1]^d$. Then we take the MC sample points to be $\{\x_i\}_{i=1}^{N/2}\cup \{ -\x_i\}_{i=1}^{N/2}$, where $\{\x_i\}_{i=1}^{N/2}$ are chosen uniformly at random. Antithetics is known to reduce the variance signifcantly when $f$ is close to affine, or more generally, close to an odd function. 

Here is the MCLS viewpoint for antithetics: the symmetry in $\{\x_i\}_{i=1}^N$ implies some columns in~\eqref{eq:lsdisplay} are orthogonal to $[1,1,\ldots,1]^T$; namely the columns corresponding to $\phi_j$ that are odd functions $\phi_j(\x)=-\phi_j(-\x)$. 
In other words, the sample points are chosen in such a way that odd function integrate exactly to $0$; thus the odd part $\frac{1}{2}(f(\x)-f(-\x))$ of $f$ does not affect the outcome of $\hat I_N$. 

We can take advantage of this as follows: since with antithetics the odd basis functions play no role on $\hat I_N$, we can remove them from the least-squares problem~\eqref{eq:lsdisplay}. This clearly reduces $n$ and hence the MCLS cost. 



\section{Discussion and future directions}
The key observation of this paper is the interpretation of MC that connects approximation theory and Monte Carlo integration, which leads to a combination that often gets the best of both worlds.

One could argue that most---or indeed all---of the algorithms presented here could have been derived without the viewpoint of function approximation. This is not incorrect; for example, taking the control variates (CV) to be low-degree polynomials (which is mentioned in~\cite[Sec.~8.9]{owenmcbook}) essentially results in MCLS (with fixed degree) in Section~\ref{sec:leg}. 
However, we believe the function approximation viewpoint offers a number of advantages. 

Most importantly, the approximation viewpoint leads immediately to MCLSA, namely the idea of improving the approximant $p\approx f$ as more samples are taken. 
In this setting the basis function(s) (i.e., control variates) are chosen adaptively depending on $N$, and this is crucial for optimizing the accuracy, for example in Figure~\ref{fig:Legpolyinterp} (where $\phi_j$: polynomials) and~\ref{fig:spgridplusMC} ($\phi_1$: sparse grids). 
By contrast, in the context of MC with control variates, the variates are chosen a priori, so as to ``correlate well'' with $f$, and it is unnatural to add or change variates as $N$ increases. 
The approximation perspective also sheds light on other classical variance reduction techniques in MC, as we described in Section~\ref{sec:discuss}. 


Similarly, the use of the optimal weighting described in Section~\ref{sec:leghigh} is perfectly natural from the approximation viewpoint, but much less so from a control variate viewpoint. 


Last but not least, we believe it would always be conceptually helpful to think about the underlying function  approximation when using Monte Carlo methods, often offering deeper understanding. 


Significant and rapid developments are ongoing in the theory and practice of high-dimensional approximation, and it is fair to assume that further fundamental progress is imminent. 
These can have a direct impact on 
Monte Carlo integration using the ideas described in this paper.

While our main focus was integration on the hypercube, MCLS can be used for integration in other domains. Detailed investigation in such cases, including an effective choices of basis functions, is left for future work. 
Other directions include a more complete analysis of the relation between QMC and (Q)MCLS, applying MCLS to stochastic differential equations and combining MCLS with multi-level Monte Carlo methods~\cite{giles2015multilevel}.

\appendix

\section{Computational aspects}\label{sec:comp}
Once the $N$ samples are taken, the main computation in MCLS(A) is in solving the least-squares problem~\eqref{eq:ls}:
$\min_{\c\in\mathbb{R}^{n+1}}\|\V\c-\f\|_2$. 
\ignore{
The standard method to solve a least-squares problem~\eqref{eq:ls}, which we express as  
\begin{equation}\nonumber 
\min_{c\in\mathbb{R}^{d}}\|Vc-f\|_2,  
\end{equation}}
The standard method for solving least-squares problems employs the QR factorization~\cite[Sec.~5.3.3]{golubbook4th}: 
\begin{enumerate}
\item Compute the QR factorization $\V=\Q\R$.
\item Solve the square, upper-triangular linear system $\R\c=\Q^T\f$ for $\c$.  
\end{enumerate}

Mathematically this computes $\c = (\V^T\V)^{-1}\V^T\f$ (assuming $\V$ has full column rank), but the use of the QR factorization is recommended for numerical stability\footnote{The conditioning of least-squares problems is by no means straightforward; in particular it is much more complicated than linear systems~\cite[Ch.~20]{Higham:2002:ASNA}.}, over solving the mathematically equivalent normal equation $\V^T\V\c=\V^T\f$. 
With $N$ samples using $n$ basis functions, 
 the cost is $O(Nn^2)$, or more precisely $2Nn^2-\frac{2}{3}n^3$ flops~\cite[Sec.~5.3.3]{golubbook4th}. The dominant cost lies in the QR factorization, especially when $N\gg n$. The linear system requires just $n^2$ flops. 

Below we discuss how to efficiently update the solution $\c$ via the QR factorization when more samples are taken. 
For solving a single problem $\min_\c\|\V\c-\f\|_2$, faster algorithms requiring $O(Nn)$ operations are discussed in Appendix~\ref{sec:CG}; updating $\c$ in this context is also discussed there. 

\subsection{Updating solution as more samples are taken}\label{sec:update}
In Monte Carlo methods, it is often of interest to track the convergence by examining the values of the approximation~\eqref{eq:mcdef} at many values $N=N_1,N_2,\ldots$ with $N_1\leq N_2\leq \cdots$. For example, one may wish to sample enough so that the confidence interval becomes smaller than a prescribed width. 
Clearly in standard Monte Carlo methods this can be done, essentially at no redundant cost. 
Thus to obtain the Monte Carlo estimates for all values $N\in \{N_1,\ldots,N_\ell\}$, one needs just $O(N_\ell)$ flops. 
Here we show that MCLS (when $n$ is fixed) inherits such efficiency when updating the estimates. 
The essence of the process is described in~\cite[\S~4.3]{stewart1}. 

The essence of the linear algebra problem is as follows. Given a solution 
to $\min\|\V_1\c_1-\f_1\|_2$ with $\V_1\in\mathbb{R}^{N_1\times n}$ ($N_1\geq n$), 
obtain the solution with the extended sample size
\begin{equation}  \label{eq:lsextend}
\min_{\c\in\mathbb{R}^{n+1}}\left\|
  \begin{bmatrix}
\V_1\\\V_2    
  \end{bmatrix}\c
-
\begin{bmatrix}
\f_1\\\f_2  
\end{bmatrix}
\right\|_2. 
\end{equation}
The standard method for solving a least-squares problem 
is $\min\|\V_1\c_1-\f_1\|_2$
To solve this, the first task is to compute the QR factorization of 
$  \begin{bmatrix}
\V_1\\\V_2    
  \end{bmatrix}$. 
Let $\V_1=\Q_1\R_1$, $\V_2=\Q_2\R_2$ be QR factorizations (no assumption is made on the size of $\V_2$ it can be tall or fat). 
Then we have 
\[
  \begin{bmatrix}
\V_1\\\V_2    
  \end{bmatrix} =   
\begin{bmatrix}
\Q_1\R_1\\\Q_2\R_2
  \end{bmatrix}  = 
  \begin{bmatrix}
\Q_1\\ & \Q_2    
  \end{bmatrix}
\begin{bmatrix}
\R_1\\\R_2
  \end{bmatrix}. 
\]
Hence, with the QR factorization $\begin{bmatrix}
\R_1\\\R_2
  \end{bmatrix} = \Q_R\R$, we obtain 
\[
  \begin{bmatrix}
\V_1\\\V_2    
  \end{bmatrix} =   
\left(  \begin{bmatrix}
\Q_1\\ & \Q_2    
  \end{bmatrix}
\Q_R\right)
\R, 
\]
which is a proper QR factorization. We can then solve~\eqref{eq:lsextend} 
via the linear system $Rc=\left(  \begin{bmatrix}
\Q_1\\ & \Q_2    
  \end{bmatrix}
\Q_R\right)^T
\begin{bmatrix}
\f_1\\\f_2  
\end{bmatrix}=\Q_R^T\begin{bmatrix}\Q_1^T\f_1\\\Q_2^T\f_2  \end{bmatrix}$. 

Assuming $\min\|\V_1\c_1-\f_1\|_2$ has been solved, we have access to 
 the QR factorization $\V_1=\Q_1\R_1$ and the vector $\Q_1^T\f$. 
Then the extra cost lies in 
\begin{enumerate}
\item The $(N_2-N_1)\times d$ QR factorization $\V_2=\Q_2\R_2$, costing $O((N_2-N_1)d^2)$ flops, 
\item The $2d\times d$ QR factorization $\begin{bmatrix}
\R_1\\\R_2  \end{bmatrix} = \Q_R\R$, costing $O(n^3)$ flops, 
\item Computing the vector $b=\Q_R^T\begin{bmatrix}\Q_1^T\f_1\\\Q_2^T\f_2  \end{bmatrix}$, costing $O((N_2-N_1)n+n^2)$ flops as $\Q_1^T\f_1$ is known. 
\item Solving the linear system $\R\c=\b$. This costs $O(n^2)$ flops. 
\end{enumerate}
Finally, the integration $\int_\Omega \sum_{j=0}^nc_j\phi_j(\x)d\x$ is performed. This is $O(n)$ flops (we assume $\int_\Omega \phi_j(\x)d\x$ is available for each $j$). Overall, the cost for updating the estimate from $N_1$ samples to $N_2$ samples is $O((N_2-N_1)n^2+n^3)$. 
 The cost for computing the Monte Carlo estimates for all values $N\in \{N_1,\ldots,N_\ell\}$ is therefore 
$O(N_\ell+n^3\ell)$ flops, the same as in standard Monte Carlo, up to the term $n^3\ell$. 


\subsection{Updating solution as more basis functions are used}\label{sec:updatebasis}
An analogous updating scheme for the QR factorization is possible for appending basis functions for MCLSA, that is, when more columns are added to $\V$. 
Briefly: given the QR factorization $\V_1=\Q_1\R_1$, we compute the QR factorization of $[\V_1\ \V_2]$ via computing $\R_{12} = \Q_1^T\V_2$ and the QR factorization $(\I-\Q_1\Q_1^T)\V_2=\Q_2\R_{22}$. Then 
$[\V_1\ \V_2]=[\Q_1\ \Q_2]\big[
\begin{smallmatrix}
  \R_1& \R_{12}\\ & \R_{22}
\end{smallmatrix}
\big]$ is the QR factorization. 

\subsection{$O(Nn)$ algorithm: conjugate gradients}\label{sec:CG}
Contrary to the $O(Nn^2)$ complexity of the QR-based methods described above, 
the least-squares problem \eqref{eq:lsdisplay} can often be solved with $O(Nn)$ instead of $O(Nn^2)$ cost. The key fact 
is that $\V^T\V$ is converging to (a multiple of) identity. This means that for $N$ large enough, $\V^T\V$ is well conditioned with high probability (recall Section~\ref{sec:leghigh}). 
Therefore, the conjugate gradient (CG) method~\cite[Sec.~11.3]{golubbook4th} (alternatively MINRES, which minimizes the residual) applied to the normal equation $\V^T\V\c = \V^T\f$ converges geometrically at the rate $(\kappa_2(\V)-1)/(\kappa_2(\V)+1)$, converging to working precision usually in a couple of iterations. 

It appears to be more difficult to update a solution in a style explained above with CG, aside from using the previous solution as the initial guess. 


\ignore{
\subsection{paragraphs to maybe use later}
Our least-squares Monte Carlo method, including standard Monte Carlo, 
can thus be thought of as being based on approximating the function $f$ by a (very) simple function, based on massively (over)sampling $f$, followed by a least-squares fit. 
We wish to mention the roles of oversampling in numerical approximation theory. 
The famous Runge phenomenon-Newton Cotes, and resolution by oversampling+least-squares fitting. Used in~\cite{chkifa2015discrete,migliorati2013approximation} 
}
\subsection*{Acknowledgments }
I have benefited tremendously from discussions with Abdul-Lateef Haji-Ali. 
Nick Trefethen provided great insights into classical quadrature rules in numerical analysis: it was a discussion with him and Karlheinz Gr\"ochenig that inspired the key interpretation of MC as a quadrature rule. 
I thank Casper Beentjes, Mike Giles, Fabio Nobile, Christoph Reisinger, Ian Sloan, Bernd Sturmfels, and Alex Townsend for their perceptive comments. Most of this work was carried out in the Mathematical Institute at the University of Oxford, where the author was supported by JSPS as an Overseas Research Fellow. 

\def\noopsort#1{}\def\l{\char32l}\def\v#1{{\accent20 #1}}
  \let\^^_=\v\def\hbk{hardback}\def\pbk{paperback}


\begin{thebibliography}{10}

\bibitem{adcock1703.06987}
B.~Adcock, S.~Brugiapaglia, and C.~G. Webster.
\newblock Compressed sensing approaches for polynomial approximation of
  high-dimensional functions, 2017.

\bibitem{beylkin2009multivariate}
G.~Beylkin, J.~Garcke, and M.~J. Mohlenkamp.
\newblock Multivariate regression and machine learning with sums of separable
  functions.
\newblock {\em SIAM J. Sci. Comp}, 31(3):1840--1857, 2009.

\bibitem{bigoni2016spectral}
D.~Bigoni, A.~P. Engsig-Karup, and Y.~M. Marzouk.
\newblock Spectral tensor-train decomposition.
\newblock {\em SIAM J. Sci. Comp}, 38(4):A2405--A2439, 2016.

\bibitem{boyd2002computing}
J.~P. Boyd.
\newblock Computing zeros on a real interval through {C}hebyshev expansion and
  polynomial rootfinding.
\newblock {\em SIAM J. Numer. Anal.}, 40(5):1666--1682, 2002.

\bibitem{bungartz2004sparse}
H.-J. Bungartz and M.~Griebel.
\newblock Sparse grids.
\newblock {\em Acta Numerica}, 13:147--269, 2004.

\bibitem{caflisch1998monte}
R.~E. Caflisch.
\newblock Monte {C}arlo and quasi-{M}onte {C}arlo methods.
\newblock {\em Acta Numerica}, 7:1--49, 1998.

\bibitem{chevreuil2015least}
M.~Chevreuil, R.~Lebrun, A.~Nouy, and P.~Rai.
\newblock A least-squares method for sparse low rank approximation of
  multivariate functions.
\newblock {\em SIAM/ASA J. Uncertain. Quantif.}, 3(1):897--921, 2015.

\bibitem{CHKIFA2015400}
A.~Chkifa, A.~Cohen, and C.~Schwab.
\newblock Breaking the curse of dimensionality in sparse polynomial
  approximation of parametric {PDE}s.
\newblock {\em Journal de Math\'ematiques Pures et Appliqu\'ees}, 103(2):400 --
  428, 2015.

\bibitem{cohen2013stability}
A.~Cohen, M.~A. Davenport, and D.~Leviatan.
\newblock On the stability and accuracy of least squares approximations.
\newblock {\em Found. Comput. Math.}, 13(5):819--834, 2013.

\bibitem{cohen2015approximation}
A.~Cohen and R.~DeVore.
\newblock Approximation of high-dimensional parametric {PDEs}.
\newblock {\em Acta Numerica}, 24:1--159, 2015.

\bibitem{cohenoptimal}
A.~Cohen and G.~Migliorati.
\newblock Optimal weighted least-squares methods.
\newblock {\em SMAI-Journal of Computational Mathematics}, 3:181--203, 2017.

\bibitem{cools1997constructing}
R.~Cools.
\newblock Constructing cubature formulae: the science behind the art.
\newblock {\em Acta Numerica}, 6:1--54, 1997.

\bibitem{davis2007methods}
P.~J. Davis and P.~Rabinowitz.
\newblock {\em Methods of Numerical Integration}.
\newblock Courier Corporation, 2007.

\bibitem{dick2013high}
J.~Dick, F.~Y. Kuo, and I.~H. Sloan.
\newblock High-dimensional integration: the {quasi-Monte Carlo} way.
\newblock {\em Acta Numerica}, 22:133--288, 2013.

\bibitem{genz1984testing}
A.~Genz.
\newblock Testing multidimensional integration routines.
\newblock In {\em Proc. of International Conference on Tools, Methods and
  Languages for Scientific and Engineering Computation}, pages 81--94. Elsevier
  North-Holland, Inc., 1984.

\bibitem{gerstner2003dimension}
T.~Gerstner and M.~Griebel.
\newblock Dimension--adaptive tensor--product quadrature.
\newblock {\em Computing}, 71(1):65--87, 2003.

\bibitem{giles2015multilevel}
M.~B. Giles.
\newblock Multilevel {Monte Carlo} methods.
\newblock {\em Acta Numerica}, 24:259, 2015.

\bibitem{glasserman2013monte}
P.~Glasserman.
\newblock {\em Monte Carlo Methods in Financial Engineering}, volume~53.
\newblock Springer Science \& Business Media, 2013.

\bibitem{glynn2002some}
P.~W. Glynn and R.~Szechtman.
\newblock Some new perspectives on the method of control variates.
\newblock In {\em Monte Carlo and Quasi-Monte Carlo Methods 2000}, pages
  27--49. Springer, 2002.

\bibitem{golubbook4th}
G.~H. Golub and C.~F. Van~Loan.
\newblock {\em Matrix {C}omputations}.
\newblock {The Johns Hopkins University Press}, 4th edition, 2012.

\bibitem{haber1968combination}
S.~Haber.
\newblock A combination of {Monte Carlo} and classical methods for evaluating
  multiple integrals.
\newblock {\em Bull. Amer. Math. Soc.}, 74(4):683--686, 1968.

\bibitem{haber1969stochastic}
S.~Haber.
\newblock Stochastic quadrature formulas.
\newblock {\em Math. Comp.}, 23(108):751--764, 1969.

\bibitem{2017arXiv170700026H}
A.-L. {Haji-Ali}, F.~{Nobile}, R.~{Tempone}, and S.~{Wolfers}.
\newblock {Multilevel weighted least squares polynomial approximation}.
\newblock {\em ArXiv e-prints}, June 2017.

\bibitem{hampton2015coherence}
J.~Hampton and A.~Doostan.
\newblock Coherence motivated sampling and convergence analysis of least
  squares polynomial {C}haos regression.
\newblock {\em Comput. Methods in Appl. Mech. Eng.}, 290:73--97, 2015.

\bibitem{hickernell2005control}
F.~J. Hickernell, C.~Lemieux, and A.~B. Owen.
\newblock Control variates for quasi-{Monte Carlo}.
\newblock {\em Statistical Science}, 20(1):1--31, 2005.

\bibitem{Higham:2002:ASNA}
N.~J. Higham.
\newblock {\em Accuracy and Stability of Numerical Algorithms}.
\newblock SIAM, Philadelphia, PA, USA, second edition, 2002.

\bibitem{spgridtoolboxklimke2007}
A.~Klimke.
\newblock {S}parse {G}rid {I}nterpolation {T}oolbox -- user's guide.
\newblock Technical Report {IANS} report 2007/017, University of Stuttgart,
  2007.

\bibitem{Klimke_spalg}
A.~Klimke and B.~Wohlmuth.
\newblock Algorithm 847: {spinterp}: Piecewise multilinear hierarchical sparse
  grid interpolation in {MATLAB}.
\newblock {\em ACM Trans. Math. Soft.}, 31(4), 2005.

\bibitem{Kloedenbook}
P.~E. Kloeden and E.~Platen.
\newblock {\em Numerical Solution of Stochastic Differential Equations}.
\newblock Springer, 1992.

\bibitem{LemieuxMC}
C.~Lemieux.
\newblock {\em Monte Carlo and Quasi-Monte Carlo Sampling}.
\newblock Springer, 2009.

\bibitem{longstaff2001valuing}
F.~A. Longstaff and E.~S. Schwartz.
\newblock Valuing {A}merican options by simulation: a simple least-squares
  approach.
\newblock {\em Review of Financial studies}, 14(1):113--147, 2001.

\bibitem{manly2006randomization}
B.~F.~J. Manly.
\newblock {\em Randomization, bootstrap and {Monte Carlo} methods in biology},
  volume~70.
\newblock CRC press, 2006.

\bibitem{murphy2012machine}
K.~P. Murphy.
\newblock {\em Machine Learning: a Probabilistic Perspective}.
\newblock MIT press, 2012.

\bibitem{nakatsukasa2013computing}
Y.~Nakatsukasa, V.~Noferini, and A.~Townsend.
\newblock Computing the common zeros of two bivariate functions via
  {B}{\'e}zout resultants.
\newblock {\em Numer. Math.}, 129:181--209, 2015.

\bibitem{narayan2017christoffel}
A.~Narayan, J.~Jakeman, and T.~Zhou.
\newblock A {C}hristoffel function weighted least squares algorithm for
  collocation approximations.
\newblock {\em Math. Comp.}, 86(306):1913--1947, 2017.

\bibitem{nevai1986geza}
P.~Nevai.
\newblock G{\'e}za {Freud}, orthogonal polynomials and {C}hristoffel functions.
  {A} case study.
\newblock {\em J. Approx. Theory}, 48(1):3--167, 1986.

\bibitem{nobile2016adaptive}
F.~Nobile, L.~Tamellini, F.~Tesei, and R.~Tempone.
\newblock An adaptive sparse grid algorithm for elliptic {PDE}s with lognormal
  diffusion coefficient.
\newblock In {\em Sparse Grids and Applications-Stuttgart 2014}, pages
  191--220. Springer, 2016.

\bibitem{nobile2015multi}
F.~Nobile and F.~Tesei.
\newblock A {Multi Level Monte Carlo} method with control variate for elliptic
  pdes with log-normal coefficients.
\newblock {\em Stochastic Partial Differential Equations: Analysis and
  Computations}, 3(3):398--444, 2015.

\bibitem{novak1996high}
E.~Novak and K.~Ritter.
\newblock High dimensional integration of smooth functions over cubes.
\newblock {\em Numer. Math.}, 75(1):79--97, 1996.

\bibitem{oseledets2011tensor}
I.~V. Oseledets.
\newblock Tensor-train decomposition.
\newblock {\em SIAM J. Sci. Comp}, 33(5):2295--2317, 2011.

\bibitem{owenmcbook}
A.~B. Owen.
\newblock {\em Monte Carlo Theory, Methods and Examples}.
\newblock 2013.

\bibitem{peherstorfer2016optimal}
B.~Peherstorfer, K.~Willcox, and M.~Gunzburger.
\newblock Optimal model management for multifidelity monte carlo estimation.
\newblock {\em SIAM J. Sci. Comp}, 38(5):A3163--A3194, 2016.

\bibitem{petrushev1998approximation}
P.~P. Petrushev.
\newblock Approximation by ridge functions and neural networks.
\newblock {\em SIAM J. Math. Anal.}, 30(1):155--189, 1998.

\bibitem{robert2004monte}
C.~P. Robert and G.~Casella.
\newblock {\em Monte Carlo Methods}.
\newblock Wiley Online Library, 2004.

\bibitem{rubinstein2016simulation}
R.~Y. Rubinstein and D.~P. Kroese.
\newblock {\em Simulation and the Monte Carlo Method}.
\newblock John Wiley \& Sons, 2016.

\bibitem{smolyak1963quadrature}
S.~Smolyak.
\newblock Quadrature and interpolation formulas for tensor products of certain
  classes of functions.
\newblock In {\em Soviet Math. Dokl.}, volume~4, pages 240--243, 1963.

\bibitem{stewart1}
G.~W. Stewart.
\newblock {\em Matrix Algorithms Volume I: Basic decompositions}.
\newblock SIAM, Philadelphia, 1998.

\bibitem{suli2003introduction}
E.~S{\"u}li and D.~F. Mayers.
\newblock {\em An Introduction to Numerical Analysis}.
\newblock Cambridge University Press, 2003.

\bibitem{Teseithesis}
F.~Tesei.
\newblock {\em Numerical Approximation of Flows in Random Porous Media}.
\newblock PhD thesis, EPFL, 2016.

\bibitem{trefethenatap}
L.~N. Trefethen.
\newblock {\em Approximation Theory and Approximation Practice}.
\newblock SIAM, Philadelphia, 2013.

\bibitem{trefethen2017multivariate}
L.~N. Trefethen.
\newblock Multivariate polynomial approximation in the hypercube.
\newblock {\em Proc. Amer. Math. Soc.}, 145:4837--4844, 2017.

\bibitem{trefethen2014exponentially}
L.~N. Trefethen and J.~A.~C. Weideman.
\newblock The exponentially convergent trapezoidal rule.
\newblock {\em SIAM Rev.}, 56(3):385--458, 2014.

\bibitem{tropp2012user}
J.~A. Tropp.
\newblock User-friendly tail bounds for sums of random matrices.
\newblock {\em Found. Comput. Math.}, 12(4):389--434, 2012.

\bibitem{wendland2004scattered}
H.~Wendland.
\newblock {\em Scattered Data Approximation}.
\newblock Cambridge University Press, 2004.

\end{thebibliography}
\end{document}